\newif\ifcdc

\ifcdc
\documentclass[journal,twoside,web]{ieeecolor}
\usepackage{generic}
\else
\documentclass{article}
\usepackage[centering,margin=2cm]{geometry}
\fi
\usepackage{cite}

\usepackage{amsmath,amssymb,amsfonts,amsthm}
\usepackage{algorithmic}
\usepackage{graphicx}
\usepackage{algorithm,algorithmic}
\usepackage{textcomp}
\def\BibTeX{{\rm B\kern-.05em{\sc i\kern-.025em b}\kern-.08em
    T\kern-.1667em\lower.7ex\hbox{E}\kern-.125emX}}
\let\labelindent\relax
\usepackage{comment}
\usepackage{stmaryrd}
\usepackage{subcaption}
\usepackage{multirow}
\usepackage{mdframed}
\usepackage{cases,xspace,enumitem}
\usepackage{version}
\usepackage[prependcaption,colorinlistoftodos]{todonotes}

\makeatletter
\let\NAT@parse\undefined
\makeatother
\usepackage{hyperref}
\hypersetup{hidelinks=true}
\hypersetup{colorlinks=true, linkcolor=blue, breaklinks=true, urlcolor=blue, citecolor=blue}

\newtheorem{definition}{Definition}
\newtheorem{thm}{Theorem}[section]

\newtheorem{remark}{Remark}

\newtheorem*{assumption*}{Assumption}
\newenvironment{rem}{\begin{remark}}{\qed\end{remark}}

\newcommand{\bd}{\begin{definition}} 
\newcommand{\ed}{\end{definition}}
\newcommand{\bt}{\begin{thm}} 
\newcommand{\et}{\end{thm}}
\let\eps\varepsilon

\newcommand{\R}{\mathbb{R}}
\newcommand{\K}{\mathcal{K}}

\def\realnonnegative{\mathbb{R}_{\geq0}}
\def\realpositive{\mathbb{R}_{>0}}

\newcommand{\bi}{\begin{itemize}} 
\newcommand{\ei}{\end{itemize}} 
\newcommand{\bds}{\begin{description}} 
\newcommand{\eds}{\end{description}} 
\newcommand{\beq}{\begin{equation}} 
\newcommand{\eeq}{\end{equation}} 

\newcommand{\norm}[1]{\|#1\|}

\newcommand{\diag}[1]{[#1]}

\newcommand{\lognorm}[2]{\mu_{#1}(#2)}

\newcommand{\until}[1]{\{1,\dots, #1\}}

\newcommand{\subscr}[2]{#1_{\textup{#2}}}

\newcommand{\setdef}[2]{\{#1 \; | \; #2\}}

\newcommand{\map}[3]{#1\colon #2 \rightarrow #3}

\DeclareMathOperator*{\argmin}{arg\,min}

\newcommand{\ds}{\displaystyle}

\newcommand{\e}{\mathrm{e}}

\newcommand{\amin}{\subscr{a}{min}}
\newcommand{\amax}{\subscr{a}{max}}

\newcommand{\relu}{\operatorname{ReLU}}

\newcommand{\0}{\mbox{\fontencoding{U}\fontfamily{bbold}\selectfont0}}

\usepackage{booktabs} 
\usetikzlibrary{shapes,arrows,positioning,calc}
\usepackage[most]{tcolorbox}

\ifcdc
\includeversion{cdc}
\excludeversion{arxiv}
\else
\excludeversion{cdc}   
\includeversion{arxiv} 
\fi
\makeatletter
\ifcdc
  
\else
  \renewenvironment{arxiv}{\ignorespaces}{\unskip\ignorespacesafterend}
\fi
\makeatother

\definecolor{theoryblue}{RGB}{0, 102, 204}
\definecolor{assumptiongreen}{RGB}{0, 102, 51}
\newtcolorbox[auto counter, number within=section]{myassumption}[2][]{%
    colback=assumptiongreen!5,
    colframe=assumptiongreen,
    coltitle=black,              
    fonttitle=\bfseries,
    title={Assumption~\thetcbcounter~(#2)},
    enhanced,
    attach title to upper,
    after title={.\ },
    arc=5pt,
    outer arc=5pt,
    boxrule=0.4pt,
    left=4pt, right=4pt, top=2pt, bottom=2pt,
    drop shadow,
    #1                           
}

\definecolor{defred}{RGB}{188, 108, 107} 
\newtcolorbox[auto counter, number within=section]{mydefinition}[2][]{%
    colback=defred!5,          
    colframe=defred,           
    coltitle=black,             
    fonttitle=\bfseries,
    title={Definition~\thetcbcounter~(#2)},
    enhanced,
    attach title to upper,
    after title={.\ },
    arc=5pt,
    outer arc=5pt,
    boxrule=0.4pt,
    left=4pt, right=4pt, top=2pt, bottom=2pt,
    drop shadow,
    #1                          
}

\definecolor{thmblue}{RGB}{70, 102, 136} 
\newtcolorbox[auto counter, number within=section]{mytheorem}[2][]{%
    colback=thmblue!5,            
    colframe=thmblue,             
    coltitle=black,               
    fonttitle=\bfseries,
    title={Theorem~\thetcbcounter~(#2)},
    enhanced,
    attach title to upper,
    after title={.\ },
    arc=5pt,
    outer arc=5pt,
    boxrule=0.4pt,
    left=4pt, right=4pt, top=2pt, bottom=2pt,
    drop shadow,
    #1                            
}

\definecolor{lemmalilac}{RGB}{112, 92, 150} 

\newtcolorbox[auto counter, number within=section]{mylemma}[2][]{%
    colback=lemmalilac!5,
    colframe=lemmalilac,
    coltitle=black,
    fonttitle=\bfseries,
    title={Lemma~\thetcbcounter~(#2)},
    enhanced,
    attach title to upper,
    after title={.\ },
    arc=5pt,
    outer arc=5pt,
    boxrule=0.4pt,
    left=4pt, right=4pt, top=2pt, bottom=2pt,
    drop shadow,
    #1
}

\definecolor{corollarygold}{RGB}{184, 134, 11} 

\newtcolorbox[auto counter, number within=section]{mycorollary}[2][]{%
    colback=corollarygold!5,
    colframe=corollarygold,
    coltitle=black,
    fonttitle=\bfseries,
    title={Corollary~\thetcbcounter~(#2)},
    enhanced,
    attach title to upper,
    after title={.\ },
    arc=5pt,
    outer arc=5pt,
    boxrule=0.4pt,
    left=4pt, right=4pt, top=2pt, bottom=2pt,
    drop shadow,
    #1
}

\begin{document}
\title{A Unified Control-Theoretic Framework for Saddle-Point Dynamics in Constrained Optimization}
\title{CUPIDO: Constrained Unified PID \\
for Optimization}
\begin{arxiv}
\title{SPPID: Saddle-Point PID for Constrained Optimization}
\author{%
Veronica Centorrino\textsuperscript{1},
Rawan Hoteit\textsuperscript{1},
Efe C. Balta\textsuperscript{1,2},
and John Lygeros\textsuperscript{1}%
\thanks{%
This work was supported as a part of NCCR Automation, a National Centre of Competence in Research, funded by the Swiss National Science Foundation (grant number 51NF40\_225155).\\
\textsuperscript{1}Automatic Control Laboratory (IfA), ETH Zürich, 8092 Zürich, Switzerland,
{\tt\small\{vcentorrino, rhoteit, lygeros\}@ethz.ch}.\\%
\textsuperscript{2}Control and Automation Group, inspire AG,
8005 Zürich, Switzerland,
{\tt\small efe.balta@inspire.ch}.\\%
}
}
\end{arxiv}

\maketitle

\begin{abstract}
This paper studies constrained optimization problems through the lens of feedback control.
Building on the interpretation of Lagrange multipliers as feedback controllers, we propose the \emph{saddle-point PID (SPPID) dynamics}: a unified proportional--integral--derivative (PID) framework for continuous-time saddle-point dynamics.
The proposed dynamics employ PID control for equality constraints and anti-windup PI control for inequality constraints.
We show that SPPID is equivalent to a preconditioned primal--dual gradient flow of the augmented Lagrangian. This equivalence reveals the distinct role of each feedback component: integral action enforces constraint satisfaction, proportional action induces the augmented Lagrangian structure, and derivative action modifies the geometry of the primal dynamics via a state-dependent Riemannian metric.
For convex problems, we establish convergence to the KKT set together with conditions for global asymptotic stability.
For equality-constrained problems, we establish local exponential convergence for nonlinear equality constraints under standard assumptions, and show that projected gradient flow emerges as the infinite derivative-gain limit.
For strongly convex problems with affine equality or inequality constraints, we establish global exponential convergence by leveraging contraction theory.
Finally, we provide various numerical examples to illustrate the utility of SPPID.
\end{abstract}
\section{Introduction}
\begin{arxiv}
We study constrained optimization problems of the form
\end{arxiv}
\beq
\begin{aligned}
\min_{x \in \R^n} \quad & f(x)\\
\text{s.t.} \quad & h(x) = \0_p, \\
& g(x) \le \0_m,
\end{aligned}
\label{eq:eq_ineq}
\eeq
where $\map{f}{\R^n}{\R}$, $\map{h}{\R^n}{\R^p}$, and $\map{g}{\R^n}{\R^m}$ are continuously differentiable.
Such nonlinear programs arise throughout engineering, science, and machine learning, where equality constraints typically encode conservation laws, system dynamics, or coupling requirements between systems, whereas inequality constraints represent safety or resource limitations.

A classical approach to solving~\eqref{eq:eq_ineq} is through continuous-time primal--dual dynamics performing gradient descent and ascent on the associated Lagrangian with respect to the primal and dual variables, respectively. 
From a control-theoretic perspective, these flows can be viewed as closed-loop dynamical systems, where the primal dynamics constitute the plant, the constraint violations are the measured outputs, and the Lagrange multipliers are the control inputs. In this light, the standard primal--dual dynamics arise as an integral feedback on the constraint residuals.
This perspective has been explored in an empirical, algorithm-design setting in~\cite{AS-JA-PA:20}, and has more recently been formalized for general equality-constrained problems in~\cite{VC-SMF-SP-DR:25}, shifting the emphasis from designing optimization algorithms to designing feedback controllers for constrained dynamical systems.
Despite the growing interest in studying the interplay between optimization and control, the role of feedback design in constrained optimization remains underexplored. Existing methods are typically derived from optimization principles rather than control-theoretic design, for example, via primal--dual flows, augmented penalties, or projected dynamics~\cite{KJA-LH-HU:58, YT-GQ-NL:20, AH-ZH-SB-GH-FD:24}. There is currently no systematic understanding of how different feedback policies influence the resulting optimization dynamics.
Once optimization is viewed as a feedback design problem, a natural question arises:
\begin{center}
\textit{How does the choice of feedback policy on the dual variables influence the optimization dynamics?}
\end{center}
We address this question by developing a unified feedback-control framework for constrained optimization based on proportional--integral--derivative (PID) control of the dual variables. We show that PID control generates a family of saddle-point dynamics associated with the augmented Lagrangian. Rather than producing a single optimization algorithm, the proposed framework reveals a common control-theoretic principle underlying a broad class of primal--dual flows.
Beyond this unification, our approach makes explicit how the individual P, I, and D actions shape the closed-loop dynamics of the optimization algorithm, providing a systematic pathway for recovering known algorithms and designing new ones.

\subsubsection*{Literature review}
Studying optimization problems via continuous-time dynamics is a classical problem dating back to 1958~\cite{KJA-LH-HU:58, RWB:91}. This perspective has gained renewed interest due to the interpretation of optimization algorithms through the lens of feedback control~\cite{LSL-JWSP-EM:21, GB-JC-JIP-EDA:22, FD-ZH-GB-SB-JL-MM:24, AH-ZH-SB-GH-FD:24}. The standard approach for solving linear equality-constrained convex optimization problems relies on primal--dual (augmented) dynamics~\cite{KJA-LH-HU:58, AC-BG-JC:17, GQ-NL:19, YT-GQ-NL:20}. Other approaches to solving nonlinear programs include projected dynamical systems~\cite{AN-DZ:12}, which project the gradient of the objective function onto the cone of feasible descent directions, and through nonlinear sign-based control of projected gradient flows ~\cite{CF-RW:20}.

In reinforcement learning,~\cite{AS-JA-PA:20} proposed a heuristic PID update rule for the Lagrange multiplier of a safety cost constraint, empirically demonstrating reduced oscillations and overshoot in constraint violations during training.
More broadly, the use of Lagrange multipliers as feedback controllers has recently been analyzed for solving optimization problems with equality~\cite{VC-SMF-SP-DR:25, RZ-AR-JS-NL:25} and inequality~\cite{VC-SMF-SP-DR:24b, RZ-AR-JS-NL:25} constraints separately, as well as for smooth nonlinear programs via control barrier functions~\cite{AA-JC:24}.
For problems involving only equality constraints, in the setting of full-rank affine constraints and strongly convex, $L$-smooth objectives, global exponential convergence has been established for both I-controlled~\cite{GQ-NL:19} and PI-controlled~\cite{VC-SMF-SP-DR:25} dynamics.
The recent work~\cite{JR-SLJ:26} shows that PI control is equivalent to discrete-time primal--dual dynamics on the augmented Lagrangian.
For problems with only full-rank affine inequality constraints and strongly convex cost, both global exponential convergence~\cite{GQ-NL:19} and contractivity~\cite{AD-VC-AG-GR-FB:23f} of augmented primal--dual gradient dynamics have been established.
For general convex constrained problems with affine constraints and strongly convex cost, semi-global exponential stability of augmented primal--dual gradient dynamics has been established under standard constraint qualification~\cite{YT-GQ-NL:20}.
More broadly, there has been a growing interest in leveraging contractive dynamics to solve optimization problems~\cite{HDN-TLV-KT-JJES:18, AD-VC-AG-GR-FB:23f, VC-AD-AG-GR-FB:24a}. This is motivated by the exponential convergence and robustness guarantees that such dynamics enjoy~\cite{FB:26-CTDS}. Recent work has established contractivity for bilinear saddle-point problems in discrete time via operator-theoretic tools~\cite{CD-MB-PDG-JL-FD:24}.

Although existing works address the control of Lagrange multipliers for equality- or inequality-constrained problems separately, a systematic design and analysis for the \emph{simultaneous} treatment, is still lacking. Moreover, it remains unclear how derivative action can be incorporated into these multiplier dynamics with provable convergence guarantees, and what advantages it offers over purely PI designs.
\subsubsection*{Contributions}
We propose a unified control-theoretic framework for saddle-point dynamics arising in nonlinear programs. Our approach builds on the closed-loop system interpretation of~\cite{VC-SMF-SP-DR:25} and extends it to handle both equality- and inequality-constraints. The main contributions are summarized as follows.
\begin{enumerate}[label=\textup{(\roman*)}]
    \item We propose the \emph{saddle-point PID (SPPID) dynamics}, a family of saddle-point flows solving~\eqref{eq:eq_ineq} via PID control on the Lagrange multipliers. Specifically, we interpret the primal gradient flow as a MIMO plant with two independently controlled input--output channels, the equality and inequality multipliers being the respective inputs, and the corresponding constraint residuals the outputs.
    \item We prove forward invariance of the nonnegative orthant under the proposed inequality controller, guaranteeing dual feasibility along all trajectories, and show that the equilibria of SPPID coincide with the KKT points of~\eqref{eq:eq_ineq}.
    \item We show that SPPID is a preconditioned primal--dual flow of the augmented Lagrangian $\subscr{L}{aug}$ associated with~\eqref{eq:eq_ineq}, with the primal dynamics corresponding to Riemannian gradient descent on $\subscr{L}{aug}$ and the dual dynamics to scaled Euclidean gradient ascent. This yields a novel geometric interpretation of derivative feedback, arising solely from the equality controller, that does not emerge from standard primal--dual formulations. For convex problems, we further establish convergence to the KKT set, global asymptotic stability under strict convexity and linear independence constraint qualification (LICQ), and exponential convergence under additional regularity assumptions.
    \item We specialize the framework to equality-only and inequality-only constrained problems, recovering several known methods as limiting or special cases, and providing new results in each setting. For equality constraints, we show that the projected gradient flow arises as the infinite-derivative-gain limit of SPPID and provide, to the best of our knowledge, the first formal derivation of the projected gradient flow as a consequence of derivative feedback in a PID framework. We also establish global exponential convergence for affine constraints and local exponential stability for non-affine constraints under standard regularity assumptions, with explicit gain and rate bounds. For affine inequality constraints, we show that SPPID recovers known primal--dual flows, and we provide an explicit contraction metric and rate that complement the implicit characterizations in~\cite{AD-VC-AG-GR-FB:23f}.
\end{enumerate}
Finally, we validate the proposed framework via three numerical examples: non-convex sequential quadratic programming, multi-period Stackelberg games, and suboptimal nonlinear model predictive control application.

\subsubsection*{Outline}
Section~\ref{sec:preliminaries} presents the necessary mathematical preliminaries. Section~\ref{sec:sp_pid} introduces and develops the unified PID framework for solving~\eqref{eq:eq_ineq}. Sections~\ref{sec:equality_constraint} and~\ref{sec:inequality_constraint} specialize the framework to equality- and inequality-constrained problems, respectively, and establish convergence results in each setting. Section~\ref{sec:numerical_example} validates the proposed approach in simulation, and Section~\ref{sec:discussion_conclusion} concludes the paper.
\section{Mathematical Preliminaries}
\label{sec:preliminaries}
The sets of natural, real and positive real numbers are denoted by $\mathbb{N}$, $\R$ and $\R_{>0}$, respectively. For $n \in \mathbb{N}$, $\R^n$ denotes the $n$-dimensional Euclidean space. Vector inequalities of the form $x \leq (\geq)~y$ are entry-wise. We denote by $\0_n \in \R^n$ the zero vector in $\R^n$, and by $I_n$ the identity matrix in $\R^{n \times n}$. Given $x \in \R^n$, we let $\relu(x) := \max(x,\0_n)$, where the maximum is taken element-wise, and $\operatorname{diag}(x) \in \mathbb{R}^{n \times n}$ the diagonal matrix with $i$-th diagonal entry $x_i$. Given $A, B \in \R^{n\times n}$ symmetric, we write $A \preceq B$ (resp. $A \prec B$) if $B-A$ is positive semidefinite (resp. definite). For a symmetric matrix $A = A^\top$, we let $\subscr{\lambda}{min}(A)$ and $\subscr{\lambda}{max}(A)$ be its minimum and maximum eigenvalue, respectively.

\subsubsection*{Norms and Logarithmic Norms}
We let $\| \cdot \|$ denote both a norm on $\R^n$ and its corresponding induced matrix norm on $\R^{n \times n}$. Given $A \in \R^{n \times n}$ the \emph{logarithmic norm} (log-norm) induced by $\| \cdot \|$ is 
\begin{arxiv}
$\mu_{\| \cdot \|}(A) := \lim_{h\to 0^+} \dfrac{\norm{I_n + h A} -1}{h}.$~
\end{arxiv}%
Given $P \succ 0$, we let $\|\cdot\|_{P}$ be the $P$-weighted $\ell_2$-norm $\|x\|_{P} := \sqrt{x^\top P x}$, $x \in \R^n$. The corresponding log-norm is $\mu_{P}(A) =\min\setdef{b \in \R}{PA + A^\top P \preceq 2bP}$~\cite[Lemma~2.7]{FB:26-CTDS}.

\subsubsection*{Calculus and Function Classes} 
We denote by $\mathcal{C}^k$ the class of $k$-times continuously differentiable functions.
Let $\map{f}{\R^n}{\R}$ be twice differentiable. Then $\nabla f(x)$ and $\nabla^2 f(x)$ denote its gradient and Hessian matrices, respectively.
For $\map{f}{\R^{n+m}}{\R}$, we let $\nabla_{x} f(x,y) := \frac{\partial f}{\partial x}(x,y)$ be its partial gradients with respect to $x$. Given $\map{h}{\R^n}{\R^m}$ continuously differentiable, $J_h(x) \in \R^{m \times n}$ denotes its Jacobian matrix.
\begin{arxiv}
Finally, we recall the following standard definitions.
\begin{mydefinition}[label=def:strongly_convex_smooth]{Strongly convex and $L$-smooth maps}
A map $\map{f}{\R^n}{\R}$ is
\begin{enumerate}[label = (\roman*)]
\item \emph{$\rho$-strongly convex} if there exists $\rho > 0$ such that the map $x \mapsto f(x) - \frac{\rho}{2}\|x\|_2^2$ is convex;
\item \emph{$L$-smooth} if it is differentiable and there exists $L >0$ such that $\nabla f$ is $L$-Lipschitz.
\end{enumerate}
\end{mydefinition}
\smallskip
\end{arxiv}%

\subsubsection*{Optimality Conditions}
\begin{arxiv}
We review basic background on optimality conditions~\cite{DPB:97}.~
\end{arxiv}%
Consider a constrained optimization problem of the form~\eqref{eq:eq_ineq}. Let $\mathcal{F} := \setdef{x \in \R^n}{g(x) \leq \0_m,\ h(x) = \0_p}$ denote its feasible set and $\mathcal{A}(x) := \setdef{i \in \{1,\dots,m\}}{g_i(x) = 0}$ denote the set of active inequality constraint indices at $x \in \mathcal{F}$. 
Necessary conditions for optimality requires appropriate constraint qualification (CQ) conditions at the point of interest. We mainly consider LICQ at $x \in \mathcal{F}$, that is, the gradients $\{\nabla g_i(x)\}_{i \in \mathcal{A}(x)} \cup \{\nabla h_j(x)\}_{j = 1, \dots p}$ are linearly independent.
\begin{arxiv}
\begin{mydefinition}[label=def:stationary_point]{KKT conditions}
A point $(x^\star, \lambda^\star, \mu^\star) \in \R^{n+p+m}$ satisfies the Karush--Kuhn--Tucker (KKT) conditions if
\begin{subequations}
\label{eq:kkt_conditions}
\begin{align}
\nabla f(x^\star) + J_h(x^\star)^\top \lambda^\star + J_g(x^\star)^\top \mu^\star &= \0_n,
\label{kkt_stationarity}\\
h(x^\star) &= \0_p,
\label{kkt_stationarity_eq}\\
g(x^\star) &\le \0_m,
\label{kkt_primal_feasibility}\\
\mu^\star &\ge \0_m,
\label{kkt_dual_feasibility}\\
\mu^\star_i g_i(x^\star) &= 0,\quad \forall i \in \{1,\dots,m\}.
\label{kkt_complementarity}
\end{align}
\end{subequations}
\end{mydefinition}
\smallskip
\end{arxiv}%
The pair $(\lambda^\star, \mu^\star)$ is the \emph{Lagrange multiplier} associated with the equality and inequality constraints, respectively. The set of KKT points is $\mathcal{K} := \setdef{(x^\star, \lambda^\star, \mu^\star)}{(x^\star, \lambda^\star, \mu^\star) \text{ satisfies~\eqref{eq:kkt_conditions}}}$.
If $x^\star$ is a local minimizer of~\eqref{eq:eq_ineq} and an appropriate CQ holds at $x^\star$, then there exist  $\lambda^\star$ and $\mu^\star$ such that $(x^\star, \lambda^\star, \mu^\star)$ satisfies~\eqref{eq:kkt_conditions}. Moreover, if LICQ holds  at $x^\star$, the multipliers $(\lambda^\star, \mu^\star)$ are unique.
When $f$ and $g$ are convex and $h$ is affine, the KKT conditions are sufficient for global optimality. 
\begin{arxiv}
Therefore, when a CQ additionally holds at $x^\star$, the KKT conditions are both necessary and sufficient for global optimality.
\end{arxiv}%

Second-order conditions sharpen the KKT conditions by characterizing local optimality through the curvature of the Lagrangian. We say that the Second Order Sufficient Condition (SOSC) holds at $(x^\star, \lambda^\star, \mu^\star)$ if for all 
$d \in\setdef{d \in \R^n}{\nabla g_i(x^\star)^\top d \le 0, i\in \mathcal{A}(x^\star)\text{ with }\mu_i^\star = 0; \nabla g_i(x^\star)^\top d = 0, i\in\mathcal{A}(x^\star)\text{ with }\mu_i^\star>0, J_h(x^\star)d=0_p}\setminus\{0_n\} $ it holds
$$
d^\top \Big( \nabla^2 f(x^\star) + \sum_{i=1}^m \mu_i^\star \nabla^2 g_i(x^\star) + \sum_{j=1}^p \lambda_j^\star \nabla^2 h_j(x^\star) \Big) d > 0.
$$
If $(x^\star, \lambda^\star, \mu^\star)$ is a KKT point at which SOSC holds, then $x^\star$ is a strict local minimizer of~\eqref{eq:eq_ineq}.
\smallskip

\subsubsection*{Contraction Theory}
Consider a dynamical system 
\beq
\label{eq:dynamical_system}
\dot{x}(t) = f\bigl(t,x(t)\bigr),
\eeq 
where $\map{f}{\R_{\geq 0} \times \R^n}{\R^n}$. We give the following~\cite{FB:26-CTDS}.
\begin{arxiv}
\begin{mydefinition}[label=def:contracting_system]{Contracting dynamics}
Given a norm $\norm{\cdot}$ with associated log-norm $\mu_{\norm{\cdot}}$, a smooth function $\map{f}{\R_{\geq 0} \times \mathcal{X}}{\R^n}$, with $\mathcal{X} \subseteq \R^n$ open, convex, and forward invariant for~\eqref{eq:dynamical_system}, and a \emph{contraction rate} $c >0$ ($c = 0$), $f$ is $c$-strongly (weakly) infinitesimally contracting on $\mathcal{X}$ if
$$
\mu_{\norm{\cdot}}\bigl(J_f(t, x)\bigr) \leq -c,\quad \text{ for all } x \in \mathcal{X} \textup{ and } t\in \R_{\geq0},
$$
where $J_f(t,x) := \frac{\partial f}{\partial x}(t,x)$.
\end{mydefinition}
\smallskip
\end{arxiv}%
The above definition immediately implies the following.
\begin{arxiv}
\begin{mylemma}[label=lem:contractivity]{Global exponential stability}
Let $c>0$, and let $\map{f}{\R_{\geq 0} \times \mathcal{X}}{\R^n}$ be $c$-strongly infinitesimally contracting on a set $\mathcal{X} \subseteq \R^n$, open, convex, and forward invariant for~\eqref{eq:dynamical_system}. Then for any two trajectories $x(\cdot)$ and $y(\cdot)$ of~\eqref{eq:dynamical_system} with initial conditions $x_0, y_0 \in \mathcal{X}$, respectively,
$$
\|x(t) - y(t)\| \leq \e^{-ct}\|x_0 -y_0\|, \textup{ for all } t \geq 0,
$$
i.e., the distance between the two trajectories converges exponentially with rate $c$.
\end{mylemma}
\smallskip
\end{arxiv}%
One of the main benefits of strong contractivity is that it ensures global exponential convergence to the unique equilibrium, along with favorable robustness properties~\cite{FB:26-CTDS}. Moreover, Euler discretization with sufficiently small step sizes preserves convergence guarantees~\cite{FB-PCV-AD-SJ:21e}.
\section{Unified Framework for Nonlinear Programs}
\label{sec:sp_pid}
\label{sec:problem_formulation}
\begin{arxiv}
Consider the constrained problem in~\eqref{eq:eq_ineq}, which we rewrite here for convenience
\begin{align*}
\min_{x \in \R^n} \quad & f(x)\\
\text{s.t.} \quad & h(x) = \0_p, \\
& g(x) \le \0_m,
\end{align*}
where $\map{f}{\R^n}{\R}$, $\map{h}{\R^n}{\R^p}$ and $\map{g}{\R^n}{\R^m}$ are continuously differentiable functions. We assume that~\eqref{eq:eq_ineq} admits at least one optimal solution at which a CQ holds, so that $\mathcal{K} \neq \emptyset$.~
\end{arxiv}%
We cast the design of continuous-time dynamics solving~\eqref{eq:eq_ineq} as a feedback control problem, formulating a closed-loop system
comprising the primal gradient-flow dynamics (the plant) and a Lagrange multiplier feedback controller.

Consider the Lagrangian associated with~\eqref{eq:eq_ineq}, that is, the map $\map{L}{\R^n \times \R^p \times \R^m}{\R}$ defined by
\begin{arxiv}
$$L(x, \lambda, \mu)=f(x) + \lambda^{\top} h(x) + \mu^\top g(x),$$
\end{arxiv}%
where $\lambda \in \R^p$ and $\mu \in \R_{\geq 0}^m$ are the Lagrange multipliers.
Treating $\lambda$, $\mu$ as control inputs to the gradient flow of $L$ with respect to $x$ yields the following system with inputs $(\lambda,\mu)$ and outputs $(y_1,y_2)$
\begin{arxiv}
\beq
\label{eq:open_system_eq_ineq}
\begin{cases}
\dot{x}(t)= - \nabla_x L(x(t), \lambda(t), \mu(t)) = - \nabla f(x(t)) - J_h(x(t))^\top\lambda(t) - J_g(x(t))^\top \mu(t) \\
y_1(t) = h(x(t)),\\
y_2(t) = g(x(t)).
\end{cases}
\eeq
\end{arxiv}%
The control objective is twofold. For the equality constraints, we seek a feedback policy for $\lambda$ that drives $y_1$ to zero. For the inequality constraints, we seek a feedback policy for $\mu$ that regulates $y_2$ to enforce dual feasibility and complementary slackness. Moreover, both control objectives need to ensure convergence of~\eqref{eq:open_system_eq_ineq} to an equilibrium satisfying~\eqref{eq:kkt_conditions}. In the next section, we design feedback policies that achieve these objectives and analyze the resulting closed-loop dynamics.
\subsection{PID Control of Nonlinear Optimization Problems}
We show that PID control on the Lagrange multipliers gives rise to a family of saddle-point flows associated with the augmented Lagrangian. For simplicity, we subsequently drop the time dependence $(t)$.
\begin{arxiv}
\begin{figure}[!h]
    \centering
    \includegraphics[width=0.48\linewidth]{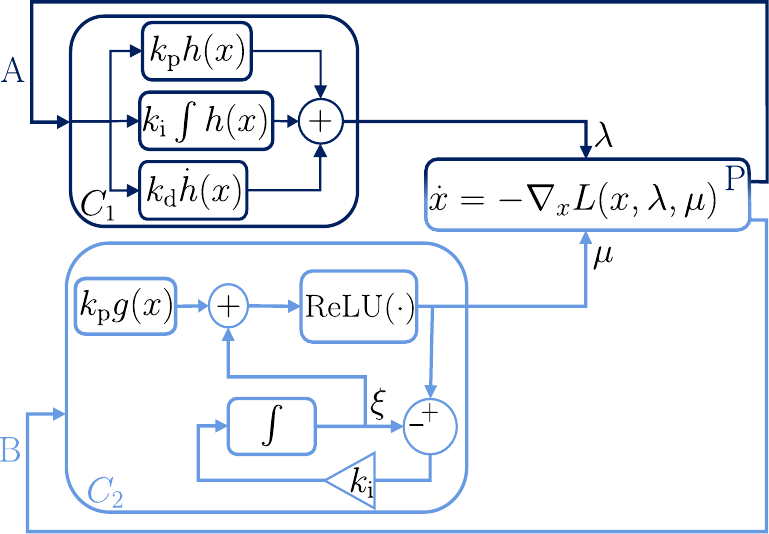}
    \caption{Equivalent SPPID Dynamics. The plant $\textup{P}$ implements the primal gradient flow with state $x(t) \in \mathbb{R}^n$, control inputs $\lambda(t) \in \mathbb{R}^p$, $\mu(t) \in \mathbb{R}^m_{\geq 0}$ and measured outputs $y_1(t) = h(x(t)) \in \mathbb{R}^p$ and $y_2(t) = g(x(t)) \in \mathbb{R}^m$. The controller $\textup{C}_1$ implements the PID law~\eqref{eq:controller_eq} (Panel A), whereas the controller $\textup{C}_2$ implements anti-windup PI law~\eqref{eq:controller_ineq} (Panel B).
    }
    \label{fig:block_diagram}
\end{figure}
\end{arxiv}%
\subsubsection*{Controller for Equality Constraints}
While~\cite{VC-SMF-SP-DR:25} employs a PI controller for the residual $y_1 = h(x)$, yielding the \emph{PI-controlled multiplier optimization} (PI-CMO), we extend this approach to a full PID scheme. Specifically, the PID controller for $\lambda$, illustrated in Figure~\ref{fig:block_diagram}, Panel A, is
\beq
\label{eq:controller_eq}
\lambda = \subscr{k}{i}^{\textup{eq}} \int_{0}^t h(x) d\tau + \subscr{k}{p}^{\textup{eq}} h(x) + \subscr{k}{d}^{\textup{eq}} J_h(x)\dot{x},
\eeq
where $\subscr{k}{i}^{\textup{eq}}$, $\subscr{k}{p}^{\textup{eq}}$, $\subscr{k}{d}^{\textup{eq}} \in\R_{\ge 0}$ are the integral, proportional, and derivative gains, respectively.
Let $\nu \in \R^p$ denote the controller state associated with integral action, satisfying $\dot{\nu} = \subscr{k}{i}^{\textup{eq}} h(x)$.
Each term in~\eqref{eq:controller_eq} plays a structurally distinct role:
\begin{arxiv}
\begin{enumerate}[label=\textup{(\roman*)}]
    \item the \emph{integral term} $\subscr{k}{i}^{\textup{eq}} \int_0^t h(x(\tau)) d \tau$ enforces constraint satisfaction by accumulating constraint violations and driving them to zero; this requires the gain to be strictly positive;
    \item the \emph{proportional term} $\subscr{k}{p}^{\textup{eq}} h(x)$ shapes the energy landscape by introducing an augmented term to the Lagrangian;
    \item the \emph{derivative term} $\subscr{k}{d}^{\textup{eq}} J_h(x)\dot{x}$ modifies the geometry of the primal dynamics by inducing a state-dependent metric.
\end{enumerate}
\end{arxiv}%

\subsubsection*{Controller for Inequality Constraints}
A pure integral controller $\dot{\xi} = \subscr{k}{i}^{\textup{in}} g(x)$ is unsuitable for the inequality constraints for two main reasons. First, whenever a constraint is strictly inactive complementary slackness requires the corresponding multiplier to vanish; however, the integrator state decrease without bound, a phenomenon known as \emph{integrator windup}. Second, the condition $\mu \geq \0_m$ is not naturally enforced.
To address these issues, we interpret dual non-negativity as a saturation constraint and employ an anti-windup PI controller in which the saturated output is fed back into the integrator
\beq
\label{eq:controller_ineq}
\mu = \relu\big(\xi + \subscr{k}{p}^{\textup{in}} g(x)\big), 
\qquad
\dot{\xi} = \subscr{k}{i}^{\textup{in}} (\mu - \xi),
\eeq
where $\xi \in \R^m$ is the internal state, $\subscr{k}{i}^{\textup{in}}, \subscr{k}{p}^{\textup{in}} \in \R_{>0}$ are the integral and proportional gains, respectively. Figure~\ref{fig:block_diagram}, Panel B, shows the resulting controller.
The considered anti-windup PI controller operates as follows.
When unsaturated, $\mu = \xi + \subscr{k}{p}^{\textup{in}} g(x)$ is a PI controller, and $\dot{\xi} = \subscr{k}{i}^{\textup{in}} \subscr{k}{p}^{\textup{in}} g(x)$ reduces to an integrator accumulating constraint violation.
When saturated, $\mu = \0_m$ and $\dot{\xi} = - \subscr{k}{i}^{\textup{in}} \xi$, which drives $\xi$ back to zero, thereby preventing indefinite drift.
A key structural property of~\eqref{eq:controller_ineq} is that it preserves dual non-negativity, as formalized next.
\begin{arxiv}
\begin{mylemma}[label=lem:forward_invariance_dual_orthant]{Forward invariance of the nonnegative orthant}
Consider the controller~\eqref{eq:controller_ineq} and let $x(t)$ be a trajectory of~\eqref{eq:open_system_eq_ineq}. If $\xi(0) \in \R^m_{\geq 0}$, then $\xi(t) \in \R^m_{\geq 0}$ for all $t \geq 0$.
\end{mylemma}
\begin{proof}
Let $\subscr{f}{c}$ denote the right-hand side of the $\xi$-dynamics in~\eqref{eq:controller_ineq}. To prove the statement we show that $\R_{\geq 0}^m$ is a forward invariant set for $\subscr{f}{c}$. By applying Nagumo's Theorem~\cite{MN:1942}, this holds if and only if
\beq
\label{eq:nagumo_pos}
\subscr{f}{c,i}(\xi) \geq 0 \quad \forall \xi \in \R_{\geq 0}^m \text{ such that } \xi_i =0.
\eeq
Consider~\eqref{eq:controller_ineq} in components $\dot{\xi}_i = - \subscr{k}{i}^{\textup{in}} \xi_i + \subscr{k}{i}^{\textup{in}} \relu\big(\xi_i + \subscr{k}{p}^{\textup{in}} g_i(x)\big) = \subscr{f}{c,i}(\xi)$, $i \in \until{m}$.
Writing condition~\eqref{eq:nagumo_pos} for $\subscr{f}{c,i}$, for all $\xi \in \R_{\geq 0}^m$ such that $\xi_i =0$ we have
$
\subscr{f}{c,i}(\xi) = \subscr{k}{i}^{\textup{in}}\relu\big(\xi_i + \subscr{k}{p}^{\textup{in}} g_i(x)\big) \geq 0,
$
since $\relu(\cdot) \geq 0$ for all arguments. This concludes the proof.
\end{proof}
\end{arxiv}%

\subsubsection*{SPPID Dynamics}
We now close the loop by applying the controllers~\eqref{eq:controller_eq}
and~\eqref{eq:controller_ineq} to the system~\eqref{eq:open_system_eq_ineq}. Since $\lambda$ depends on $\dot x$, and $\mu$ is non-smooth at the saturation boundary, we express the closed loop through the integral states $(\nu,\xi)$, yielding
\begin{arxiv}
\beq
\label{eq:sp_pid_eq_ineq}
\begin{cases}
\dot{x} = - M(x)^{-1} \Bigl( \nabla f(x) + J_h(x)^{\top} \big(\nu + \subscr{k}{p}^{\textup{eq}}h(x)\big) + J_g(x)^{\top} \relu\big(\xi + \subscr{k}{p}^{\textup{in}} g(x)\big)\Bigr), \\
\dot{\nu} = \subscr{k}{i}^{\textup{eq}} h(x), \\
\dot{\xi} = \subscr{k}{i}^{\textup{in}} \left( \relu\big(\xi + \subscr{k}{p}^{\textup{in}} g(x)\big) - \xi \right).
\end{cases}
\eeq
\end{arxiv}%
where $M(x):= I_n + \subscr{k}{d}^{\textup{eq}} J_h(x)^{\top} J_h(x)$ is symmetric and positive definite.
We refer to the dynamics~\eqref{eq:sp_pid_eq_ineq} as the \emph{saddle-point PID (SPPID) dynamics}.

\begin{remark}[Derivative action and inequality constraints]
\label{rem:kd_inequality}
In the equality-constrained case, introducing a derivative term through $J_g(x)\dot x$ allows for the explicit inversion of the matrix $M(x)$, inducing a state-dependent Riemannian metric. 
In the inequality case, however, the derivative term interacts non-smoothly with the $\relu (\cdot)$ operator enforcing dual feasibility, generating a non-smooth implicit algebraic loop. Resolving this algebraic loop and, in turn, designing a PID scheme for inequality-constrained problems, remains an open problem.
\end{remark}

The role of proportional action can be formalized via the augmented Lagrangian. Given $\subscr{k}{p}^{\textup{eq}} \geq 0$ and $\subscr{k}{p}^{\textup{in}} > 0$, consider the augmented Lagrangian associated with~\eqref{eq:eq_ineq}~\cite{DPB:97}
\begin{arxiv}
\beq
\label{eq:aug_lagrangian_eq_ineq}
\subscr{L}{aug}(x, \nu, \xi) = f(x) + \nu^{\top} h(x) + \frac{\subscr{k}{p}^{\textup{eq}}}{2} \|h(x)\|^2 + \frac{1}{2\subscr{k}{p}^{\textup{in}}} \sum_{i=1}^m \left( \left[\relu\big(\xi_i + \subscr{k}{p}^{\textup{in}} g_i(x)\big) \right]^2 - \xi_i^2 \right).
\eeq
The inequality penalty term encodes complementary slackness: the $i$-th summand penalizes the constraint only when it is active or violated, and vanishes when $g_i(x) < 0$ and $\xi_i = 0$.
\end{arxiv}%

The following theorem establishes that equilibria of SPPID coincide with KKT points of~\eqref{eq:eq_ineq} and that the dynamics admit a primal--dual gradient-flow interpretation.
\begin{arxiv}
\begin{mytheorem}
[label=thm:sp_pid_eq_ineq]{Equilibria and saddle-flow structure of SPPID}
Consider the nonlinear program~\eqref{eq:eq_ineq} and the SPPID dynamics~\eqref{eq:sp_pid_eq_ineq} with $\subscr{k}{i}^{\textup{eq}}, \subscr{k}{i}^{\textup{in}}, \subscr{k}{p}^{\textup{in}} >0$. Then,
\begin{enumerate}[label=\textup{(\roman*)}]
\item
\label{thm_pi_eq_ineq:item1}
A point $(x^\star, \nu^\star, \xi^\star)$ is an equilibrium of~\eqref{eq:sp_pid_eq_ineq} if and only if there exist multipliers $\lambda^\star \in \R^p$ and $\mu^\star \in \R^m_{\ge 0}$ such that $(x^\star, \lambda^\star, \mu^\star)$ satisfies the KKT conditions for~\eqref{eq:eq_ineq}, and $\lambda^\star = \nu^\star$ and $\mu^\star = \xi^\star$.
\item 
\label{thm_pi_eq_ineq:item2}
The dynamics~\eqref{eq:sp_pid_eq_ineq} are a preconditioned primal--dual gradient flow of the augmented Lagrangian~\eqref{eq:aug_lagrangian_eq_ineq}, namely
\beq
\label{eq:pd_gradient_flow_eq_ineq}
\begin{cases}
\dot{x} = - M(x)^{-1}\nabla_x \subscr{L}{aug}(x, \nu, \xi), \\
\dot{\nu} = \subscr{k}{i}^{\textup{eq}} \nabla_\nu \subscr{L}{aug}(x, \nu, \xi), \\
\dot{\xi} = \subscr{k}{i}^{\textup{in}} \subscr{k}{p}^{\textup{in}} \nabla_\xi \subscr{L}{aug}(x, \nu, \xi).
\end{cases}
\eeq
\end{enumerate}
\end{mytheorem}
\end{arxiv}%
\begin{proof}
Let $z := (x, \nu, \xi)$ and let $z^\star := (x^\star, \nu^\star, \xi^\star)$ be an equilibrium of~\eqref{eq:sp_pid_eq_ineq}. The condition $\dot{\nu} = \0_p$ gives $h(x^\star) = \0_p$, establishing~\eqref{kkt_stationarity_eq}. Then, evaluating the equality controller~\eqref{eq:controller_eq} at equilibrium yields $\lambda^{\star} = \nu^{\star}$.
The condition $\dot{\xi} = \0_m$ gives $\xi^\star = \relu \bigl(\xi^\star + \subscr{k}{p}^{\textup{in}} g(x^\star)\bigr) = \mu^\star$. By non-negativity of the $\relu$ operator, this fixed-point condition automatically satisfies dual feasibility~\eqref{kkt_dual_feasibility}. We now show this condition is also equivalent to primal feasibility~\eqref{kkt_primal_feasibility} and complementary slackness~\eqref{kkt_complementarity} by analyzing it component-wise for each $i \in \{1, \dots, m\}$. If $\mu_i^\star = 0$, then substituting into the fixed-point condition gives $0 = \relu\big(\subscr{k}{p}^{\textup{in}} g_i(x^\star)\big)$, which implies $g_i(x^\star) \le 0$. Thus,~\eqref{kkt_primal_feasibility} and~\eqref{kkt_complementarity} hold. If instead $\mu_i^\star > 0$, then $\mu_i^\star = \mu_i^\star + \subscr{k}{p}^{\textup{in}} g_i(x^\star)$, which implies $g_i(x^\star) = 0$. Thus,~\eqref{kkt_primal_feasibility} and~\eqref{kkt_complementarity} hold.
Finally, $\dot{x} = \0_n$ and $M(x^{\star})\succ 0$ yield $\nabla f(x^\star) +  J_h(x^\star)^{\top} \lambda^\star  + J_g(x^\star)^{\top} \mu^\star = \0_n$, that is, stationarity~\eqref{kkt_stationarity}. Therefore, $z^\star$ is a KKT point. Moreover $\nu^\star = \lambda^\star$ and $\xi^\star = \mu^\star$. The converse follows by noticing that any KKT point $z^\star$ with $\lambda^\star = \nu^\star$ and $\mu^\star = \xi^\star$ annihilates all three vector fields in~\eqref{eq:sp_pid_eq_ineq}. This concludes the proof of item~\ref{thm_pi_eq_ineq:item1}.

Next, since $\subscr{L}{aug}(z)$ is continuously differentiable in all
arguments~\cite{DPB:97}, we compute its partial gradients directly.
Recalling that the derivative of $\phi(s) = \frac{1}{2} [\relu(s)]^2$ is $\phi'(s) = \relu(s)$, the chain rule yields $\nabla_x \subscr{L}{aug}(z) = \nabla f(x) + J_h(x)^{\top} \big(\nu + \subscr{k}{p}^{\textup{eq}}h(x)\big) + J_g(x)^{\top} \relu\bigl(\xi + \subscr{k}{p}^{\textup{in}} g(x)\bigr)$.
Since $M(x)\succ 0$, the equation $M(x)\dot x = -\nabla_x \subscr{L}{aug}(z)$ admits the unique solution $\dot x = - M(x)^{-1}\nabla_x \subscr{L}{aug}(z)$.
Moreover, $\nabla_{\nu} \subscr{L}{aug}(z) = h(x)$ and 
$
\nabla_\xi \subscr{L}{aug}(z) = \frac{1}{\subscr{k}{p}^{\textup{in}}} \bigl(\relu\big(\xi + \subscr{k}{p}^{\textup{in}} g(x)\big) - \xi \bigr).
$
Substituting these expressions back into~\eqref{eq:pd_gradient_flow_eq_ineq} recovers the SPPID dynamics~\eqref{eq:sp_pid_eq_ineq}. This concludes the proof.
\end{proof}
From the perspective of the dual variables, it is important to note that along transient trajectories, $\nu$ and $\xi$ do not in general coincide with the classic Lagrange multipliers $\lambda$ and $\mu$ of problem~\eqref{eq:eq_ineq}. Instead, $\nu$ and $\xi$ represent shifted versions of these multipliers (with $\xi$ additionally undergoing a saturation) and play the role of the dual variables for $\subscr{L}{aug}$. At equilibrium, however, this transient discrepancy vanishes, and the states precisely recover the original Lagrange multipliers.

Theorem~\ref{thm:sp_pid_eq_ineq} provides the conceptual foundation of the proposed framework. Besides establishing that the equilibria of SPPID coincide with the KKT points of the original problem, it shows that the proposed dynamics admit a saddle-point flow interpretation. In particular, item~\ref{thm_pi_eq_ineq:item2} shows that the primal dynamics correspond to Riemannian gradient descent on $\subscr{L}{aug}$ with respect to the metric induced by $M(x)$, whereas the dual dynamics correspond to a preconditioned Euclidean gradient ascent. Notably, this Riemannian structure arises solely from the derivative action in the equality controller~\eqref{eq:controller_eq}, providing a geometric interpretation of derivative feedback that does not emerge from standard primal--dual formulations of $\subscr{L}{aug}$. The remainder of the paper builds on this characterization to establish convergence properties and study special cases.

\subsection{Convergence Analysis of SPPID in Convex Settings}
\label{sec:convergence_convex}
We study the convergence properties of SPPID in convex settings. Namely, we assume that the objective function and the inequality constraints are convex, and the equality constraints are affine, i.e., $h(x) := Ax - b$, where $A \in \R^{p \times n}$ and $b \in \R^p$.
\begin{arxiv}
The SPPID~\eqref{eq:sp_pid_eq_ineq} becomes, with $M := I_n + \subscr{k}{d}^{\textup{eq}} A^{\top} A$,
\beq
\label{eq:sp_pid_eq_ineq_convex}
\begin{cases}
\dot{x} = - M^{-1} \Bigl( \nabla f(x) + A^{\top} \big(\nu + \subscr{k}{p}^{\textup{eq}}\bigl(Ax - b\bigr)\big) + J_g(x)^{\top} \relu\bigl(\xi + \subscr{k}{p}^{\textup{in}} g(x)\bigr)\Bigr), \\
\dot{\nu} = \subscr{k}{i}^{\textup{eq}} \bigl(Ax - b\bigr), \\
\dot{\xi} =  - \subscr{k}{i}^{\textup{in}} \xi + \subscr{k}{i}^{\textup{in}} \relu\bigl(\xi + \subscr{k}{p}^{\textup{in}} g(x)\bigr).
\end{cases}
\eeq
\end{arxiv}%

The convergence analysis relies on the following standard assumptions on $f$, $g$, and $A$.
\begin{arxiv}
\begin{myassumption}{Convexity and constraint structure}
\label{ass:1}
Given $\map{f}{\R^n}{\R}$, $\map{g}{\R^n}{\R^m}$, $A \in \R^{p \times n}$, $b \in \R^{p}$, assume
\begin{enumerate}[label=\textup{(\roman*)}]
\item
\label{ass:1_f_g}
the functions $f$ and $g_i$, $i \in \{1, \dots, m\}$, are convex and continuously differentiable;
\item
\label{ass:3_slater}
there exists a vector $\bar{x}$ such that $A\bar x = b$ and $g(\bar{x}) < \0_m$ (Slater's condition).
\end{enumerate}
\end{myassumption}
\smallskip
\end{arxiv}%
Under Assumption~\ref{ass:1}, problem~\eqref{eq:eq_ineq} is a convex program.
Moreover, Slater's condition~\ref{ass:3_slater} guarantees that the KKT conditions are necessary and sufficient for optimality and that $\mathcal{K}$ is nonempty~\cite{DPB:97}.
The next result characterizes the convergence properties of SPPID under the above assumptions.
\begin{arxiv}
\begin{mytheorem}[label=thm:gas_eq_ineq]{Convergence to the KKT set}
Under Assumption~\ref{ass:1}, consider SPPID~\eqref{eq:sp_pid_eq_ineq}. Then, for every $\subscr{k}{d}^{\textup{eq}} \geq 0$, $\subscr{k}{p}^{\textup{eq}}, \subscr{k}{i}^{\textup{eq}}$, $\subscr{k}{p}^{\textup{in}}$, $\subscr{k}{i}^{\textup{in}} > 0$, the following hold
\begin{enumerate}[label=\textup{(\roman*)}]
\item
\label{thm_conv_eq_ineq:item1}
every trajectory converges to the set of KKT points $\mathcal{K}$ of problem~\eqref{eq:eq_ineq}.
\item
\label{thm_conv_eq_ineq:item2}
if additionally $f$ is strictly convex and LICQ holds, then~\eqref{eq:eq_ineq} has a unique KKT point $z^\star := (x^\star, \nu^\star, \xi^\star)$ and $z^\star$ is globally asymptotically stable.
\end{enumerate}
\end{mytheorem}
\begin{proof}
Let $z := (x, \nu, \xi)$ and let $z^\star$ be a KKT point of~\eqref{eq:eq_ineq} guaranteed to exist by Assumption~\ref{ass:1}. Consider the following Lyapunov function candidate
$$
V (x, \nu, \xi) = \frac{1}{2}\norm{x - x^\star}^2_{M} + \frac{1}{2 \subscr{k}{i}^{\textup{eq}}}\norm{\nu - \nu^\star}^2 + \frac{1}{2 \subscr{k}{i}^{\textup{in}} \subscr{k}{p}^{\textup{in}}}\norm{\xi - \xi^\star}^2.
$$
Clearly, $V(x, \nu, \xi)$ is positive definite, radially unbounded, and vanishes only at $(x^\star, \nu^\star, \xi^\star)$. Throughout, we write $r(x,\xi) := \relu\bigl(\xi + \subscr{k}{p}^{\textup{in}} g(x)\bigr)$ for brevity.

\textbf{Computing $\dot{V}$.}
Differentiating along trajectories of~\eqref{eq:sp_pid_eq_ineq_convex} gives
\begin{align*}
\dot{V}(x, \nu, \xi) &= (x - x^\star)^{\top} M \dot x + \frac{1}{\subscr{k}{i}^{\textup{eq}}} (\nu - \nu^\star)^{\top} \dot \nu + \frac{1}{\subscr{k}{i}^{\textup{in}} \subscr{k}{p}^{\textup{in}}} (\xi - \xi^\star)^{\top} \dot \xi \\
&= (x - x^\star)^{\top} \Bigl(- \nabla f(x) - J_g(x)^{\top} r(x, \xi) - A^{\top} \nu - \subscr{k}{p}^{\textup{eq}}A^{\top} h(x)\Bigr)\\ 
&\quad + \frac{1}{\subscr{k}{i}^{\textup{in}} \subscr{k}{p}^{\textup{in}}} (\xi - \xi^\star)^{\top} \Bigl( - \subscr{k}{i}^{\textup{in}} \xi + \subscr{k}{i}^{\textup{in}} r(x, \xi) \Bigr)+ (\nu - \nu^\star)^{\top} h(x)\\
&= - (x - x^\star)^{\top} \nabla f(x) - (x - x^\star)^{\top} J_g(x)^{\top} r(x, \xi) - (x - x^\star)^{\top} A^{\top} \bigl(\nu + \subscr{k}{p}^{\textup{eq}}h(x)\bigr)\\
&\quad + (\nu - \nu^\star)^{\top} h(x) - \frac{1}{\subscr{k}{p}^{\textup{in}}} (\xi - \xi^\star)^{\top} \xi +  \frac{1}{\subscr{k}{p}^{\textup{in}}} (\xi - \xi^\star)^{\top} r(x, \xi).
\end{align*}
Using the stationarity condition~\eqref{kkt_stationarity}, $\nabla f(x^\star) = - J_g(x^\star)^{\top} \xi^\star  - A^{\top} \nu^\star $, and adding and subtracting $(x-x^\star)^{\top} \nabla f(x^\star)$ to $\dot{V}(x, \nu, \xi)$ yields
\begin{align}
\dot{V}(x, \nu, \xi) &= - (x - x^\star)^{\top} \bigl(\nabla f(x) - \nabla f(x^\star)\bigr) + (x - x^\star)^{\top} J_g(x^\star)^{\top} \xi^\star + (x - x^\star)^{\top} A^{\top} \nu^\star \nonumber\\
&\quad - (x - x^\star)^{\top} J_g(x)^{\top} r(x, \xi) - \frac{1}{\subscr{k}{p}^{\textup{in}}} (\xi - \xi^\star)^{\top} \xi +  \frac{1}{\subscr{k}{p}^{\textup{in}}} (\xi - \xi^\star)^{\top} r(x, \xi) \nonumber \\
& \quad - (x - x^\star)^{\top} A^{\top} \bigl(\nu + \subscr{k}{p}^{\textup{eq}}h(x)\bigr) + (\nu - \nu^\star)^{\top} h(x) \nonumber\\
&= - (x - x^\star)^{\top} \bigl(\nabla f(x) - \nabla f(x^\star)\bigr) + (x - x^\star)^{\top} \bigl(J_g(x^\star)^{\top} \xi^\star - J_g(x)^{\top} r(x, \xi)\bigr) - \frac{1}{\subscr{k}{p}^{\textup{in}}} (\xi - \xi^\star)^{\top} \bigl(\xi - r(x, \xi)\bigr) \nonumber \\
& \quad - (x - x^\star)^{\top} A^{\top} \bigl(\nu -\nu^\star + \subscr{k}{p}^{\textup{eq}}h(x)\bigr) + (\nu - \nu^\star)^{\top} h(x) \nonumber\\ 
&= - (x - x^\star)^{\top} \bigl(\nabla f(x) - \nabla f(x^\star)\bigr) + (x - x^\star)^{\top} \bigl(J_g(x^\star)^{\top} \xi^\star - J_g(x)^{\top} r(x, \xi)\bigr) - \frac{1}{\subscr{k}{p}^{\textup{in}}} (\xi - \xi^\star)^{\top} \bigl(\xi - r(x, \xi)\bigr) \nonumber \\
& \quad - \subscr{k}{p}^{\textup{eq}}(x - x^\star)^{\top} A^{\top} A\bigl(x - x^\star\bigr) \nonumber \\
&= - (x - x^\star)^{\top} \bigl(\nabla f(x) - \nabla f(x^\star)\bigr) + (x - x^\star)^{\top} \bigl(J_g(x^\star)^{\top} \xi^\star - J_g(x)^{\top} r(x, \xi)\bigr) \nonumber\\
& \quad - \frac{1}{\subscr{k}{p}^{\textup{in}}} (\xi - \xi^\star)^{\top} \bigl(\xi - r(x, \xi)\bigr) - \subscr{k}{p}^{\textup{eq}}\norm{h(x)}^2, \label{v_ineq_1_}
\end{align}
where in the last equalities we used the fact that $h(x) = Ax -b = Ax - A x^\star = A(x - x^\star)$.

\textbf{Bounding the Jacobian cross-term.}
Since $g$ is convex and continuously differentiable, the first-order inequality $g(z) - g(y) \geq J_g(y)(z-y)$ holds for all $z,y$.
Evaluating the above inequality at $(z, y) = (x, x^\star)$ and $(z, y) = (x^\star, x)$, respectively, yields
\begin{align*}
g(x) - g(x^\star) \geq J_g(x^\star) (x - x^\star) \quad &\implies \quad (g(x) - g(x^\star))^{\top} \xi^\star \geq \big(J_g(x^\star)(x - x^\star)\big)^{\top} \xi^\star, \\
g(x^\star) - g(x) \geq J_g(x) (x^\star - x) \quad &\implies \quad -\big(g(x) - g(x^\star)\bigr)^{\top} r(x, \xi) \geq -\big(J_g(x)(x - x^\star)\big)^{\top} r(x, \xi),
\end{align*}
where the inequalities are preserved under multiplication since $\xi^\star \geq \0_m$ by KKT conditions and $r(x, \xi) \geq \0_m$ by the definition of the $\relu$ operator.
Substituting the above inequalities into~\eqref{v_ineq_1_} gives
\begin{align*}
\dot{V}(x, \xi) &\leq - (x - x^\star)^{\top} \bigl(\nabla f(x) - \nabla f(x^\star)\bigr)+ \bigl(g(x) - g(x^\star)\bigr)^{\top} \bigl(\xi^\star - r(x, \xi)\bigr) - \frac{1}{\subscr{k}{p}^{\textup{in}}} (\xi - \xi^\star)^{\top} \bigl(\xi - r(x, \xi)\bigr)- \subscr{k}{p}^{\textup{eq}}\norm{h(x)}^2.
\end{align*}

\textbf{Applying properties of $\relu$.}
The firmly non-expansive property of the $\relu$ operator allows us to establish the following bound on the constraint interaction terms:
\begin{align}
\bigl(g(x) - g(x^\star)\bigr)^{\top}& \bigl(\xi^\star - r(x, \xi)\bigr) \leq - \frac{1}{\subscr{k}{p}^{\textup{in}}} \norm{r(x, \xi) - \xi^\star}^2 + \frac{1}{\subscr{k}{p}^{\textup{in}}} \bigl(\xi - \xi^\star \bigr)^{\top} \bigl(r(x, \xi) - \xi^\star \bigr). \label{eq:relu_nonexpansive_bound_}
\end{align}
To see this, recall that firmly non-expansiveness of the $\relu$ implies 
$\norm{\relu(z) - \relu(y)}^2 \leq \bigl(z - y\bigr)^{\top}\bigl(\relu(z) - \relu(y) \bigr)$, for all $z,y \in \R^m$.
Inequality~\eqref{eq:relu_nonexpansive_bound_} follows directly by selecting $z = \xi + \subscr{k}{p}^{\textup{in}} g(x)$ and $y = \xi^\star + \subscr{k}{p}^{\textup{in}} g(x^\star)$.  

\textbf{Combining the bounds.}
Then, we can bound $\dot{V}(x, \nu, \xi)$ obtaining
\begin{align*}
- &(x - x^\star)^{\top} \bigl(\nabla f(x) - \nabla f(x^\star)\bigr) - \frac{1}{\subscr{k}{p}^{\textup{in}}} \norm{r(x, \xi) - \xi^\star}^2 + \frac{1}{\subscr{k}{p}^{\textup{in}}} \bigl(\xi - \xi^\star \bigr)^{\top} \bigl(r(x, \xi) - \xi^\star \bigr) - \frac{1}{\subscr{k}{p}^{\textup{in}}} (\xi - \xi^\star)^{\top} \bigl(\xi - r(x, \xi)\bigr)\\
&= - (x - x^\star)^{\top} \bigl(\nabla f(x) - \nabla f(x^\star)\bigr) + \frac{1}{\subscr{k}{p}^{\textup{in}}} \Bigl(- \norm{r(x, \xi)}^2 - \norm{\xi^\star}^2 + 2r(x, \xi)^{\top} \xi^\star + 
\xi^{\top} r(x, \xi) - r(x, \xi)^{\top} \xi^\star - \xi^{\top} \xi^\star + \norm{\xi^\star}^2\\
&\quad \quad \quad \quad \quad\quad \quad \quad \quad \quad\quad \quad \quad \quad \quad \quad \quad - \norm{\xi}^2 + \xi^{\top} \xi^\star + \xi^{\top} r(x, \xi) - r(x, \xi)^{\top} \xi^\star\Bigr)\\
&= - (x - x^\star)^{\top} \bigl(\nabla f(x) - \nabla f(x^\star)\bigr) + \frac{1}{\subscr{k}{p}^{\textup{in}}} \Bigl(- \norm{r(x, \xi)}^2 - \norm{\xi}^2 + 2 \xi^{\top} r(x, \xi) \Bigr) \\
&= - (x - x^\star)^{\top} \bigl(\nabla f(x) - \nabla f(x^\star)\bigr) - \frac{1}{\subscr{k}{p}^{\textup{in}}} \norm{r(x, \xi) - \xi}^2 \leq 0, 
\end{align*}
where the last inequality follows from $(x - x^\star)^{\top} \bigl(\nabla f(x) - \nabla f(x^\star)\bigr) \geq 0$, by convexity of $f$.
In summary, we have shown that
\begin{equation}
\label{eq:Vdot_final_}
\dot{V}(x, \nu, \xi) \leq - (x - x^\star)^{\top} \bigl(\nabla f(x) - \nabla f(x^\star)\bigr) - \frac{1}{\subscr{k}{p}^{\textup{in}}} \norm{r(x, \xi) - \xi}^2 - \subscr{k}{p}^{\textup{eq}}\norm{h(x)}^2 \leq 0.
\end{equation}
Hence, $V$ is non-increasing along trajectories of~\eqref{eq:sp_pid_eq_ineq_convex}, every trajectory is bounded, and $(x^\star,\nu^\star, \xi^\star)$ is stable.

\textbf{LaSalle's invariance principle.}
Let $\mathcal{E} := \setdef{(x,\nu,\xi)\in\R^{n+p+m}}{\dot{V}(x,\nu,\xi)=0}$ and let $\Omega\subseteq\mathcal{E}$ denote the largest invariant subset of $\mathcal{E}$. Consider a trajectory $(x(t),\nu(t),\xi(t))\in\Omega$. From~\eqref{eq:Vdot_final_}, we have
\begin{align}
&(x-x^\star)^{\top} \bigl(\nabla f(x) -\nabla f(x^\star)\bigr) = 0, \label{eq:vdot_primal_}\\
&\relu\bigl(\xi+\subscr{k}{p}^{\textup{in}}g(x)\bigr) - \xi = \0_m, \label{eq:vdot_dual_1_}\\
&h(x) = \0_p. \label{eq:vdot_dual_2_}
\end{align}
Equations~\eqref{eq:vdot_dual_1_} and~\eqref{eq:vdot_dual_2_} imply $\dot{\xi} = \0_m$ and $\dot{\nu} = \0_p$, respectively. Therefore $\xi(t) = \bar{\xi}$ and $\nu(t) = \bar{\nu}$, for all $t \geq 0$.
Moreover, the fixed-point condition~\eqref{eq:vdot_dual_1_} and~\eqref{eq:vdot_dual_2_} imply primal feasibility, dual feasibility, and complementary slackness~\eqref{kkt_primal_feasibility}--\eqref{kkt_complementarity} in $(x, \bar{\nu}, \bar{\xi})$, namely, $h(x) = \0_p$, $g(x) \leq \0_m$, $\bar{\xi} \geq \0_m$, and $\bar{\xi}^{\top} g(x) = 0$.
Then, the primal dynamics reduce to $\dot{x} = - M^{-1}\Bigl(\nabla f(x) + J_h(x)^{\top} \bar{\nu} + J_g(x)^{\top} \bar{\xi}\Bigr)$. Define the Lagrangian $L(x,\nu,\xi) = f(x) + h(x)^{\top} \nu + g(x)^{\top} \xi$, which is convex due to $f$, $g$ and $h$ convex.
The primal dynamics read $\dot{x} = - M^{-1}\nabla L(x,\bar{\nu},\bar{\xi})$. To understand the invariant set for the condition~\eqref{eq:vdot_primal_}, consider the first term of the Lyapunov candidate $\phi(x) = 1/2 \|x - x^\star\|^2_M$. Along trajectories in $\Omega$, we have
\begin{align*}
\dot{\phi}(x) &= (x - x^\star)^{\top} M \dot{x} = -(x - x^\star)^{\top} \nabla L(x,\bar{\nu},\bar{\xi})  \leq L(x^\star,\bar{\nu},\bar{\xi}) - L(x,\bar{\nu},\bar{\xi}) \leq f(x^\star) - f(x) \leq 0,
\end{align*}
where the first inequality follows from convexity of $L(\cdot,\bar{\nu},\bar{\xi})$, the second uses $L(x^\star,\bar{\nu},\bar{\xi}) \leq f(x^\star)$ due to $h(x^\star)=\0_p$, $\bar{\xi}\geq \0_m$ and $g(x^\star)\leq \0_m$, and $L(x,\bar{\nu},\bar{\xi}) = f(x)$ due to complementary slackness and~\eqref{eq:vdot_dual_2_}. The last inequality follows from the optimality of $x^\star$ and the feasibility of $x$.
So $\phi(x)$ is nonincreasing along trajectories in $\Omega$, and the feasible set is forward invariant, implying all solutions converge to the largest invariant subset of $\{\dot{\phi}=0\}$.
On this invariant set, $0=\dot{\phi} \leq f(x^\star)-f(x) \leq 0$, which implies $f(x^\star) = f(x)$. Thus, $x$ minimizes $L(\cdot, \bar{\nu}, \bar\xi)$ by convex sufficiency, giving $\nabla f(x) + J_h(x)^{\top}\bar{\nu} + J_g(x)^{\top}\bar{\xi} = \0_n$, implying $\dot{x} = \0_n$.
Therefore, every point in $\Omega$ satisfies the KKT conditions and, by LaSalle's invariance principle, the $\omega$-limit set of every solution is contained in $\Omega$, and we conclude $\operatorname{dist}((x, \nu, \xi), \mathcal{K}) \to 0$.
This proves item~\ref{thm_conv_eq_ineq:item1}.

Item~\ref{thm_conv_eq_ineq:item2} follows by noticing that strict convexity of $f$ and LICQ together imply that the KKT set reduces to the singleton  $\{(x^\star, \nu^\star, \xi^\star)\}$. Global asymptotic stability then follows from LaSalle's invariance principle combined with radial unboundedness of $V$. This concludes the proof.
\end{proof}
\end{arxiv}%
With additional regularity at $z^\star$, convergence is exponential after an initial phase of linearly bounded decay.
\begin{arxiv}
\begin{mytheorem}[label=thm:lin_exp_eq_ineq]{Linear-exponential convergence of SPPID}
Let Assumption~\ref{ass:1} hold and let $z^\star = (x^\star, \nu^\star, \xi^\star)$ be a KKT point of~\eqref{eq:eq_ineq}. If strict complementarity, LICQ, and SOSC hold at $z^\star$, then $z^\star$ unique and, for every $\subscr{k}{d}^{\textup{eq}} \geq 0$ and $\subscr{k}{p}^{\textup{eq}}$, $\subscr{k}{i}^{\textup{eq}}$, $\subscr{k}{p}^{\textup{in}}$, $\subscr{k}{i}^{\textup{in}} > 0$, every trajectory of~\eqref{eq:sp_pid_eq_ineq} linearly-exponentially converges to $z^\star$.
\end{mytheorem}
\begin{proof}
Since~\eqref{eq:eq_ineq} is convex under Assumption~\ref{ass:1}, and strict complementarity, LICQ, and SOSC hold at $z^\star$, $x^\star$ is a strict local minimizer of~\eqref{eq:eq_ineq}~\cite{DPB:97}. Moreover, since a strict local minimizer of a convex program is its unique global minimizer, and LICQ ensures unique multipliers, $z^\star$ is the unique KKT point of~\eqref{eq:eq_ineq}, i.e., $\K = \{z^\star\}$.
We show that SPPID is globally weakly contracting (non-expansive), and locally strongly contracting in a neighborhood of $z^\star$.
Define $P = \operatorname{diag}\bigl(M, \frac{1}{\subscr{k}{i}^{\textup{eq}}} I_p, \frac{1}{\subscr{k}{i}^{\textup{in}} \subscr{k}{p}^{\textup{in}}} I_m\bigr) \succ 0$.
Moreover, let $y := \xi + \subscr{k}{p}^{\textup{in}} g(x)$ and define $G(y) := \nabla \relu(y)$ where it exists. By Rademacher's Theorem, the Jacobian $J(z)$ of the SPPID~\eqref{eq:sp_pid_eq_ineq_convex} exists almost everywhere and is given by
$$
J(z)= -
\begin{bmatrix}
M^{-1}\left( \nabla^2 f(x) + \subscr{k}{p}^{\textup{eq}}A^{\top} A + \sum_{j=1}^m \relu(y_j)\nabla^2 g_j(x) +\subscr{k}{p}^{\textup{in}}J_g(x)^{\top} G(y) J_g(x) \right) & M^{-1}A^{\top} & M^{-1}J_g(x)^{\top} G(y) \\
-\subscr{k}{i}^{\textup{eq}}A & 0 & 0 \\
-\subscr{k}{i}^{\textup{in}}\subscr{k}{p}^{\textup{in}}G(y)J_g(x) & 0 & \subscr{k}{i}^{\textup{in}}(I_m-G(y))
\end{bmatrix}.
$$
To prove globally weak contractivity, it suffices to show that $\mu_{P}(J(z)) \leq 0$ or equivalently $\frac{1}{2}\left(P J(z) + J(z)^{\top} P\right) \preceq 0$ for all $z$ for which $J(z)$ exists.
A direct calculation gives
\begin{align*}
\frac{1}{2}\left(P J(z) + J(z)^{\top} P\right) =
\begin{bmatrix} 
-\nabla^2 f(x) - \subscr{k}{p}^{\textup{eq}} A^{\top} A - \sum_{j=1}^m\relu(y_j)\nabla^2g_j(x) - \subscr{k}{p}^{\textup{in}}J_g(x)^{\top} G(y)J_g(x) & 0 & 0 \\
0 & 0 & 0 \\
0 & 0 & -\frac{1}{\subscr{k}{p}^{\textup{in}}}(I_m - G(y))
\end{bmatrix},
\end{align*}
which is block-diagonal. The LMI then follows by noticing that convexity of $f$ and $g_j$ implies $\nabla^2 f(x) \succeq 0$ and $\nabla^2 g_j(x) \succeq 0$, for all $x$, and $0 \preceq G(y) \preceq I_m$, for all $y$.

Next, we prove local strong contraction around $z^\star$. By strict complementarity, for $j \in \mathcal{A}$, $g_j(x^\star) = 0$ and $\xi^\star_j > 0$ imply $y^\star_j > 0$ and thus $G_j(y^\star) = 1$. For $j \notin \mathcal{A}$, $g_j(x^\star) < 0$ and $\xi^\star_j = 0$ imply $y^\star_j < 0$, hence $G_j(y^\star) = 0$.
In particular, $y^\star_j \neq 0$ for every $j$ and, by continuity, the activation pattern remains constant in a neighborhood of $z^\star$. The inactive multiplier dynamics satisfy $\dot{\xi}_{\mathcal{I}} = -\subscr{k}{i}^{\textup{in}} \xi_{\mathcal{I}}$, and therefore converge exponentially. The remaining active subsystem $(x, \nu, \xi_{\mathcal{A}})$ has Jacobian
$$
J_{\mathcal{A}}(z^\star)=
\begin{bmatrix}
-M^{-1}\left( \nabla^2 f(x^\star) + \subscr{k}{p}^{\textup{eq}}A^{\top} A + \sum_{j=1}^m \xi^\star_j\nabla^2 g_j(x^\star) +\subscr{k}{p}^{\textup{in}}J_g^{\mathcal{A}}(x^\star)^{\top} J_g^{\mathcal{A}}(x^\star) \right) & -M^{-1}A^{\top} & -M^{-1}J_g^{\mathcal{A}}(x^\star)^{\top} \\
\subscr{k}{i}^{\textup{eq}}A & 0 & 0 \\
\subscr{k}{i}^{\text{in}}\subscr{k}{p}^{\textup{in}}J_g^{\mathcal{A}}(x^\star) & 0 & 0
\end{bmatrix}.
$$
By LICQ, the matrix $C := \begin{bmatrix} A\\ J_g^{\mathcal A}(x^\star) \end{bmatrix}$ has full row rank. By SOSC at $(x^\star, \xi^\star)$, we have $H := \nabla^2 f(x^\star) + \sum_{j=1}^m \xi^\star_j\nabla^2 g_j(x^\star) \succ 0$. 
In fact, SOSC implies
$d^\top H d > 0$ for all $d \neq \0_n$ satisfying $Ad=0$, and $J_g^{\mathcal A}(x^\star)d=0$, that is, for all nonzero $d\in\ker(C) \setminus \{\0_n\}$. For any $d\neq \0_n$,
$
\begin{aligned}
d^\top J_1(z^\star)d &= d^\top H d + \subscr{k}{p}^{\textup{eq}}\|Ad\|^2 + \subscr{k}{p}^{\textup{in}} \|J_g^{\mathcal A}(x^\star)d\|^2.
\end{aligned}
$
If $d\in\ker(C)$, then the penalty terms vanish and the SOSC implies $d^\top J_1(z^\star)d = d^\top H d >0$.
Otherwise, $Cd\neq0$, and we have
$
\subscr{k}{p}^{\textup{eq}}\|Ad\|^2 + \subscr{k}{p}^{\textup{in}} \|J_g^{\mathcal A}(x^\star)d\|^2 >0.
$
Moreover, $H\succeq0$ by convexity and $\xi_j^\star\ge0$. Hence,
$
d^\top J_1(z^\star)d>0
$
for every $d\neq0$, and therefore
$
J_1(z^\star)\succ0 .
$
Then, $J_{\mathcal{A}}(z^\star)$ has the saddle form in Lemma~\ref{lemma:saddle-matrices_general}, which in turn guarantees that the $z^\star$ is locally exponentially stable with rate 
$$
c_{\mathcal{A}}= \frac{1}{4} \min \left\{ \frac{1}{\lambda_{\max}(J_1(z^\star))}, \tau \frac{\lambda_{\min}(J_1(z^\star))} {\sigma_{\max}^2(C)} \right\} \frac{\sigma_{\min}^2(C)}{1 + \subscr{k}{d}^{\textup{eq}} \subscr{a}{max}},
$$

where $\tau := \min \bigl(\subscr{k}{i}^{\textup{eq}}, \subscr{k}{i}^{\textup{in}}\subscr{k}{p}^{\textup{in}}\bigr)^{-1}$.
Then,~\eqref{eq:sp_pid_eq_ineq} is strongly infinitesimally contracting in a neighborhood of $z^\star$~\cite{TS:75}. The statement then follows by~\cite[Theorem 2]{VC-AD-AG-GR-FB:24a}. This concludes the proof.
\end{proof}
\end{arxiv}%
The proof of Theorem~\ref{thm:lin_exp_eq_ineq} reveals how the problem geometry determines the local exponential convergence rate. Specifically, the inactive constraint dynamics satisfy $\dot{\xi}_{\mathcal{I}} = -\subscr{k}{i}^{\textup{in}} \xi_{\mathcal{I}}$,  yielding the exponential rate $c_{\mathcal{I}} = \subscr{k}{i}^{\textup{in}}$. For the active subsystem $(x, \nu, \xi_{\mathcal{A}})$, the local active rate is $c_{\mathcal{A}}$. A local exponential convergence rate is then $\subscr{c}{loc} = \min(c_{\mathcal{I}}, c_{\mathcal{A}})$.

\subsection{Relation to Primal--Dual Gradient Dynamics}
\label{sec:comparison}
The Augmented Primal--Dual Gradient Dynamics (Aug-PDGD) proposed in~\cite{YT-GQ-NL:20} for convex nonlinear programs coincide with SPPID~\eqref{eq:sp_pid_eq_ineq} when $\subscr{k}{p}^{\textup{eq}} = \subscr{k}{d}^{\textup{eq}} = 0$, $\subscr{k}{i}^{\textup{eq}} = \subscr{k}{i}^{\textup{in}} = 1$, and the penalty parameter $\rho$ in~\cite{YT-GQ-NL:20} equals $\subscr{k}{p}^{\textup{in}}$. 
From our control-theoretic perspective, the Aug-PDGD framework employs only an integral control to the equality Lagrange multiplier. Extending the observations in~\cite{VC-SMF-SP-DR:25} for equality-constrained settings, we show that the absence of proportional feedback in the equality channel can preclude local stability in nonconvex problems.
\begin{arxiv}
\begin{mylemma}[label=lem:comparison]{SPPID local stability and Aug-PDGD instability}
There exist nonconvex programs for which the Aug-PDGD dynamics of~\cite{YT-GQ-NL:20} diverge for every choice of parameters, whereas the SPPID dynamics~\eqref{eq:sp_pid_eq_ineq_convex} are locally asymptotically stable for suitable gain selections.
\end{mylemma}
\end{arxiv}%
\begin{proof}
Consider the nonlinear program
\begin{arxiv}
\beq
\label{eq:comparison_problem}
\min_{x\in\R^2} \frac{w_1}{2}x_1^2 + \frac{w_2}{2}x_2^2 \quad \quad \textup{s.t. } a x_2 = 0, \; x_1 - d \leq 0,
\eeq
\end{arxiv}
with $w_1 > 0$, $w_2 < 0$, $a \neq 0$, and $d > 0$. The KKT conditions have the unique solution $z^\star = (\0_2, 0, 0)$.
We show that $z^\star$ is unstable under Aug-PDGD, while it is locally asymptotically stable under SPPID for any $\subscr{k}{i}^{\textup{eq}}, \subscr{k}{i}^{\textup{in}}, \subscr{k}{p}^{\textup{in}} > 0$ and $\subscr{k}{p}^{\textup{eq}} > \frac{|w_2|}{a^2}$.
Since $g(x^\star) = -d < 0$, there exists a neighborhood of $x^\star$ where the inequality constraint remains inactive. Thus, $\relu(\xi + \subscr{k}{p}^{\textup{in}}g(x)) = 0$ locally, and the $\xi$-subsystem reduces to the exponentially stable system $\dot{\xi} = -\subscr{k}{i}^{\textup{in}}\xi$. Moreover, $\dot{x}_1 = -w_1 x_1$ in the same local set, and this dynamics is exponentially stable since $w_1>0$. Local stability is then characterized by the $(x_2,\nu)$ subsystem.
For Aug-PDGD, the $(x_2,\nu)$ dynamics around $x^\star$ are $\dot{z} = \subscr{A}{aug}z$, where $z = [x_2, \nu]^{\top}$ and $\subscr{A}{aug} = \begin{bmatrix} -w_2 & - a \\ a  & 0 \end{bmatrix}$. Since the determinant of $\subscr{A}{aug}$ is positive and its trace is negative, both eigenvalues of $\subscr{A}{aug}$ have positive real part, and thus $x^\star$ is unstable under Aug-PDGD for every parameter choice.
For SPPID the linearized subsystem is $\dot{z} = \subscr{A}{PI}z$ with $\subscr{A}{PI} = \begin{bmatrix} -\gamma(w_2 + a^2\subscr{k}{p}^{\textup{eq}}) & -\gamma a \\ a\subscr{k}{i}^{\textup{eq}} & 0 \end{bmatrix}$, where $\gamma := \bigl(1+\subscr{k}{d}^{\textup{eq}}a^2\bigr)^{-1}>0$.
For any $\subscr{k}{i}^{\textup{eq}} > 0$ and $\subscr{k}{p}^{\textup{eq}} > |w_2|/a^2$, the trace $-\gamma(w_2+a^2\subscr{k}{p}^{\textup{eq}})$ and determinant $\gamma a^2 \subscr{k}{i}^{\textup{eq}}$ of $\subscr{A}{PI}$ are strictly negative and positive, respectively. Local asymptotic stability then follows by the Routh--Hurwitz criterion.
\begin{arxiv}
We refer to Figure~\ref{fig:convergence_comparison_} for numerical illustrations of this behavior.~
\end{arxiv}%
This concludes the proof.
\end{proof}
\begin{arxiv}
\begin{figure}[!h]
    \centering
    \includegraphics[width=.73\linewidth]{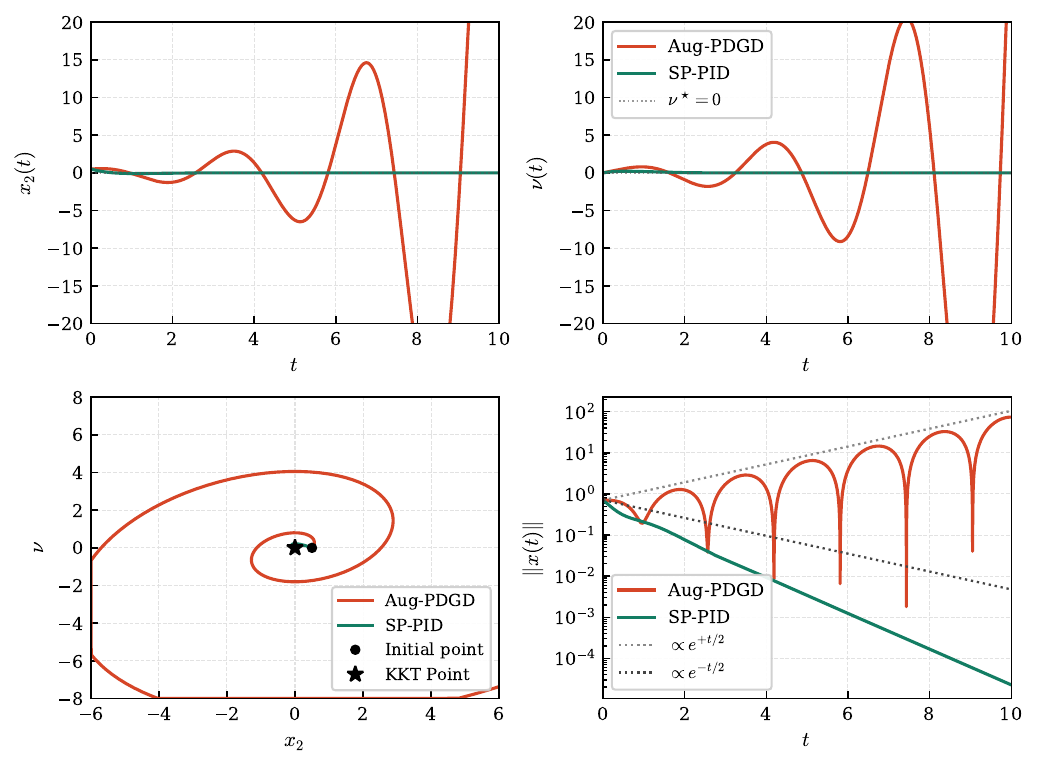}
    \caption{Numerical comparison of Aug-PDGD and the SPPID under the non-convex quadratic program for $w_1 = 1$, $w_2 = -1$, $a = 2$, $d = 1$, $\subscr{k}{p}^{\textup{eq}} = 1$, $\subscr{k}{i}^{\textup{eq}} = 1$, $\subscr{k}{d}^{\textup{eq}} = 0$, $\subscr{k}{p}^{\textup{in}} = 1$, $\subscr{k}{i}^{\textup{in}} = 1$. (a) Time trajectory of the equality state variable $x_2(t)$ and (b) of the dual multiplier $\nu(t)$, demonstrating the exponential divergence of Aug-PDGD against the asymptotic stabilization achieved by the SPPID. (c) Phase portrait in the $(x_2, \nu)$-plane reflecting mirror-symmetric closed-loop eigenvalues across the imaginary axis. (d) Transient evolution of the primal norm $\|x(t)\|$ relative to analytical exponential bounds.}
    \label{fig:convergence_comparison_}
\end{figure}
\end{arxiv}%
\section{Equality-Constrained Optimization Problems}
\label{sec:equality_constraint}
We specialize the general framework of Section~\ref{sec:sp_pid} to equality-constrained optimization problems
\beq
\begin{aligned}
    \min_{x \in \R^n} & f(x)\\
    \text{s.t. } & h(x) = \0_p,
\end{aligned}
\label{eq:eq_constrained}
\eeq
where $\map{f}{\R^n}{\R}$ and $\map{h}{\R^n}{\R^p}$ are continuously differentiable. 
\begin{arxiv}
We assume that~\eqref{eq:eq_constrained} admits at least one optimal solution.~
\end{arxiv}
By omitting the inequality dynamics, the corresponding SPPID closed-loop dynamics reduce from~\eqref{eq:sp_pid_eq_ineq} to
\begin{arxiv}
\beq
\label{eq:sp_pid_eq}
\begin{cases}
\dot{x}(t)= -\nabla f(x) - J_h(x)^\top \left(\nu + \subscr{k}{p}^{\textup{eq}} h(x) + \subscr{k}{d}^{\textup{eq}} J_h(x)\dot{x}\right) \\
\dot{\nu} = \subscr{k}{i}^{\textup{eq}} h(x),
\end{cases}
\iff
\begin{cases}
\dot{x} = - M(x)^{-1}\bigl(\nabla f(x) + J_h(x)^\top \nu + \subscr{k}{p}^{\textup{eq}} J_h(x)^\top h(x)\bigr),\\
\dot{\nu} = \subscr{k}{i}^{\textup{eq}} h(x),
\end{cases}
\eeq
\end{arxiv}%
where $\nu\in\R^p$ is the internal (integral) state, and $M(x):= I_n + \subscr{k}{d}^{\textup{eq}} J_h(x)^\top J_h(x)$. 
As in the general case, the proportional action admits an equivalent augmented-Lagrangian interpretation. Specifically, given $\subscr{k}{p}^{\textup{eq}}\ge 0$, consider the augmented Lagrangian associated with~\eqref{eq:eq_constrained}, that is 
\beq
\label{eq:augmented_lag_eq}
\subscr{L}{aug}^{\textup{eq}}(x,\nu) := f(x) + \nu^\top h(x) + \frac{\subscr{k}{p}^{\textup{eq}}}{2} \norm{h(x)}^2.
\eeq
This setting inherits the structural properties of Theorem~\ref{thm:sp_pid_eq_ineq}. 
\begin{arxiv}
\begin{mycorollary}
[label=thm:unified_pid_framework]{Equilibria and saddle-flow structure of equality-constrained SPPID}
Consider problem~\eqref{eq:eq_constrained}, where $f \in \mathcal{C}^1$ and $h \in \mathcal{C}^2$, and SPPID~\eqref{eq:sp_pid_eq}. Given $\subscr{k}{i}^{\textup{eq}} > 0$ and $\subscr{k}{p}^{\textup{eq}}, \subscr{k}{d}^{\textup{eq}} \geq 0$, the following statements hold
\begin{enumerate}[label=\textup{(\roman*)}]
\item \label{thm_unified:item1}
A point $(x^\star,\nu^\star)$ is an equilibrium of~\eqref{eq:sp_pid_eq} if and only if there exists $\lambda^\star \in \R^p$ such that $(x^\star,\lambda^\star)$ satisfies the KKT conditions of~\eqref{eq:eq_constrained}, and $\lambda^\star=\nu^\star$.
\item \label{thm_unified:item3}
The closed-loop dynamics~\eqref{eq:sp_pid_eq} are a preconditioned primal--dual gradient flow of the augmented Lagrangian~\eqref{eq:augmented_lag_eq}, namely
\beq
\label{eq:riemannian_saddle_explicit}
\begin{cases}
\dot x = - M(x)^{-1} \nabla_x \subscr{L}{aug}^{\textup{eq}}(x,\nu), \\
\dot\nu = \subscr{k}{i}^{\textup{eq}} \nabla_\nu \subscr{L}{aug}^{\textup{eq}}(x,\nu).
\end{cases}
\eeq 
Consequently, the primal dynamics correspond to gradient descent of $\subscr{L}{aug}^{\textup{eq}}$ with respect to the Riemannian metric induced by $M(x)$, while the dual dynamics correspond to Euclidean gradient ascent scaled by the factor $\subscr{k}{i}^{\textup{eq}}$.
\end{enumerate}
\end{mycorollary}
\begin{proof}
Let $z:=(x,\nu)$, and let $z^\star:=(x^\star,\nu^\star)$ be an equilibrium of~\eqref{eq:sp_pid_eq}. The condition $\dot\nu=\0_p$ gives $h(x^\star)=\0_p$, and since $M(x^\star)\succ0$, the condition $\dot x=\0_n$ gives $\nabla f(x^\star) + J_h(x^\star)^\top \nu^\star = \0_n$.
This in turn implies, $\nu^\star := \lambda^\star - \subscr{k}{p}^{\textup{eq}} h(x^\star) - \subscr{k}{d}^{\textup{eq}} J_h(x^\star)\dot{x}^\star = \lambda^\star$. Those two conditions are exactly the KKT stationarity and primal feasibility conditions~\eqref{kkt_stationarity},~\eqref{kkt_stationarity_eq} for~\eqref{eq:eq_constrained}. Conversely, if $(x^\star,\lambda^\star)$ satisfies the KKT conditions of~\eqref{eq:eq_constrained}, then $z^\star:=(x^\star,\nu^\star)$ with $\nu^\star=\lambda^\star$ annihilates both vector fields in~\eqref{eq:sp_pid_eq}, so $z^\star$ is an equilibrium. This proves item~\ref{thm_unified:item1}.
Next, we show item~\ref{thm_unified:item3}. Note that $\nabla_x \subscr{L}{aug}^{\textup{eq}}(x,\nu) = \nabla f(x) + J_h(x)^\top \nu + \subscr{k}{p}^{\textup{eq}} J_h(x)^\top h(x)$ and $\nabla_{\nu} \subscr{L}{aug}^{\textup{eq}}(x,\nu) = h(x)$.
Moreover, $M(x)\succ 0$ implies that the equation $M(x)\dot x = -\nabla_x \subscr{L}{aug}^{\textup{eq}}(x,\nu)$ admits the unique solution $\dot x = - M(x)^{-1}\nabla_x \subscr{L}{aug}^{\textup{eq}}(x,\nu).$ The dual equation follows directly from $\dot \nu = \subscr{k}{i}^{\textup{eq}} h(x) = \subscr{k}{i}^{\textup{eq}} \nabla_\nu \subscr{L}{aug}^{\textup{eq}}(x,\nu)$. This proves that the closed-loop dynamics~\eqref{eq:sp_pid_eq} coincides with the saddle-point flow~\eqref{eq:riemannian_saddle_explicit}.
\end{proof}

The following lemma clarifies the relationship between The saddle-point PI (SPPI) dynamics for equality constraints (that is,~\eqref{eq:sp_pid_eq} with $\subscr{k}{d}^{\textup{eq}}=0$) and the PI-CMO dynamics~\cite{VC-SMF-SP-DR:25}.
\begin{mylemma}
[label=thm_unified:item2]{Global diffeomorphism}
Consider the equality-constrained problem~\eqref{eq:eq_constrained}, where $f \in \mathcal{C}^1$ and $h \in \mathcal{C}^1$. Given $\subscr{k}{p}^{\textup{eq}} \geq0$ the coordinate transformation
\beq
\label{eq:diffeo}
\map{T}{(x,\lambda)}{(x,\nu) := (x,\lambda - \subscr{k}{p}^{\textup{eq}} h(x))}
\eeq
is a smooth global diffeomorphism.
\end{mylemma}
\begin{proof}
Since $h$ is smooth, the map $T$ is smooth. Its inverse is explicitly given by $T^{-1}(x,\nu) = \bigl(x,\nu + \subscr{k}{p}^{\textup{eq}} h(x)\bigr)$, which is also smooth.
Moreover, its Jacobian is given by
\[
J_T(x,\lambda) =
\begin{bmatrix}
I_n & 0 \\
- \subscr{k}{p}^{\textup{eq}} J_h(x) & I_p
\end{bmatrix},
\]
which is block triangular with determinant equal to $1$. Hence $J_T(x,\lambda)$ is nonsingular for all $(x,\lambda) \in \R^{n+p}$. Therefore $T$ is a smooth bijection with smooth inverse, and thus a global diffeomorphism.
\end{proof}
Lemma~\ref{thm_unified:item2} establishes a smooth equivalence between the SPPI and the PI-CMO dynamics in~\cite{VC-SMF-SP-DR:25}. In particular, the transformation $T$ defines a relation between the two dynamical systems, preserving trajectories in the primal variable while inducing a smooth reparameterization of the dual variables.
\end{arxiv}%

\begin{rem}[Comparison with~\cite{JR-SLJ:26}]
The recent work~\cite{JR-SLJ:26} adopts an optimization perspective and considers a discrete-time PI control, establishing connections with augmented Lagrangian methods. In contrast, our analysis is formulated in continuous time and explicitly characterizes the induced dynamics as a saddle-point flow within a control-theoretic framework. Furthermore, the inclusion of the derivative term generalizes these results and reveals a fundamentally different structure: the derivative action induces a state-dependent metric, yielding a Riemannian saddle-point flow.
\end{rem}

The following lemma focuses on the geometric role of $\subscr{k}{d}^{\textup{eq}}$.
\begin{arxiv}
\begin{mylemma}[label=lem:kd_limit]{SPPID as projected gradient flow}
Consider the dynamics~\eqref{eq:sp_pid_eq}, with $h \in \mathcal{C}^2$. Let
$$
\Pi(x) := I_n - J_h(x)^{\top} \bigl(J_h(x)J_h(x)^{\top}\bigr)^{-1}J_h(x)
$$
{denote the orthogonal projector onto $\ker J_h(x)$, defined wherever $J_h(x)$ has full row rank,} let $F$ denote the vector field of SPPID~\eqref{eq:sp_pid_eq}, and define
$
F_\infty(x,\nu) :=
\begin{bmatrix}
-\Pi(x)\nabla f(x) \\
\subscr{k}{i}^{\textup{eq}} h(x)
\end{bmatrix}.
$
Then
\begin{enumerate}[label=\textup{(\roman*)}]
\item \label{lem:kd:item1}
For every $(x,\nu) \in \R^{n+p}$ at which $J_h(x)$ has full row rank, $F(x,\nu) \to F_\infty(x,\nu)$ as $\subscr{k}{d}^{\textup{eq}} \to \infty$. Moreover
$$
\norm{M(x)^{-1} - \Pi(x)} \leq \frac{\Sigma_{\max}} {\subscr{k}{d}^{\textup{eq}}\,\Sigma_{\min}^2},
$$
where $\Sigma_{\max}$ and $\Sigma_{\min}$ are the maximum and minimum eigenvalues of $J_h(x)J_h(x)^{\top}$, respectively.
\end{enumerate}
{Assume, in addition, that LICQ holds on an open neighborhood of $\mathcal{M} := \setdef{x \in \R^n}{h(x) = \0_p}$. Then:}
\begin{enumerate}[label=\textup{(\roman*)}, start=2]
\item \label{lem:kd:item2}
The manifold $\mathcal{M}$ is forward invariant under $F_\infty$.
\item \label{lem:kd:item3}
A point $(x^{\star},\nu^{\star}) \in \mathcal{M} \times \R^p$ is an equilibrium of $F_\infty$ if and only if $x^{\star}$ satisfies the KKT conditions for~\eqref{eq:eq_constrained}, for an arbitrary $\nu^{\star} \in \R^p$. Moreover, the primal and dual dynamics decouple on $\mathcal{M}$, and along any such equilibrium the associated Lagrange multiplier is uniquely given by 
$$
\lambda^{\star} = -\bigl(J_h(x^{\star})J_h(x^{\star})^{\top}\bigr)^{-1} J_h(x^{\star})\nabla f(x^{\star}).
$$
\end{enumerate}
\end{mylemma}
\begin{proof}
Applying the Woodbury formula~\cite{WWH-89} to $M(x)$ yields
$$
M(x)^{-1} = I_n - J_h(x)^{\top} \Bigl(\frac{1}{\subscr{k}{d}^{\textup{eq}}}I_p + \Sigma(x) \Bigr)^{-1}J_h(x),
$$
where $\Sigma(x) := J_h(x)J_h(x)^{\top}\succ 0$, since $J_h(x)$ has full row rank by hypothesis. As $\subscr{k}{d}^{\textup{eq}}\to\infty$, $M(x)^{-1} \to \Pi(x)$.
Since $\Pi(x)J_{h}(x)^{\top} = 0$, the first component of $F(x,\nu)$ satisfies
$$
-M(x)^{-1}\bigl(\nabla f(x) + J_{h}(x)^{\top}(\nu + \subscr{k}{p}^{\textup{eq}}h(x))\bigr) \to -\Pi(x)\nabla f(x).
$$
The second component of $F$ is independent of $\subscr{k}{d}^{\textup{eq}}$, so $F(x,\nu)\to F_{\infty}(x,\nu)$ as $\subscr{k}{d}^{\textup{eq}} \to\infty$.
Next, we compute
\begin{align}
\|M(x)^{-1} - \Pi(x)\| &= \Big\|J_h(x)^{\top}\Bigl(\Sigma(x)^{-1} - \bigl(\frac{1}{\subscr{k}{d}^{\textup{eq}}}I_p + \Sigma(x)\bigr)^{-1}\Bigr) J_h(x)\Big\| \nonumber\\
&= \Big\|\frac{1}{\subscr{k}{d}^{\textup{eq}}}J_h(x)^{\top} \Sigma(x)^{-1}\Bigl(\frac{1}{\subscr{k}{d}^{\textup{eq}}} I_p+\Sigma(x)\Bigr)^{-1}J_h(x)\Big\| \label{eq:identity_inverse}\\
&\leq \frac{\Sigma_{\max}} {\subscr{k}{d}^{\textup{eq}}\,\Sigma_{\min}^2}, \label{ineq:inverse}
\end{align}
where in~\eqref{eq:identity_inverse} we applied the identity $A^{-1} - B^{-1} = A^{-1}(B-A)B^{-1}$ with $A = \Sigma$ and $B = \frac{1}{\subscr{k}{d}^{\textup{eq}}} I_p + \Sigma$, while in~\eqref{ineq:inverse} we used $\|(\frac{1}{\subscr{k}{d}^{\textup{eq}}}I_p+\Sigma)^{-1}\| \leq 1/\Sigma_{\min}$ and $\|J_h^\top\|\|J_h\| = \Sigma_{\max}$.

To show item~\ref{lem:kd:item2}, let $x(0) \in \mathcal{M}$. Along trajectories of $F_\infty$, we have
$
\frac{d}{dt}h(x(t)) = J_h(x)\dot{x} = -J_h(x)\Pi(x)\nabla f(x) = \0_p,
$
since $J_h(x)\Pi(x) = 0$. Hence $x(t) \in \mathcal{M}$ for all $t \geq 0$.

Finally, we prove item~\ref{lem:kd:item3}. A point $(x^{\star},\nu^{\star}) \in \mathcal{M} \times \R^p$ is an equilibrium of $F_\infty$ if and only if $\Pi(x^{\star})\nabla f(x^{\star}) = \0_n$ and $h(x^{\star}) = \0_p$. The second condition holds by definition of $\mathcal{M}$.
The first condition is equivalent to $\nabla f(x^{\star}) \in (\ker J_h(x^{\star}))^{\perp} = \operatorname{range}(J_h(x^{\star})^{\top})$, since $J_h(x^{\star})$ has full row rank. Hence there exists a unique $\lambda^{\star}$ satisfying the KKT stationarity condition~\eqref{kkt_stationarity}. Together with $h(x^{\star}) = \0_p$, this establishes the KKT conditions for~\eqref{eq:eq_constrained}. Finally, since $h(x^{\star}) = \0_p$, the dual equation $\subscr{k}{i}^{\textup{eq}}h(x^{\star}) = \0_p$ is satisfied for any $\nu^{\star} \in \R^p$, establishing the decoupling. This concludes the proof.
\end{proof}
\end{arxiv}%
Item~\ref{lem:kd:item3} reveals a structural difference between the finite-$\subscr{k}{d}^{\textup{eq}}$ and limiting regimes. For finite $\subscr{k}{d}^{\textup{eq}}$, stationarity of $\dot x$ together with LICQ forces $\nu^{\star} = \lambda^{\star}$ uniquely (Corollary~\ref{thm:unified_pid_framework}). In the limit $\subscr{k}{d}^{\textup{eq}} \to \infty$, the primal and dual dynamics decouple on $\mathcal{M}$, the equilibrium value of $\nu$ becomes arbitrary, and $\lambda$ is recovered directly from the primal optimality condition.

\begin{arxiv}
Table~\ref{tab:pid_special_cases} summarizes the saddle-point flows induced by SPPID~\eqref{eq:sp_pid_eq} for different choices of gains.
\begin{table}[!h]
\centering
\resizebox{.6\columnwidth}{!}{%
\begin{tabular}{c c c c}
\hline
Controller & $\subscr{k}{p}^{\textup{eq}}$ & $\subscr{k}{d}^{\textup{eq}}$ & Resulting Dynamics \\
\hline\\[-8pt]
I   & $0$ & $0$ & Arrow--Hurwicz--Uzawa flow \\[4pt]

PI  & $>0$ & $0$ & Augmented Lagrangian primal--dual \\[4pt]

PID & $\geq 0$ & $>0$ & Riemannian saddle-point flow \\[4pt]

PID & $\geq 0$ & $\to\infty$ & Projected gradient flow \\[2pt]
\hline
\end{tabular}
}
\caption{Saddle-point flows induced by SPPID~\eqref{eq:sp_pid_eq} for different gain choices. To ensure feasibility, we let $\subscr{k}{i}^{\textup{eq}} > 0$.}
\label{tab:pid_special_cases}
\end{table}
\end{arxiv}

\subsection{Convergence Analysis for Affine Equality Constraints}
\label{sec:convergence}
We study the convergence properties of SPPID in the setting of affine constraints, i.e., we consider $h(x) := Ax - b$, with $A \in \R^{p \times n}$, $b \in \R^p$. In this case, the Jacobian of $h$ is constant, which simplifies the geometry of the dynamics and allows us to establish global exponential convergence using contraction theory.
\begin{arxiv}
The equality-constrained optimization problem then becomes
\begin{align*}
\min_{x \in \R^n} & f(x)\\
\text{s.t. } & Ax = b.
\end{align*}
The SPPID~\eqref{eq:sp_pid_eq} becomes
\beq
\label{eq:sp_pid_eq_affine}
\begin{cases}
\dot x = {-} \left(I_n {+} \subscr{k}{d}^{\textup{eq}} A^\top A\right)^{-1} \left(\nabla f(x) {+} A^\top \nu {+} \subscr{k}{p}^{\textup{eq}} A^\top(Ax {-} b)\right), \\
\dot\nu = \subscr{k}{i}^{\textup{eq}} \left(Ax - b\right).
\end{cases}
\eeq
\end{arxiv}%
We work under the following assumptions.
\begin{arxiv}
\begin{myassumption}{Strong convexity and constraint structure}
\label{ass:1_eq}
For the dynamics~\eqref{eq:sp_pid_eq_affine}, assume:
\begin{enumerate}[label=\textup{($A$\arabic*)}, leftmargin=3em]
    \item \label{ass:1_eq_f} the function $\map{f}{\R^n}{\R}$ is $\rho$-strongly convex and $L$-smooth;
    \item \label{ass:2_A_eq} 
    the matrix $A \in \R^{p \times n}$ satisfies $\amin I_p \preceq AA^\top \preceq \amax I_p$ for $\amin,\amax \in \R_{>0}$.
\end{enumerate}
\end{myassumption}
\smallskip
\end{arxiv}%
Assumption~\ref{ass:1_eq} guarantees the existence of a unique global minimum and a unique KKT point of~\eqref{eq:eq_constrained}. These assumptions are standard in the literature when establishing global convergence (see, e.g.,~\cite{GQ-NL:19, AD-VC-AG-GR-FB:23f, VC-SMF-SP-DR:25}).
The next theorem shows that, under this Assumption, SPPID~\eqref{eq:sp_pid_eq} is globally contracting.
\begin{arxiv}
\begin{mytheorem}[label=thm:riemannian_saddle_affine-contractivity]{Contractivity of~\eqref{eq:sp_pid_eq_affine}}
Under Assumptions~\ref{ass:1_eq_f} and~\ref{ass:2_A_eq}, for any $\subscr{k}{p}^{\textup{eq}} \geq0 $, $\subscr{k}{i}^{\textup{eq}} > 0$, and $\subscr{k}{d}^{\textup{eq}} \geq 0$, the PID saddle-point flow~\eqref{eq:sp_pid_eq_affine} is strongly infinitesimally contracting with respect to $\|\cdot\|_{P}$ with rate $\subscr{c}{eq} > 0$, where
\begin{align}
P &=
\begin{bmatrix}
M & \alpha A^\top \\
\alpha A & \ds \frac{1}{\subscr{k}{i}^{\textup{eq}}} I_p
\end{bmatrix}
\succ 0,  \label{thm:def:P}\\
\alpha&=\frac{1}{2}\min\Bigg\{\frac{1}{L + \subscr{k}{p}^{\textup{eq}}\amax}, \frac{\rho}{\subscr{k}{i}^{\textup{eq}} \amax}\Bigg\}, \label{thm:def:alpha}\\
\subscr{c}{eq} &= \frac{1}{2}\alpha\subscr{k}{i}^{\textup{eq}}\frac{\amin}{1 + \subscr{k}{d}^{\textup{eq}}\amax}. \label{thm:def:c}
\end{align}
\end{mytheorem}
\begin{proof}
Let $z = (x,\nu) \in \R^{n + p}$, and $\dot{z} = \bar{\operatorname{F}}(z)$ denote~\eqref{eq:sp_pid_eq_affine}.
Assumption~\ref{ass:1_eq}.\ref{ass:1_eq_f} implies that $\nabla f$ is Lipschitz, $\nabla^2 f$ exists almost everywhere by Rademacher's theorem, and $\rho I_n \preceq \nabla^2 f(x) \preceq L I_n$, for all $x$ for which the Hessian exists. Hence the Jacobian of~\eqref{eq:sp_pid_eq_affine} exists almost everywhere and is given by
$$
\subscr{J}{F}(z)
:=
\begin{bmatrix}
- M^{-1} \left(\nabla^2 f(x) + \subscr{k}{p}^{\textup{eq}} A^\top A\right) & - M^{-1} A^\top\\
\subscr{k}{i}^{\textup{eq}} A & 0 
\end{bmatrix}.
$$
To prove strong infinitesimal contractivity it suffices to show that $\mu_{P}(\subscr{J}{F}(z)) \leq -\subscr{c}{eq}$ for all $z$ for which $\subscr{J}{F}(z)$ exists, with $P$ and $c$ as in~\eqref{thm:def:P},~\eqref{thm:def:c}  and~\eqref{thm:def:alpha}, respectively. The assumption of $\rho$-strong convexity and $L$-smoothness of $f$ imply the inequalities $\rho I_n \preceq \nabla^2 f(x) \preceq L I_n$, for all $x$ for which the Hessian exists. Moreover,
\begin{equation*}
\sup_{z}\mu_{P}(\subscr{J}{F}(z)) \leq \max_{\rho I_n \preceq B \preceq L I_n} \mu_{P}\left(
\begin{bmatrix}
- M^{-1} \left(B + \subscr{k}{p}^{\textup{eq}} A^\top A\right)  & - M^{-1} A^\top \\
\subscr{k}{i}^{\textup{eq}} A & 0
\end{bmatrix}
\right), 
\end{equation*}
where the $\sup$ is over all points where $\subscr{J}{F}(z)$ exists. 
The above matrix has the scaled saddle structure of Lemma~\ref{lemma:saddle-matrices_general} in Appendix~\ref{apx:log_norm_of_saddle_matrix} with $B' := B + \subscr{k}{p}^{\textup{eq}} A^\top A$, $\tau := \frac{1}{\subscr{k}{i}^{\textup{eq}}}$, and $M = I_n + \subscr{k}{d}^{\textup{eq}} A^\top A$.
The result then follows from Lemma~\ref{lemma:saddle-matrices_general} with $\subscr{m}{min}= 1$, $\subscr{m}{max}= 1 + \subscr{k}{d}^{\textup{eq}} \amax$, $\subscr{b}{min} = \rho$, $\subscr{b}{max} = L + \subscr{k}{p}^{\textup{eq}}\amax$, and $\Theta = \subscr{k}{i}^{\textup{eq}} I_p$ (i.e., $\subscr{\theta}{min}=\subscr{\theta}{max}=\subscr{k}{i}^{\textup{eq}}$).
\end{proof}%
Theorem~\ref{thm:riemannian_saddle_affine-contractivity} proves that~\eqref{eq:sp_pid_eq_affine} is contracting for any gains $\subscr{k}{i}^{\textup{eq}}>0$, $\subscr{k}{p}^{\textup{eq}}\geq 0$, and $\subscr{k}{d}^{\textup{eq}}\geq0$. In particular, no additional tuning conditions are required to guarantee global exponential convergence. The gains only affect the contraction rate $\subscr{c}{eq}$ in~\eqref{thm:def:c}, and therefore influence the speed of convergence but not stability. A direct consequence of Theorem~\ref{thm:riemannian_saddle_affine-contractivity} is that all trajectories of~\eqref{eq:sp_pid_eq_affine} converge exponentially to a unique equilibrium. Contractivity also ensures incremental stability and robustness with respect to perturbations of the vector field~\cite{FB:26-CTDS}.

The proven worst-case contraction rate $\subscr{c}{eq}$ in \eqref{thm:def:c} decreases monotonically with $\subscr{k}{d}^{\textup{eq}}$, so a larger derivative gain does not by itself certify faster convergence. This raises the question of what role $\subscr{k}{d}^{\textup{eq}}$ plays beyond the geometric interpretation of Lemma~\ref{lem:kd_limit}. We point out one such role, related to the numerical conditioning of flow, on a simple example. A general characterization is left to future work.

Let $f(x) = \frac{1}{2}\|x\|^2$ and let $b \in \R^p$ and $A\in \R^{p \times n}$ be diagonal with $0 < a_1 \le \dots \le a_p$. Then $M(x) = \operatorname{diag}(m_1,\dots,m_n)$, where $m_j := 1 + \subscr{k}{d}^{\textup{eq}} a_j^2$ for $j = 1,\dots, p$ and $m_j := 1$ for $j = p+1,\dots, n$. The unique equilibrium is $x_j^\star = b_j/a_j$, $\nu_j^\star = -x_j^\star/a_j$ for $j \le p$, and $x_j^\star = 0$ for $j > p$. Writing $e_{x,j} := x_j - x_j^\star$, $e_{\nu,j} := \nu_j - \nu_j^\star$, the error dynamics is
\begin{align}
m_j \dot e_{x,j} &= -\big(1 + \subscr{k}{p}^{\textup{eq}} a_j^2\big) e_{x,j} - a_j e_{\nu,j},
\label{eq:diag-ex}\\
\dot e_{\nu,j} &= \subscr{k}{i}^{\textup{eq}} a_j e_{x,j}. \label{eq:diag-enu}
\end{align}
Differentiating \eqref{eq:diag-ex} and substituting \eqref{eq:diag-enu} gives
$$
m_j \ddot e_{x,j} + \big(1 + \subscr{k}{p}^{\textup{eq}} a_j^2\big) \dot e_{x,j} + \subscr{k}{i}^{\textup{eq}} a_j^2 e_{x,j} = 0, \quad j = 1,\dots,n.
$$
For $j \le p$, this is a damped harmonic oscillator with frequency
$
\omega_j\big(\subscr{k}{d}^{\textup{eq}}\big) = a_j \sqrt{\dfrac{\subscr{k}{i}^{\textup{eq}}}{1 + \subscr{k}{d}^{\textup{eq}} a_j^2}},
$
while for $j > p$ the error evolves independently of $\subscr{k}{d}^{\textup{eq}}$ as $\dot e_{x,j} = -e_{x,j}$.
At $\subscr{k}{d}^{\textup{eq}} = 0$, the frequency spread across the $p$ constrained components is $\omega_p/\omega_1 = a_p/a_1 =: \kappa(A)$, the condition number of $A$, which grows unboundedly as $A$ becomes ill-conditioned. As $\subscr{k}{d}^{\textup{eq}} \to \infty$, however, $\omega_j(\subscr{k}{d}^{\textup{eq}}) \sim \sqrt{\subscr{k}{i}^{\textup{eq}}/\subscr{k}{d}^{\textup{eq}}}$ uniformly in $j$, so $\omega_p/\omega_1 \to 1$. This simplified (diagonal, quadratic) scenario makes explicit how $\subscr{k}{d}^{\textup{eq}}$ compresses the spread of time scales induced by an ill-conditioned $A$. Since numerical integration of a stiff continuous-time flow requires a step size inversely proportional to the fastest mode, this suggests that $\subscr{k}{d}^{\textup{eq}}$ may also mitigate stiffness in a discretized implementation of SPPID.

We verify this beyond the idealized diagonal setting in Figures~\ref{fig:role_of_kd1} and~\ref{fig:role_of_kd2}. For a dense, coupled quadratic program, Fig.~\ref{fig:role_of_kd1} confirms the same frequency-spread compression and a corresponding growth in the admissible explicit-Euler step size as $\subscr{k}{d}^{\textup{eq}}$ increases, and Fig.~\ref{fig:role_of_kd2} shows that this translates directly into convergence of a fixed-step Euler discretization of SPPID that would otherwise diverge for ill-conditioned $A$.
\begin{figure}[t]
  \centering
  \includegraphics[width=.9\linewidth]{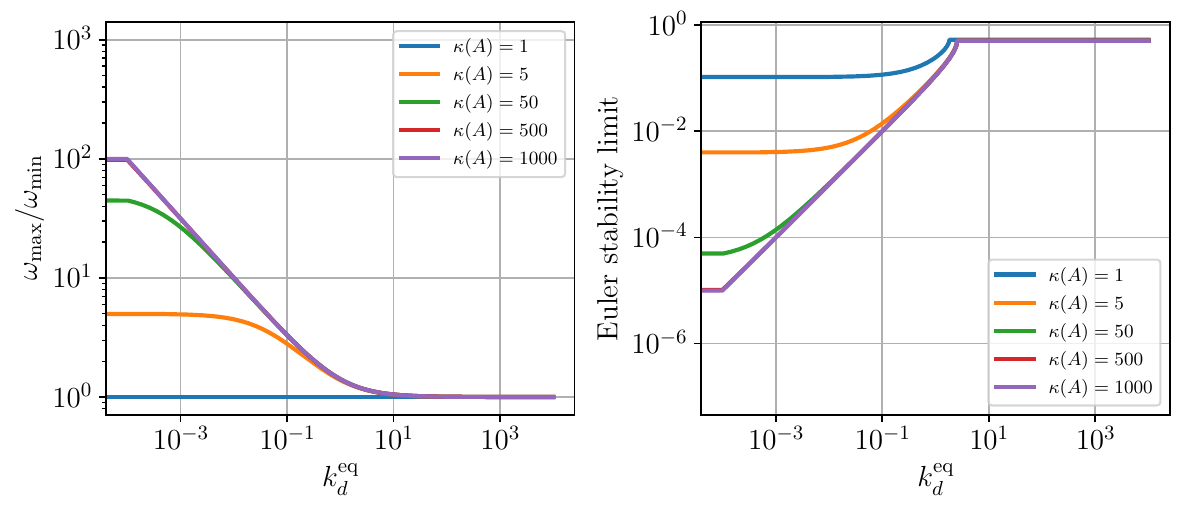}
  \caption{Numerical validation of the frequency- and stiffness-compression effects in a dense, coupled setting. For $\kappa(A) \in \{1,5,50,500,1000\}$, $A = U\Sigma V^\top \in \R^{20\times 40}$ is built from singular values geometrically spaced between $1$ and $\kappa(A)$, embedded in a single pair of fixed random orthogonal bases $U,V$ shared across all $\kappa$ so that only the conditioning of $A$ varies. We generalize the scalar example by considering a quadratic cost $f(x) =\tfrac{1}{2} x^\top Q x$ with $Q = I + 0.25 A^\top A$. Although $Q$ and $A^\top A$ commute by construction $A$ itself is dense in the standard coordinates used to integrate~\eqref{eq:sp_pid_eq_affine}, so this experiment tests the predicted phenomenon without the axis-aligned structure of the diagonal case.
  \emph{Left:} the closed-form ratio $\omega_{\max}/\omega_{\min}$ evaluated at the singular values of $A$ as a function of $\subscr{k}{d}^{\textup{eq}}$. Consistent with the analysis, the spread equals $\kappa(A)$ at $\subscr{k}{d}^{\textup{eq}}=0$ and collapses toward unity as $\subscr{k}{d}^{\textup{eq}}\to\infty$, uniformly across all five condition numbers.
   \emph{Right:} the maximum step size admitted by explicit Euler integration computed from the eigenvalues of the linearized error Jacobian at the equilibrium. The stability limit grows by several orders of magnitude as $\subscr{k}{d}^{\textup{eq}}$ increases, showing that the stiffness-reducing benefit of derivative feedback predicted analytically for the decoupled case persists in this dense, coupled setting.}
\label{fig:role_of_kd1}
\end{figure}

\begin{figure}[t]
  \centering
  \includegraphics[width=.9\linewidth]{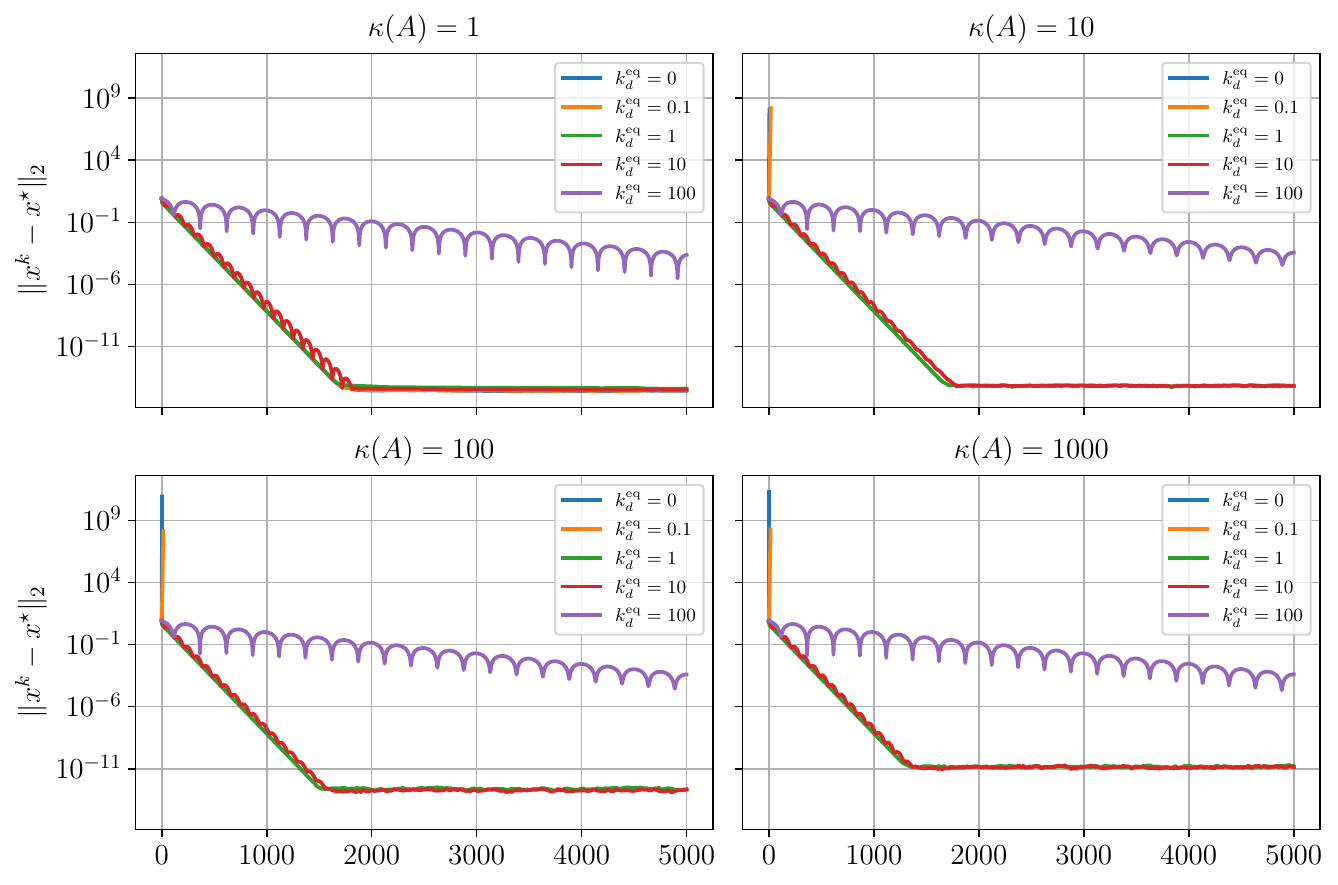}
  \caption{Effect of $\subscr{k}{d}^{\textup{eq}}$ on a fixed-step explicit Euler implementation of SPPID. Each panel corresponds to one condition number $\kappa(A)\in\{1,10,100,1000\}$, constructed as in Fig.~\ref{fig:role_of_kd1}, with right-hand side $b = A \subscr{x}{feas}$ chosen so that a consistent equilibrium $(x^\star,\nu^\star)$ exists, obtained by directly solving the KKT system. Curves show the primal error $\|x^k-x^\star\|_2$ over $5000$ Euler steps of fixed size $h = 0.02$, from a common initial condition $(x^0,\nu^0=0)$, for $\subscr{k}{d}^{\textup{eq}} \in \{0,\,10^{-1},\,1,\,10,\,10^2\}$. For well-conditioned $A$ ($\kappa=1$), all choices of $\subscr{k}{d}^{\textup{eq}}$ converge. As $\kappa(A)$ grows, trajectories with $\subscr{k}{d}^{\textup{eq}}=0$ (and small $\subscr{k}{d}^{\textup{eq}}$) diverge because $h$ exceeds the Euler stability limit of the fastest mode, while larger $\subscr{k}{d}^{\textup{eq}}$ keeps the same fixed-step iteration stable and convergent.}
  \label{fig:role_of_kd2}
\end{figure}
\end{arxiv}

\subsection{Convergence Analysis for Non-Affine Constraints}
Inspired by the proof of Theorem~\ref{thm:riemannian_saddle_affine-contractivity}, we now give a local stability result of SPPID under non-affine constraints. 
Specifically, we show that under suitable assumptions, for a sufficiently large proportional gain, SPPID recovers local exponential stability around a KKT point.

\begin{arxiv}
\begin{mytheorem}[label=thm:local_exp_stability]{Local exponential stability under SOSC and LICQ}
Consider the dynamics~\eqref{eq:sp_pid_eq} with $f, h \in \mathcal{C}^2$. 
Let $z^\star = (x^\star, \nu^\star)$ be a KKT point of~\eqref{eq:eq_constrained} at which LICQ holds. Assume SOSC holds, that is, $v^\top \nabla^2_x L(z^\star) v > 0$ for all $v \in \ker J_h(x^\star) \setminus \{\0_n\}$, where $L(x, \nu) := f(x) + \nu^\top h(x)$. Moreover, let
\begin{align*}
&\subscr{\lambda}{sosc} := \min_{\substack{v \in \ker J_h(x^\star) \\ \|v\| = 1}} v^\top \nabla^2_x L(z^\star) v > 0, \text{ and }\\
&\bar{k} := \max \left(0, \frac{2\| \nabla^2_x L(z^\star)\|^2 / \subscr{\lambda}{sosc}  - \lambda_{\min}( \nabla^2_x L(z^\star))} {\sigma_{\min}^2(J_h(x^\star))}\right).
\end{align*}
Then, for every $\subscr{k}{i}^{\textup{eq}} > 0$, $\subscr{k}{d}^{\textup{eq}} \geq 0$, and $\subscr{k}{p}^{\textup{eq}} > \bar{k}$, the equilibrium $z^\star$ is locally exponentially stable for SPPID~\eqref{eq:sp_pid_eq}.
\end{mytheorem}
\end{arxiv}
\begin{proof}
Let $z = (x,\nu)$ and let $F$ denote the right-hand side of~\eqref{eq:sp_pid_eq}, and write the primal component as $F_1(z) := - M(x)^{-1} q(z)$, where $q(z) := \nabla f(x) + J_h(x)^\top\bigl(\nu + \subscr{k}{p}^{\textup{eq}} h(x)\bigr)$. The Jacobian of $F_1$ with respect to $x$ at $z^{\star}$ satisfies
$$
\frac{\partial F_1}{\partial x}\bigg|_{z^\star} = -\frac{\partial M(x)^{-1}}{\partial x}\bigg|_{z^\star} q(z^\star) - M(x^\star)^{-1} \frac{\partial q}{\partial x}\bigg|_{z^\star}.
$$
By KKT stationarity, $q(z^\star) = \0_n$, so the first term vanishes. For the second term, using $h(x^\star) = \0_p$ and the chain rule
\begin{align*}
\frac{\partial q}{\partial x}\bigg|_{z^\star} &= \nabla^2 f(x^\star) + \sum_{i=1}^p \nu^\star_i \nabla^2 h_i(x^\star) + \subscr{k}{p}^{\textup{eq}} J_h(x^\star)^\top J_h(x^\star) \\
&= \nabla^2_x L(z^\star) + \subscr{k}{p}^{\textup{eq}} J_h(x^\star)^\top J_h(x^\star) =: B^\star.
\end{align*}
Next, we show $B^\star \succ 0$.
Let $v \in \R^n \setminus \{\0_n\}$ and decompose $v = v_0 + v_1$ with $v_0 \in \ker J_h(x^\star)$ and $v_1 \in \mathrm{range}(J_h(x^\star)^\top)$. Since $J_h(x^\star) v_0 = \0_p$, we have $\|J_h(x^\star) v\|^2 = \|J_h(x^\star)v_1\|^2$. By SOSC, $v_0^\top \nabla^2_x L(z^\star) v_0 \geq \subscr{\lambda}{sosc}\|v_0\|^2 > 0$. Applying Cauchy–Schwarz and Young's inequality with parameter $\gamma = \subscr{\lambda}{sosc}/2$ yields
\begin{align*}
2 v_0^\top \nabla^2_x L(z^\star) v_1 &\geq - 2 \norm{\nabla^2_x L(z^\star)} \norm{v_0} \norm{v_1}  \\
&\geq - \frac{\subscr{\lambda}{sosc}}{2}\|v_0\|^2 - \frac{2}{\subscr{\lambda}{sosc}}\|\nabla^2_x L(z^\star)\|^2\|v_1\|^2.
\end{align*}
By LICQ, $\|J_h(x^\star)v_1\|^2 \geq \sigma_{\min}^2(J_h(x^\star))\|v_1\|^2$. Therefore
\begin{arxiv}
\begin{align*}
v^\top B^\star v &\geq \frac{\subscr{\lambda}{sosc}}{2}\|v_0\|^2
+ \Bigl(\lambda_{\min}(\nabla^2_x L(z^\star)) + \subscr{k}{p}^{\textup{eq}}\sigma_{\min}^2(J_h(x^\star)) - \frac{2\|\nabla^2_x L(z^\star)\|^2}{\subscr{\lambda}{sosc}} \Bigr) \|v_1\|^2.
\end{align*}
\end{arxiv}
For $\subscr{k}{p}^{\textup{eq}} > \bar{k}$, both coefficients are strictly positive. Since $v_0$ and $v_1$ are perpendicular, $\|v_0\|^2 + \|v_1\|^2 = \|v\|^2 > 0$, so $v_0$ and $v_1$ cannot both vanish. Hence $v^\top B^\star v > 0$, and $B^\star \succ 0$.

Let $M^\star := M(x^\star)$. The Jacobian of SPPID at $z^\star$ is
$$
J_F(z^\star) = \begin{bmatrix}
- {M^\star}^{-1} B^\star & 
-{M^\star}^{-1} J_h(x^\star)^\top \\
\subscr{k}{i}^{\textup{eq}} J_h(x^\star) & \0
    \end{bmatrix},
$$
which has the scaled saddle structure of  Lemma~\ref{lemma:saddle-matrices_general} with matrices $(M^{\star}, B^{\star}, J_h(x^{\star}))$ and scaling $\Theta = \subscr{k}{i}^{\textup{eq}} I_p$.
Since $B^\star \succ 0$ and $M^\star \succ 0$, then $\mu_{P^\star}(J_F(z^\star)) \leq -c < 0$ for an explicit $ c > 0$ and $P^\star \succ 0$.
By continuity of $z \mapsto J_F(z)$ and of the log-norm, there exists a neighborhood $\mathcal{U}$ of $z^\star$ such that $\mu_{P^\star}(J_F(z)) \leq -c/2 < 0$ for all $z \in \mathcal{U}$. Hence SPPID is strongly infinitesimally contracting in $\mathcal{U}$~\cite{TS:75}, from which local exponential stability of $z^\star$ follows.
\end{proof}
\section{Inequality-Constrained Problems}
\label{sec:inequality_constraint}
We specialize the general framework of Section~\ref{sec:sp_pid} to inequality-constrained problems. 
\beq
\begin{aligned}
\min_{x \in \R^n} \quad & f(x)\\
\text{s.t.} \quad & g(x) \le \0_m,
\end{aligned}
\label{eq:ineq_constrained}
\eeq
where $\map{f}{\R^n}{\R}$ and $\map{g}{\R^n}{\R^m}$ are continuously differentiable.
In this scenario, the controller is an anti-windup PI controller acting on the inequality constraints. The corresponding closed-loop dynamics are
\beq
\label{eq:pi_spf_ineq}
\begin{cases}
\dot{x} = - \nabla f(x) - J_g(x)^\top \relu\big(\xi + \subscr{k}{p}^{\textup{in}} g(x)\big), \\
\dot{\xi} = - \subscr{k}{i}^{\textup{in}} \xi + \subscr{k}{i}^{\textup{in}} \relu\big(\xi + \subscr{k}{p}^{\textup{in}} g(x)\big),
\end{cases}
\eeq
which we refer to as the \emph{saddle-point PI dynamics (SPPI) for inequality-constrained problems}. The primal dynamics follows the negative gradient of $f$ augmented by a constraint-weighted penalty: this term vanishes when constraints are strictly inactive and acts as a restoring force when they are active or violated. The dual update drives $\xi$ toward the projected output $\mu$, penalizing constraint violation through the integral state.

Given $\subscr{k}{p}^{\textup{in}} > 0$, consider the augmented Lagrangian~\cite{DPB:97}
\begin{equation}
\label{eq:aug_lagrangian_ineq}
\subscr{L}{aug}^{\textup{in}}(x,\xi) = f(x) + \frac{1}{2\subscr{k}{p}^{\textup{in}}} \sum_{i=1}^m \left( \left[\relu\big(\xi_i + \subscr{k}{p}^{\textup{in}} g_i(x)\big) \right]^2 - \xi_i^2 \right).
\end{equation}

The next result follows directly from Theorem~\ref{thm:sp_pid_eq_ineq}.
\begin{arxiv}
\begin{mycorollary}
[label=thm:pi_ineq_framework]{SPPI closed-loop dynamics as saddle flows}
Consider the inequality-constrained problem~\eqref{eq:ineq_constrained}, where $\map{f}{\R^n}{\R}$ and $\map{g}{\R^n}{\R^m}$ are continuously differentiable. Given $\subscr{k}{p}^{\textup{in}} > 0$ and $\subscr{k}{i}^{\textup{in}} > 0$, consider the dynamics~\eqref{eq:pi_spf_ineq} and the augmented Lagrangian~\eqref{eq:aug_lagrangian_ineq}. Then:
\begin{enumerate}[label=\textup{(\roman*)}]
\item \label{thm_pi_ineq:item1}
A point $(x^\star,\xi^\star)$ is an equilibrium of~\eqref{eq:pi_spf_ineq} if and only if there exists $\mu^\star \in \R^m_{\geq 0}$ such that $(x^\star,\mu^\star)$ satisfies the KKT conditions of~\eqref{eq:ineq_constrained}, and $\mu^\star = \xi^\star$.
\item \label{thm_pi_ineq:item2}
The SPPI~\eqref{eq:pi_spf_ineq} coincides with the following primal--dual gradient flow of $\subscr{L}{aug}^{\textup{in}}$
\beq
\begin{cases}
\dot{x} = - \nabla_x \subscr{L}{aug}^{\textup{in}}(x, \xi), \\
\dot{\xi} = (\subscr{k}{i}^{\textup{in}} \subscr{k}{p}^{\textup{in}}) \nabla_\xi \subscr{L}{aug}^{\textup{in}}(x, \xi).
\end{cases}
\label{eq:pd_gradient_flow_ineq}
\eeq
\end{enumerate}
\end{mycorollary}
\begin{proof}
Let $(x^\star,\xi^\star) \in \R^{n+m}$ be an equilibrium of~\eqref{eq:pi_spf_ineq}. The condition $\dot{\xi} = \0_m$ implies $\xi^\star = \relu\big(\xi^\star + \subscr{k}{p}^{\textup{in}} g(x^\star)\big) = \mu^\star$, where the last equality follows by~\eqref{eq:controller_ineq}. By non-negativity of the $\relu$ operator, this fixed point condition automatically implies dual feasibility~\eqref{kkt_dual_feasibility}. We now show it is also equivalent to primal feasibility~\eqref{kkt_primal_feasibility} and complementary slackness~\eqref{kkt_complementarity}.
To see this, we analyze the condition component-wise for each $i \in \{1, \dots, m\}$. If $\mu_i^\star = 0$, then substituting into the fixed-point condition gives $0 = \relu\big(\subscr{k}{p}^{\textup{in}} g_i(x^\star)\big)$, which implies $g_i(x^\star) \le 0$. Thus,~\eqref{kkt_primal_feasibility} and~\eqref{kkt_complementarity} hold. If instead $\mu_i^\star > 0$, then $\mu_i^\star = \mu_i^\star + \subscr{k}{p}^{\textup{in}} g_i(x^\star)$, which implies $g_i(x^\star) = 0$. This also satisfies~\eqref{kkt_primal_feasibility} and~\eqref{kkt_complementarity}.
Finally, the equilibrium condition $\dot{x} = \0_n$ yields $\nabla f(x^\star) + J_g(x^\star)^\top \mu^\star = \0_n$. Thus, $(x^\star, \xi^\star)$ is a KKT point. Moreover $\xi^\star = \mu^\star$. Conversely, if $(x^\star,\mu^\star)$ satisfies the KKT conditions for~\eqref{eq:ineq_constrained}, then $\xi^\star := \mu^\star$ satisfies the fixed-point condition $\xi^\star = \relu(\xi^\star + \subscr{k}{p}^{\textup{in}} g(x^\star))$ by the same component-wise argument above, so $(x^\star,\xi^\star)$ annihilates both vector fields in~\eqref{eq:pi_spf_ineq} and is therefore an equilibrium. This concludes the proof of item~\ref{thm_pi_ineq:item1}.
Next, to prove~\ref{thm_pi_ineq:item2}, we note that $\subscr{L}{aug}^{\textup{in}}(x, \xi)$ in~\eqref{eq:aug_lagrangian_ineq} is continuously differentiable with respect to both arguments~\cite{DPB:97}. Recalling that the derivative of the squared projection function $\phi(s) = \frac{1}{2} [\relu(s)]^2$ is $\phi'(s) = \relu(s)$, applying the chain rule directly yields
\begin{align*}
    \nabla_x \subscr{L}{aug}^{\textup{in}}(x, \xi) &= \nabla f(x) + J_g(x)^\top \relu\big(\xi + \subscr{k}{p}^{\textup{in}} g(x)\big), \\
    \nabla_\xi \subscr{L}{aug}^{\textup{in}}(x, \xi) &= \frac{1}{\subscr{k}{p}^{\textup{in}}} \left(\relu\big(\xi + \subscr{k}{p}^{\textup{in}} g(x)\big) - \xi \right).
\end{align*}
Substituting these gradients back into~\eqref{eq:pd_gradient_flow_ineq} recovers the SPPI dynamics~\eqref{eq:pi_spf_ineq}. This concludes the proof.
\end{proof}
\end{arxiv}
While the product $\subscr{k}{i}^{\textup{in}}\subscr{k}{p}^{\textup{in}}$ appears as a single scaling factor in~\eqref{eq:pd_gradient_flow_ineq}, the two gains play distinct roles. The proportional gain shapes the augmented Lagrangian and controls how aggressively feasibility violations are penalized in the primal dynamics.
The integral gain governs the timescale of the dual dynamics, i.e., how rapidly $\xi$ tracks $\mu$.

\begin{remark}[Comparison with~\cite{VC-SMF-SP-DR:24b}]
The work~\cite{VC-SMF-SP-DR:24b} proposes a PI-based approach to inequality-constrained optimization sharing the same control-theoretic motivation, but differs from the SPPI in two key respects. First, their dual update injects a state-velocity term $\subscr{k}{p}^{\textup{in}} J_g(x)\dot{x}$, introduced to avoid a discontinuity in the natural PI extension. By contrast, the SPPI is derived directly from the back-calculation anti-windup scheme, thereby preserving a clearer control-theoretic interpretation. Second, their framework is restricted to affine constraints, whereas the SPPI applies to any $\mathcal{C}^1$ function.\end{remark}
\subsection{Convergence Analysis of SPPI in Convex Settings}
\label{sec:convergence_convex_ineq}
\begin{arxiv}
We analyze the convergence properties of SPPI~\eqref{eq:pi_spf_ineq} for
convex problems under Assumption~\ref{ass:1}, where Slater's condition is
interpreted for the inequality-constrained setting. Namely, there exists
$\bar{x} \in \R^n$ such that $g(\bar{x}) < \0_m$.
\begin{mycorollary}[label=thm:gas_ineq]{Convergence to the KKT set}
Under Assumption~\ref{ass:1}, consider the SPPI~\eqref{eq:pi_spf_ineq}. For every $\subscr{k}{p}^{\textup{in}}$, $\subscr{k}{i}^{\textup{in}} > 0$,
\begin{enumerate}[label=\textup{(\roman*)}]
\item every trajectory converges to the KKT set $\mathcal{K}$ of~\eqref{eq:ineq_constrained}.
\label{thm_conv_ineq:item1}
\item if additionally $f$ is strictly convex and LICQ holds, then~\eqref{eq:ineq_constrained} has a unique KKT point $z^{\star} = (x^{\star},\xi^{\star})$, which is globally asymptotically stable.
\label{thm_conv_ineq:item2}
\end{enumerate}
\end{mycorollary}
\begin{proof}
Let $(x^\star, \xi^\star)$ be any KKT point of problem~\eqref{eq:ineq_constrained}, which exists by  Assumption~\ref{ass:1}~\cite{DPB:97}.
Throughout, we write $r(x,\xi) := \relu\bigl(\xi + \subscr{k}{p}^{\textup{in}} g(x)\bigr)$ for brevity.
Consider the following Lyapunov function candidate
$$
V (x, \xi) = \frac{1}{2}\norm{x - x^\star}^2 + \frac{1}{2 \subscr{k}{i}^{\textup{in}} \subscr{k}{p}^{\textup{in}}}\norm{\xi - \xi^\star}^2.
$$
Clearly, $V(x, \xi)$ is positive definite, radially unbounded, and vanishes only at $(x^\star,\xi^\star)$.
 
\textbf{Computing $\dot{V}$.}
Differentiating along trajectories of~\eqref{eq:pi_spf_ineq} gives
\begin{align}
\dot{V}(x, \xi) &= (x - x^\star)^\top \dot x + \frac{1}{\subscr{k}{i}^{\textup{in}} \subscr{k}{p}^{\textup{in}}} (\xi - \xi^\star)^\top \dot \xi \nonumber\\
&= (x - x^\star)^\top \Bigl(- \nabla f(x) - J_g(x)^\top r(x, \xi)\Bigr) + \frac{1}{\subscr{k}{i}^{\textup{in}} \subscr{k}{p}^{\textup{in}}} (\xi - \xi^\star)^\top \Bigl( - \subscr{k}{i}^{\textup{in}} \xi + \subscr{k}{i}^{\textup{in}} r(x, \xi) \Bigr) \label{eq:vdot_lasalle}\\
&= - (x - x^\star)^\top \nabla f(x) - (x - x^\star)^\top J_g(x)^\top r(x, \xi) - \frac{1}{\subscr{k}{p}^{\textup{in}}} (\xi - \xi^\star)^\top \xi +  \frac{1}{\subscr{k}{p}^{\textup{in}}} (\xi - \xi^\star)^\top r(x, \xi). \nonumber
\end{align}
Using the stationarity condition~\eqref{kkt_stationarity}, $\nabla f(x^\star) = - J_g(x^\star)^\top \xi^\star$, and adding and subtracting $(x-x^\star)^\top \nabla f(x^\star)$ to $\dot{V}(x, \xi)$ yields
\begin{align}
\dot{V}(x, \xi) &= - (x - x^\star)^\top \bigl(\nabla f(x) - \nabla f(x^\star)\bigr) + (x - x^\star)^\top J_g(x^\star)^\top \xi^\star - (x - x^\star)^\top J_g(x)^\top r(x, \xi)\\
&\quad \quad - \frac{1}{\subscr{k}{p}^{\textup{in}}} (\xi - \xi^\star)^\top \xi +  \frac{1}{\subscr{k}{p}^{\textup{in}}} (\xi - \xi^\star)^\top r(x, \xi) \nonumber \\
&= - (x - x^\star)^\top \bigl(\nabla f(x) - \nabla f(x^\star)\bigr) + (x - x^\star)^\top \bigl(J_g(x^\star)^\top \xi^\star - J_g(x)^\top r(x, \xi)\bigr) - \frac{1}{\subscr{k}{p}^{\textup{in}}} (\xi - \xi^\star)^\top \bigl(\xi - r(x, \xi)\bigr). \label{v_ineq_1}
\end{align}

\textbf{Bounding the Jacobian cross-term.}
Since $g$ is convex and continuously differentiable, the first-order inequality $g(z) - g(y) \geq J_g(y)(z-y)$ holds for all $z,y$~\cite{DPB:97}.
Evaluating the above inequality at $(z, y) = (x, x^\star)$ and $(z, y) = (x^\star, x)$, respectively, yields
\begin{align*}
g(x) - g(x^\star) \geq J_g(x^\star) (x - x^\star) \quad &\implies \quad (g(x) - g(x^\star))^\top \xi^\star \geq \big(J_g(x^\star)(x - x^\star)\big)^\top \xi^\star, \\
g(x^\star) - g(x) \geq J_g(x) (x^\star - x) \quad &\implies \quad -\big(g(x) - g(x^\star)\bigr)^\top r(x, \xi) \geq -\big(J_g(x)(x - x^\star)\big)^\top r(x, \xi),
\end{align*}
where the inequalities are preserved under multiplication since $\xi^\star \geq \0_m$ by KKT conditions and $r(x, \xi) \geq \0_m$ by the definition of the $\relu$ operator.
Substituting the above inequalities into~\eqref{v_ineq_1} gives
\begin{align*}
\dot{V}(x, \xi) &\leq - (x - x^\star)^\top \bigl(\nabla f(x) - \nabla f(x^\star)\bigr)+ \bigl(g(x) - g(x^\star)\bigr)^\top \bigl(\xi^\star - r(x, \xi)\bigr) - \frac{1}{\subscr{k}{p}^{\textup{in}}} (\xi - \xi^\star)^\top \bigl(\xi - r(x, \xi)\bigr).
\end{align*}

\textbf{Applying properties of $\relu$.}
The firmly non-expansive property of the $\relu$ operator allows us to establish the following bound on the constraint interaction terms:
\begin{align}
\bigl(g(x) - g(x^\star)\bigr)^\top& \bigl(\xi^\star - r(x, \xi)\bigr) \leq - \frac{1}{\subscr{k}{p}^{\textup{in}}} \norm{r(x, \xi) - \xi^\star}^2 + \frac{1}{\subscr{k}{p}^{\textup{in}}} \bigl(\xi - \xi^\star \bigr)^\top \bigl(r(x, \xi) - \xi^\star \bigr).
\label{eq:relu_nonexpansive_bound}
\end{align}
To see this, recall that firmly non-expansiveness of the $\relu$ implies 
$\norm{\relu(z) - \relu(y)}^2 \leq \bigl(z - y\bigr)^\top\bigl(\relu(z) - \relu(y) \bigr)$, for all $z,y \in \R^m$.
Inequality~\eqref{eq:relu_nonexpansive_bound} follows directly by selecting $z = \xi + \subscr{k}{p}^{\textup{in}} g(x)$ and $y = \xi^\star + \subscr{k}{p}^{\textup{in}} g(x^\star)$. 

\textbf{Combining the bounds.}
Substituting~\eqref{eq:relu_nonexpansive_bound} into~\eqref{v_ineq_1} and simplifying the dual terms, yields $-\tfrac{1}{\subscr{k}{p}^{\textup{in}}}\|r-\xi^\star\|^2 + \tfrac{1}{\subscr{k}{p}^{\textup{in}}}(\xi-\xi^\star)^\top(r-\xi^\star) - \tfrac{1}{\subscr{k}{p}^{\textup{in}}}(\xi-\xi^\star)^\top(\xi-r) = -\tfrac{1}{\subscr{k}{p}^{\textup{in}}}\|r-\xi\|^2$, where the last equality follows by expanding and cancelling all terms involving $\xi^\star$.
Therefore,
\begin{equation}
\label{eq:Vdot_final}
\dot{V}(x,\xi) \leq - (x-x^\star)^\top \bigl(\nabla f(x) -\nabla f(x^\star)\bigr) - \frac{1}{\subscr{k}{p}^{\textup{in}}}\|r(x,\xi)-\xi\|^2 \leq 0,
\end{equation}
where the last inequality uses convexity (hence monotonicity) of $\nabla f$. Hence, $V$ is nonincreasing along trajectories of~\eqref{eq:pi_spf_ineq}, every trajectory is bounded, and $(x^\star,\xi^\star)$ is stable.

\textbf{LaSalle's invariance principle.}
Let $\mathcal{E} := \setdef{(x,\xi)\in\R^{n+m}}{\dot{V}(x,\xi)=0}$ and let $\Omega\subseteq\mathcal{E}$ denote the largest invariant subset of $\mathcal{E}$. Consider a trajectory $(x(t),\xi(t))\in\Omega$. From~\eqref{eq:Vdot_final}, we have
\begin{align}
&(x-x^\star)^\top \bigl(\nabla f(x) -\nabla f(x^\star)\bigr) = 0, \label{eq:vdot_primal}\\
&\relu\bigl(\xi+\subscr{k}{p}^{\textup{in}}g(x)\bigr) - \xi = \0_m \label{eq:vdot_dual}.
\end{align}
Equation~\eqref{eq:vdot_dual} implies $\dot{\xi} = -\subscr{k}{i}^{\textup{in}}\xi + \subscr{k}{i}^{\textup{in}}r(x,\xi) = \0_m$, and therefore $\xi(t) = \bar{\xi}$, for all $t \geq 0$.
Moreover, the fixed-point condition~\eqref{eq:vdot_dual} implies primal feasibility, dual feasibility, and complementary slackness~\eqref{kkt_primal_feasibility}--\eqref{kkt_complementarity} in $(x, \bar{\xi})$, namely, $g(x) \leq \0_m$, $\bar{\xi} \geq \0_m$, 
and $\bar{\xi}^\top g(x) = 0$.
Since $\bar{\xi}$ is constant on $\Omega$, the primal dynamics reduce to $\dot{x} = - \nabla f(x) - J_g(x)^\top \bar{\xi}$. Define the Lagrangian $L(x,\xi) = f(x) + g(x)^\top \xi$, which is convex due to $f$ and $g$ convex.
The primal dynamics read $\dot{x} = -\nabla L(x,\bar{\xi})$. To understand the invariant set for the condition~\eqref{eq:vdot_primal}, consider the first term of the Lyapunov candidate $\phi(x) = 1/2 \|x - x^\star\|^2$. Along trajectories in $\Omega$, we have
\begin{align*}
\dot{\phi}(x) &= (x - x^\star)^\top\dot{x} = -(x - x^\star)^\top \nabla L(x,\bar{\xi}) \leq L(x^\star,\bar{\xi}) - L(x,\bar{\xi}) \leq f(x^\star) - f(x) \leq 0,
\end{align*}
where the first inequality follows from convexity of $L(\cdot,\bar{\xi})$, the second uses $L(x^\star,\bar{\xi}) \leq f(x^\star)$ due to $\bar{\xi}\geq \0_m$ and $g(x^\star)\leq \0_m$, and $L(x,\bar{\xi}) = f(x)$ due to complementary slackness.
The last inequality follows from the optimality of $x^\star$ and the feasibility of $x$.
So $\phi(x)$ is nonincreasing along trajectories in $\Omega$, and the feasible set is forward invariant, implying all solutions converge to the largest invariant subset of $\{\dot{\phi}=0\}$.
On this invariant set, $0=\dot{\phi} \leq f(x^\star)-f(x) \leq 0$, which implies $f(x^\star) = f(x)$.
Thus, $x$ minimizes $L(\cdot,\bar\xi)$ by convex sufficiency, giving $\nabla f(x) + J_g(x)^\top\bar{\xi} = \0_n$, implying $\dot{x} = \0_n$.
Therefore, every point in $\Omega$ satisfies the KKT conditions and, by LaSalle's invariance principle, the $\omega$-limit set of every solution is contained in $\Omega$, and we conclude $\operatorname{dist}((x, \xi), \mathcal{K}) \to 0$.
This proves item~\ref{thm_conv_ineq:item1}.

Item~\ref{thm_conv_ineq:item2} follows by noticing that strict convexity  of $f$ and LICQ together imply that the KKT set reduces to the singleton  $\{(x^\star, \xi^\star)\}$. Global asymptotic stability then follows from LaSalle's invariance principle combined with radial unboundedness of $V$. This concludes the proof.
\end{proof}
\end{arxiv}
Under additional regularity assumptions on $z^\star$, we can establish exponential convergence following an initial transient phase in which the decay is only linearly bounded.
\begin{arxiv}
\begin{mycorollary}
[label=thm:lin_exp_ineq]{Linear-exponential convergence of SPPI}
Let Assumption~\ref{ass:1} hold and let $z^\star = (x^\star, \xi^\star)$ be a KKT point of~\eqref{eq:ineq_constrained}. Suppose that strict complementarity, LICQ, and SOSC hold at $z^\star$. Then $z^\star$ unique and, for every $\subscr{k}{p}^{\textup{in}}, \subscr{k}{i}^{\textup{in}} > 0$, every trajectory of~\eqref{eq:pi_spf_ineq} linearly-exponentially converges to $z^\star$.
\end{mycorollary}
\begin{proof}
We show that SPPI is globally weakly contracting and locally strongly contracting around $z^\star$.
Let $z = (x, \xi) \in \R^{n+m}$ denote the full state vector, $\dot{z} = F(z)$ the SPPI~\eqref{eq:pi_spf_ineq}, and let $P = \begin{bmatrix} I_n & 0 \\ 0 & \frac{1}{\subscr{k}{i}^{\textup{in}} \subscr{k}{p}^{\textup{in}}} I_m \end{bmatrix} \succ 0$.
Moreover, let $y := \xi + \subscr{k}{p}^{\textup{in}} g(x)$ and define $G(y) := \nabla \relu(y)$ where it exists. By Rademacher's theorem, the Jacobian of the SPPI~\eqref{eq:pi_spf_ineq} exists almost everywhere and is given by
$$
J(z)=
\begin{bmatrix}
-\nabla^2 f(x) - \sum_{j=1}^m \relu\big(\xi_j + \subscr{k}{p}^{\textup{in}} g_j(x)\big) \nabla^2 g_j(x) - \subscr{k}{p}^{\textup{in}} J_g(x)^\top G(y) J_g(x) & -J_g(x)^\top G(y) \\ 
\subscr{k}{i}^{\textup{in}}\subscr{k}{p}^{\textup{in}}G(y) J_g(x) & - \subscr{k}{i}^{\textup{in}} (I_m - G(y))
\end{bmatrix}.
$$
To prove globally weak contractivity, it suffices to show that for all $z$ for which $J(z)$ exists, we have $\mu_{P}(J(z)) \leq 0$ or equivalently $\frac{1}{2}\left(P J(z) + J(z)^\top P\right) \preceq 0$ for all $z \in \R^{n+m}$.
We have 
\begin{align*}
\frac{1}{2}\left(P J(z) + J(z)^\top P\right) =
\begin{bmatrix} 
-\nabla^2 f(x) - \sum_{j=1}^m \relu \big(\xi_j + \subscr{k}{p}^{\textup{in}} g_j(x)\big) \nabla^2 g_j(x) - \subscr{k}{p}^{\textup{in}} J_g(x)^\top G(y) J_g(x) & 0 \\ 
0 & -\frac{1}{\subscr{k}{p}^{\textup{in}}}(I_m - G(y)) 
\end{bmatrix},
\end{align*}
which is block-diagonal. The LMI then follows by convexity of $f$ and $g_j$, and the condition $0 \preceq G(y) \preceq I_m$, for almost every $y$.

Next, we prove local strong contraction around $z^\star$. Write $\mathcal{A} := \mathcal{A}(x^\star)$ for the active set and $\mathcal{I} := \{1,\dots,m\}\setminus\mathcal{A}$ for its complement. By strict complementarity, for $j \in \mathcal{A}$, $g_j(x^\star) = 0$ and $\xi^\star_j > 0$ imply $y^\star_j > 0$ and thus $G_j(y^\star) = 1$. For $j \notin \mathcal{A}$, $g_j(x^\star) < 0$ and $\xi^\star_j = 0$ imply $y^\star_j < 0$, hence $G_j(y^\star) = 0$.
In particular, $y^\star_j \neq 0$ for all $j$ and, by continuity, in a sufficiently small local neighborhood $\mathcal{B}(z^\star)$ of $z^\star$, the activation pattern $G$ remains constant. The inactive multiplier dynamics decouple as a Hurwitz linear subsystem $\dot{\xi}_{\mathcal{I}} = -\subscr{k}{i}^{\textup{in}} \xi_{\mathcal{I}}$. The remaining active subsystem coordinates $(x, \xi_{\mathcal{A}})$ have the local Jacobian:
$$
\subscr{J}{loc} = \begin{bmatrix} - \nabla^2 f(x^\star) - \sum_{j \in \mathcal{A}} \xi^\star_j \nabla^2 g_j(x^\star) - \subscr{k}{p}^{\textup{in}} J_g^{\mathcal{A}}(x^\star)^\top J_g^{\mathcal{A}}(x^\star) & -J_g^{\mathcal{A}}(x^\star)^\top \\ \subscr{k}{i}^{\textup{in}} \subscr{k}{p}^{\textup{in}} J_g^{\mathcal{A}}(x^\star) & 0 \end{bmatrix}.
$$
By SOSC, $H := \nabla^2 f(x^\star) + \sum_{j \in \mathcal{A}} \xi^\star_j \nabla^2 g_j(x^\star) \succ 0$ on $\ker J_g^{\mathcal{A}}(x^\star)$, and by LICQ, $B := J_g^{\mathcal{A}}(x^\star)$ has full row rank. By the same argument as in the proof of Theorem~\ref{thm:lin_exp_eq_ineq} (with the equality-constraint terms absent), $H + \subscr{k}{p}^{\textup{in}} B^\top B \succ 0$. Hence $\subscr{J}{loc}$ has the saddle form of Lemma~\ref{lemma:saddle-matrices_general} and is Hurwitz for any $\subscr{k}{p}^{\textup{in}}, \subscr{k}{i}^{\textup{in}} > 0$.
Because the entire Jacobian $J(z^\star)$ is Hurwitz, the unique equilibrium is locally exponentially stable. Then, the system~\eqref{eq:pi_spf_ineq} is strongly infinitesimally contracting in a neighborhood of the equilibrium point $z^\star$~\cite{TS:75}. The statement then follows by applying~\cite[Theorem 2]{VC-AD-AG-GR-FB:24a}. This concludes the proof.
\end{proof}
\end{arxiv}
\subsection{Convergence Analysis for Affine Constraints}
\label{sec:convergence_ineq}
We specialize SPPI to the affine setting $g(x) := Cx - d$, where $C \in \R^{m \times n}$ and $d \in \R^m$, giving
\beq
\label{eq:pi_spf_ineq_affine}
\begin{cases}
\dot{x} = - \nabla f(x) - C^\top \relu\big(\xi + \subscr{k}{p}^{\textup{in}} (Cx - d)\big), \\
\dot{\xi} = - \subscr{k}{i}^{\textup{in}} \xi + \subscr{k}{i}^{\textup{in}} \relu\big(\xi + \subscr{k}{p}^{\textup{in}} (Cx - d)\big).
\end{cases}
\eeq
Under precise parameter identifications, SPPI~\eqref{eq:pi_spf_ineq_affine} recovers two primal--dual systems independently studied in the literature for solving $\min\setdef{f(x)}{Cx \leq d}$: the augmented primal--dual gradient dynamics (Aug-PDGD)~\cite{GQ-NL:19} and the proximal augmented Lagrangian primal--dual dynamics (Prox-PDGD)~\cite{AD-VC-AG-GR-FB:23f}.
The novelty of our formulation relative to~\cite{GQ-NL:19, AD-VC-AG-GR-FB:23f} lies not in the vector field~\eqref{eq:pi_spf_ineq_affine} itself, but in its derivation, interpretation, and generality. Both references derive these dynamics top-down, by differentiating a nonsmooth augmented penalty function. In contrast, we derive the same dynamics bottom-up from classical control principles. This control-theoretic perspective formally bridges controller design and augmented Lagrangian methods for inequality-constrained optimization. Importantly, unlike~\cite{GQ-NL:19, AD-VC-AG-GR-FB:23f}, our framework is not limited to affine constraints: the SPPI applies to any constraint function $g \in \mathcal{C}^1$, as developed in Section~\ref{sec:inequality_constraint}.

Under Assumption~\ref{ass:1_eq} on $f$ and $C$, global exponential stability and strong infinitesimal contractivity of~\eqref{eq:pi_spf_ineq_affine} therefore follow directly from~\cite[Theorem~2]{GQ-NL:19} and~\cite[Theorem~11]{AD-VC-AG-GR-FB:23f}, respectively. However, the contractivity certificate of~\cite[Theorem 11]{AD-VC-AG-GR-FB:23f} offers limited analytical insight, since both the contraction metric $P$ and the rate $\subscr{c}{in}$ are defined implicitly through a nonlinear program. We address this limitation by providing closed-form expressions for $P$ and $\subscr{c}{in}$ directly in terms of the gains.
Notably, direct numerical comparisons show that $\subscr{c}{in}$ exceeds the rate in~\cite[Theorem~2]{GQ-NL:19} by several orders of magnitude across every tested combination of gains and problem parameters. Under the analogous identification with Prox-PDGD, $\subscr{c}{in}$ exceeds the numerically-optimized $c^\star$ of~\cite[Theorem~11]{AD-VC-AG-GR-FB:23f} across most tested regimes, with the gap narrowing as $\amax/\amin$ grows large. We refer to the extended version's appendix~\cite{VC-RH-ECB-JL:26} for plots showing these comparisons.
\begin{arxiv}
\begin{mytheorem}[label=thm:pi_spf_ineq_affine-contractivity]{Contractivity of~\eqref{eq:pi_spf_ineq_affine}}
Under Assumption~\ref{ass:1_eq} on $f$ and $C$, for any $\subscr{k}{p}^{\textup{in}}$, $\subscr{k}{i}^{\textup{in}} > 0$, SPPI~\eqref{eq:pi_spf_ineq_affine} is $\subscr{c}{in}$-strongly infinitesimally contracting with respect to $\|\cdot\|_{P}$, where
\beq
\label{eq:P_generalized_saddle}
P =
\begin{bmatrix}
\gamma I_n & C^\top \\
C & \dfrac{\gamma}{\subscr{k}{i}^{\textup{in}}\subscr{k}{p}^{\textup{in}}} I_m
\end{bmatrix}
\succ 0,
\qquad
\subscr{c}{in} := \frac{3\subscr{k}{i}^{\textup{in}}\subscr{k}{p}^{\textup{in}}\subscr{\textup{c}}{min}}{8\gamma},
\eeq
and $\gamma>0$ is given by $\gamma := \max\left(\subscr{k}{i}^{\textup{in}}, \subscr{k}{p}^{\textup{in}}\subscr{\textup{c}}{max},
\dfrac{2}{3} L + \dfrac{3\subscr{k}{i}^{\textup{in}}\subscr{k}{p}^{\textup{in}}\subscr{\textup{c}}{min}}{8 \rho} +\dfrac{2}{3} \dfrac{x(x+2L)}{\rho}\right) + \eps$, with $\eps > 0$ and $x := \subscr{k}{p}^{\textup{in}} \subscr{\textup{c}}{max} + \dfrac{3\subscr{k}{i}^{\textup{in}}\subscr{\textup{c}}{min}}{4 \subscr{\textup{c}}{max}} + \subscr{k}{i}^{\textup{in}}\sqrt{\dfrac{\subscr{\textup{c}}{max}}{\subscr{\textup{c}}{min}}}$.
\end{mytheorem}
\end{arxiv}%
\begin{proof}
Let $z = (x,\xi) \in \R^{n + m}$, and $\dot{z} = \bar{\operatorname{F}}(z)$ denote~\eqref{eq:pi_spf_ineq_affine}. Moreover, let $y := \xi + \subscr{k}{p}^{\textup{in}} (Cx - d)$ and define $G(y) := \nabla \relu (y)$ where it exists. By Rademacher's theorem, the Jacobian of the SPPI~\eqref{eq:pi_spf_ineq_affine} exists almost everywhere and is
$$
{J}_{\bar{\operatorname{F}}}(z) =
\begin{bmatrix}
- \nabla^2 f(x) - \subscr{k}{p}^{\textup{in}} C^\top G(y) C & - C^\top G(y) \\
\subscr{k}{i}^{\textup{in}}\subscr{k}{p}^{\textup{in}}G(y)C & - \subscr{k}{i}^{\textup{in}} (I_m - G(y))
\end{bmatrix}.
$$
Assumption~\ref{ass:1_eq}.\ref{ass:1_eq_f} implies $\rho I_n \preceq \nabla^2 f(x) \preceq L I_n$, moreover $0 \preceq G(y) \preceq I_m$.
\begin{arxiv}
Then,
\begin{multline*}
{\sup_{z} \mu_{P}({J}_{\bar{\operatorname{F}}}(z)) \leq \max_{\substack{0 \preceq G \preceq I_m \\
\rho I_n \preceq B \preceq L I_n}} \mu_{P}\left(
\begin{bmatrix}
- B - \subscr{k}{p}^{\textup{in}} C^\top G C & - C^\top G \\
\subscr{k}{i}^{\textup{in}}\subscr{k}{p}^{\textup{in}}GC & \subscr{k}{i}^{\textup{in}} (G -I_m)
\end{bmatrix}
\right),}
\end{multline*}
\end{arxiv}%
The result then follows from Lemma~\ref{lemma:generalized-saddle_ineq} with $\subscr{b}{min} = \rho$, $\subscr{b}{max} = L$, $\alpha=\subscr{k}{p}^{\textup{in}}$, and $\beta=\subscr{k}{i}^{\textup{in}}$.
\end{proof}
As in the equality-constrained case (Theorem~\ref{thm:riemannian_saddle_affine-contractivity}), contractivity in Theorem~\ref{thm:pi_spf_ineq_affine-contractivity} holds for every choice of gains $\subscr{k}{p}^{\textup{in}}, \subscr{k}{i}^{\textup{in}} > 0$. However, the gains affect the rate. In fact, for $\subscr{k}{i}^{\textup{in}}$ fixed, $\gamma$ is non-decreasing in $\subscr{k}{p}^{\textup{in}}$, and grows unboundedly as $\subscr{k}{p}^{\textup{in}}\to\infty$, so $\subscr{c}{in}$ vanishes in this limit. This mirrors a structural difference between the two controllers: the equality controller~\eqref{eq:controller_eq} is linear in $(x,\nu)$, whereas the inequality controller~\eqref{eq:controller_ineq} saturates through $\relu(\cdot)$, and an overly aggressive proportional gain repeatedly drives the anti-windup mechanism into saturation, degrading the achievable rate in practice.
\section{Numerical Examples}
\label{sec:numerical_example}
We now illustrate the application of SPPID in three settings: (i) executing sequential quadratic programming (SQP) on the Rosenbrock-Suzuki problem~\cite{HW-KS:80}, (ii) solving a Stackelberg game, and (iii) implementing closed-loop nonlinear model predictive control (NMPC) of an inverted pendulum. In each case, SPPID is implemented with a discrete-time step $\Delta t$ tuned to the specific formulation.
The intuition underlying the choice of SPPID gains is as follows:
\begin{itemize}
    \item Increasing the proportional gains generally speeds up SPPID. However, an excessively large $k_\text{p}^\text{in}$ saturates the anti-windup integrator and degrades convergence.
    \item The integral gain governs the time-scale separation between the primal and dual dynamics; together with $\Delta t$, it sets the effective step size of the primal--dual update.
    \item Starting from an infeasible point, increasing $k_\text{d}^\text{eq}$ promotes trajectories that are tangential to the constraint set, manifesting as oscillations. Starting from a feasible point, increasing $\subscr{k}{d}^{\textup{eq}}$ accelerates convergence toward the optimum. In discrete-time, a nonzero $\subscr{k}{d}^{\textup{eq}}$ reduces the stiffness of the resulting dynamics, enabling larger integration steps. 
\end{itemize}
\subsection{Sequential Quadratic Programming}
SPPID is embedded within an SQP framework as shown in Figure~\ref{fig:SQP Diagram}. The upper-level step size $\alpha_k$ is determined via a backtracking line search.
\begin{arxiv}
\begin{figure}[h]
    \centering
    \includegraphics[width=0.5\linewidth]{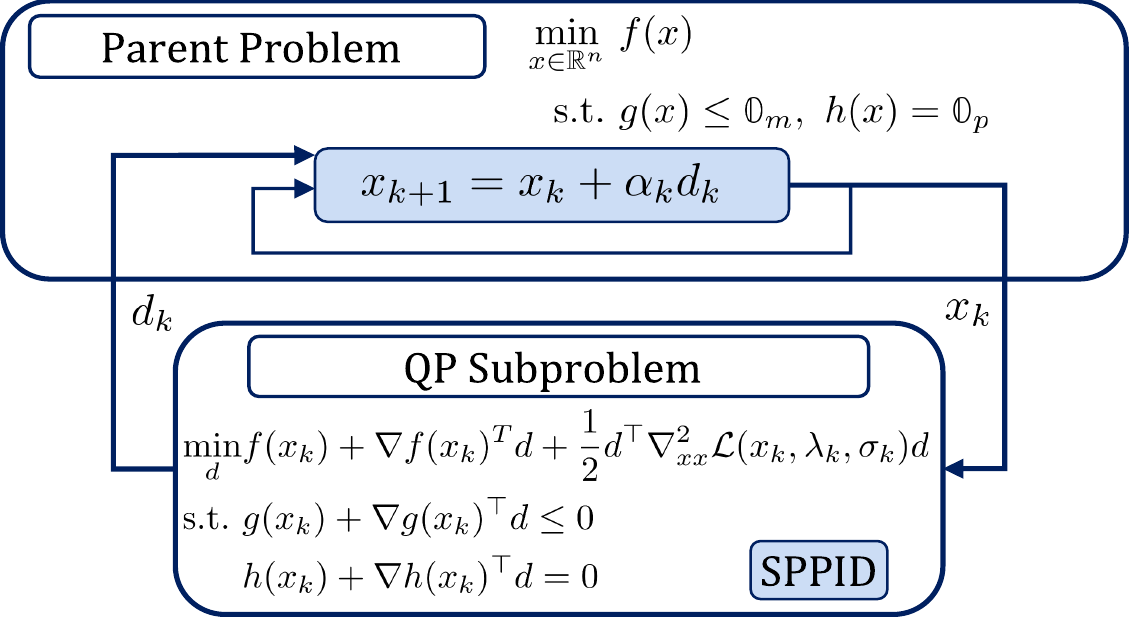}
    \caption{Architecture of an SQP optimizer incorporated with SPPID which determines the search direction}
    \label{fig:SQP Diagram}
\end{figure}
\end{arxiv}
We demonstrate this framework by solving the nonconvex Rosenbrock--Suzuki problem with solution $x^\star = (0,1,2,-1)$~\cite{HW-KS:80}:
\begin{equation*}\label{OP:Rosenbrock-Suzuki}
\begin{aligned}
\min_{x\in\R^4} & f(x) = x_1^2 + x_2^2 + 2x_3^2 + x_4^2 - 5x_1 - 5x_2 - 21x_3 + 7x_4, \\
\text{s.t. }&g_1(x) = x_1^2 + x_2^2 + x_3^2 + x_4^2 + x_1 - x_2 + x_3 - x_4\le 8,\\
&g_2(x) = x_1^2 + 2x_2^2 + x_3^2 + 2x_4^2 - x_1 - x_4 \le 10,\\
&g_3(x) = 2x_1^2 + x_2^2 + x_3^2 + 2x_1 - x_2 - x_4 \le 5.
\end{aligned}
\end{equation*} 
We sample $N = 10$ random initial points with components drawn uniformly over $[-10,10]$, and study the effect of the SPPID gains $\subscr{k}{p}^\text{in}$ and $\subscr{k}{i}^\text{in}$ on convergence. As shown in Figure~\ref{fig:RS-SQP}, increasing the value of  both $\subscr{k}{p}^\text{in}$ and $\subscr{k}{i}^\text{in}$ corresponds to an improvement in the performance of SPPID. Specifically, having non-zero $\subscr{k}{p}^\text{in}$ results in quicker convergence of the iterates obtained from the framework in Figure~\ref{fig:SQP Diagram} to $x^\star$.
Furthermore, we compare SPPID against the feedback linearization-based algorithm with momentum (FL-M) of~\cite{RZ-AR-JS-NL:25}, which exhibits a discrete time update analogous to that of an SQP. As illustrated in Figure~\ref{fig:RS-SQP}, the SQP formulation with SPPID outperforms FL-M by several orders of magnitude across the tested values of $\subscr{{k}^\text{in}}{p}$ and $\subscr{{k}^\text{in}}{i}$.
\begin{arxiv}
\begin{figure}[h]
    \centering
    \includegraphics[width =0.6\linewidth, page = 1]{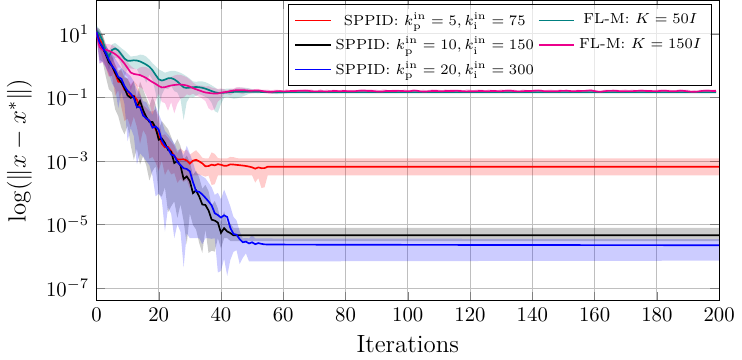}
    \caption{Convergence profiles for the Rosenbrock--Suzuki optimization problem solved utilizing SQP and SPPID for different controller gains. The central trajectory depicts the mean value evaluated over various initial conditions and the shaded area bounds the minimum and maximum values across iterations.}
    \label{fig:RS-SQP}
\end{figure}
\end{arxiv}
\subsection{Bilevel Optimization of Stackelberg Games}
\begin{arxiv}
Bilevel optimization problems have the form
\begin{equation}
\begin{aligned}
\min_{x \in \R^n} & f\bigl(x, y^\star(x)\bigr)\\
\text{s.t. } & y^\star(x) \in \argmin_{y \in \R^m} g(x,y),
\end{aligned}
\label{op:bilevel-opt}
\end{equation}
where the functions $\map{f}{\R^{n}\times\R^m}{\R}$ and $\map{g}{\R^{n}\times\R^m}{\R}$ are smooth.
\end{arxiv}%
We consider a multi-period Stackelberg game in which the leader exercises price-setting authority. Let $x = [x_1, \dots,x_T]^\top$ denote the leader’s price per quantity strategy over $T$ periods. The follower, observing $x$, solves a local minimization problem over each period $k$, choosing the quantity to produce $y = [y_1,\dots,y_T]$ to minimize the deficit $d_k= D_k- y_k$ and the total cost, where $D_k$ corresponds to the demand. Bilevel optimization is commonly used to model Stackelberg games, and the multi-period Stackelberg game can be cast into
\begin{equation}
\label{op:1L1F Bilevel}
\begin{aligned}
&\min_{x \in \R^T}  f_\textup{L}\bigl(x, y^\star(x)\bigr) := \sum_{k=1}^T-\gamma_\textup{L}^k x_ky_{k}^\star+c_1 {y_{k}^\star} ^2+c_2y_{k}^\star\\
& \begin{aligned}
\text{s.t. } & x\geq\0_T, \\
&y^\star(x) \in \argmin_{y \in \R^T} f_\textup{F}(x,y)\!:=\!\sum_{k=1}^T\!\gamma_\textup{F}^kx_k y_{k} + a_k({D}_{k}-y_k)^2\!,
\end{aligned}
\end{aligned}
\end{equation}
where $c_1,c_2>0$, and $a_k > 0$ for all $k$, ensuring the lower-level problem is strongly convex in $y$.
The costs for the leader and the followers are weighted by discount factors $\gamma_\textup{L}$ and $\gamma_\textup{F}$, respectively, to prioritize the minimization of near-term costs over future ones. 
The gradient of the lower-level objective with respect to $y$ is $h(x,y) := \nabla_y f_\textup{F}(x,y) = A_1x + A_2 y + b$, where 
\begin{arxiv}
    \begin{equation*}
    A_1 = 
    \begin{bmatrix} 
        \gamma_\textup{F} & 0 & \dots & 0 \\ 
        0 & \gamma_\textup{F}^2 & \dots & 0 \\ 
        \vdots & \vdots & \ddots & \vdots \\ 
        0 & 0 & \dots & \gamma_\textup{F}^T 
    \end{bmatrix}, \
    A_2 = \begin{bmatrix} 
       2a_1 & 0 & \dots & 0 \\ 
        0 & 2a_2 & \dots & 0 \\ 
        \vdots & \vdots & \ddots & \vdots \\ 
        0 & 0 & \dots & 2a_T 
    \end{bmatrix},
    \text{and }
    b = -\begin{bmatrix}
        2a_1 & 0 & \dots & 0 \\ 
        0 & 2a_2 & \dots & 0 \\ 
        \vdots & \vdots & \ddots & \vdots \\ 
        0 & 0 & \dots & 2a_T 
    \end{bmatrix}
    \begin{bmatrix}
        D_1\\
        D_2\\
        \vdots\\
        D_T
    \end{bmatrix}.
\end{equation*}
\end{arxiv}%
We compactly write $h(x,y) = A' [x,y]^\top+b$, with $A^\prime := [A_1 \quad A_2]$.
Using the first-order optimality condition for the strongly convex lower-level problem, we obtain the equivalent single-level constrained optimization problem
\begin{equation*} 
\label{op:bilevel_constrained}
\begin{aligned}
\min_{x \in \R^T, y \in \R^T}&  f_\textup{L}(x,y)\\
\text{s.t. } & g(x,y) := -x \leq \0_T, \ h(x,y) := \nabla_y f_\textup{F}(x,y) = \0_T,
\end{aligned}
\end{equation*}
which fits the structure in~\eqref{eq:eq_ineq}.
Numerically, we set $\gamma_\textup{L} = 0.99$ and $\gamma_\textup{F} = 0.95$, and choose arbitrary values for $a_1,\dots,a_T$, and $D$. Moreover, we set $c_1=10$ and $c_2 = 10$, and $\Delta t = 0.01$.

For visualization purposes, we first consider a single period ($T=1$) game. 
Figure~\ref{fig:BiOPT-GT-SINGLE} shows the contour plot of the objective function and the constraint set, as well as the trajectory of the iterates for various values of $k_\textup{d}^\text{eq}$, initialized from a feasible and an infeasible point. Starting from an infeasible point, as the derivative gain increases, the resulting trajectories exhibit stronger oscillations and evolve tangentially along the equality constraints toward the optimal point. We also compare against the PDGD, i.e., SPPID~\eqref{eq:sp_pid_eq_ineq} with $\subscr{k}{p}^{\textup{eq}} = \subscr{k}{d}^{\textup{eq}} = 0$, which fails to converge and oscillates without converging to the optimal point. Starting from a feasible point, increasing $k_\textup{d}^\text{eq}$ promotes trajectories that are parallel to the constraint manifold, leading to faster convergence of SPPID towards the optimum.
\begin{arxiv}
\begin{figure}[!h]
    \centering
    \begin{subfigure}[b]{\columnwidth}
        \centering
        \includegraphics[width=0.6\linewidth, page=1]{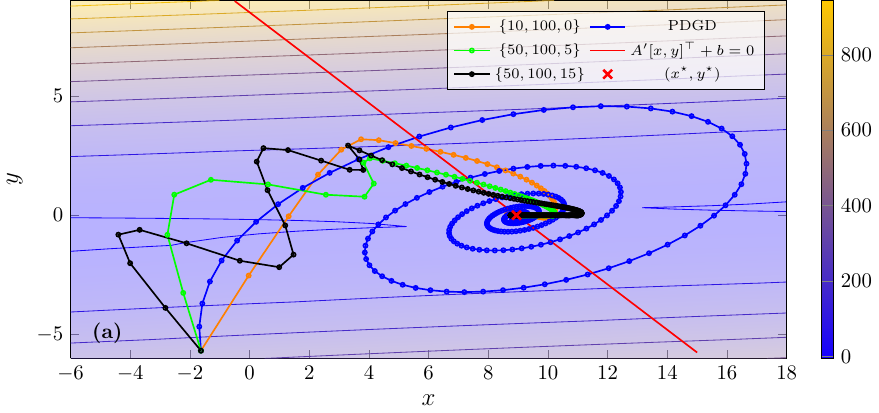}
        \label{fig:BiOPT-GT-MULTI-infeasible}
    \end{subfigure}
    \begin{subfigure}[b]{\columnwidth}
        \centering
        \includegraphics[width=0.6\linewidth, page=2]{figures/tac/Bilevel_2D.pdf} 
        \label{fig:BiOPT-GT-MULTI-feasible}
    \end{subfigure}
    \caption{SPPID convergence profiles for the solution of~\eqref{op:1L1F Bilevel} with different controller gains: (a) starting from an infeasible initial point, and (b) starting from a feasible initial point. The SPPID gains are represented in the format $\{k^\textup{eq}_\text{p},k^\textup{eq}_\text{i},k^\textup{eq}_\text{d}\}$ and $\Delta t = 0.01$, $k^\textup{in}_\text{p} = k^\textup{eq}_\text{p}$, and $k^\textup{in}_\text{i} = 100$.}
    \label{fig:BiOPT-GT-SINGLE}
\end{figure}
\end{arxiv}

Moreover, we consider a game with $T = 10$ periods and analyze 
the effect of increasing $k^\textup{in}_\text{p}$. We sample for $10$ different initial conditions, with each component sampled independently and uniformly from the interval $[-20, 20]$. 
Figure~\ref{fig:BiOPT-GT-MULTI} illustrates the logarithm of the distance error. For the considered gains, increasing $k_\text{p}^\text{in}$ and $k_\text{p}^\text{eq}$ while keeping the other gains fixed, leads to faster convergence of SPPID towards the optimum. However, only increasing $k_\text{p}^\text{in}$ too much saturates the PI controller within the anti-windup scheme, thereby inhibiting convergence, as evident when $k^\textup{in}_\text{p}$ increases from $45$ to $60$.
\begin{arxiv}
\begin{figure}[!h]
    \centering
    \includegraphics[width = .7\linewidth, page = 1]{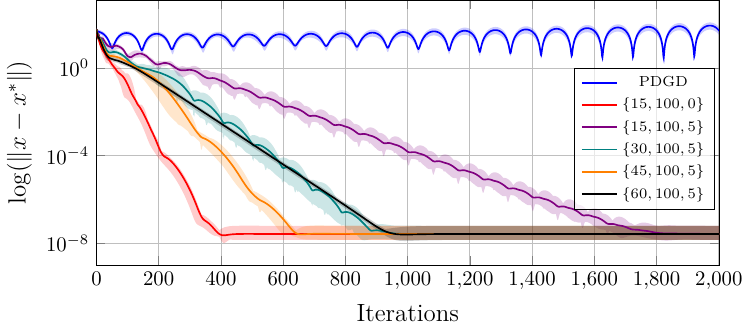}
    \caption{SPPID convergence profiles for the solution of~\eqref{eq:eq_ineq} with different controller gains. The central trajectory depicts the mean value evaluated over various initial conditions and the shaded area bounds the minimum and maximum values across iterations. The SPPID gains are represented in the format $\{k^\textup{eq}_\text{p},k^\textup{eq}_\text{i},k^\textup{eq}_\text{d}\}$ and $\Delta t = 0.01$, $k^\textup{in}_\text{p} = k^\textup{eq}_\text{p}$, and  $k^\textup{in}_\text{i} = 100$.}
    \label{fig:BiOPT-GT-MULTI}
\end{figure}
\end{arxiv}
\subsection{Closed-Loop Nonlinear Suboptimal MPC}
The dynamics of the inverted pendulum are governed by the state vector $x = [\theta, \dot\theta]^\top \in \mathbb{R}^2$, where $\theta$ denotes the angular position and $\dot\theta$ represents the angular velocity. The nonlinear equations of motion are expressed as $\dot x = f_\text{p}(x,u) := \begin{bmatrix}\dot\theta, \ \frac{g}{l}\sin(\theta) - \frac{b}{m l^2}\dot\theta + \frac{1}{m l^2}u\end{bmatrix}^\top$,
where $m$ is the mass, $l$ is the pendulum length, $b$ is the viscous friction coefficient, $g$ is the gravitational acceleration, and $u$ is the control input. 
To stabilize the pendulum at its upright equilibrium point, a nonlinear model predictive control (NMPC) scheme is employed. The resulting optimal control problem is solved via SPPID in a finite number of iterations $\subscr{N}{iter}$ and is formulated as 
\begin{equation*}
    \begin{aligned}
     \min_{\mathbf{x}, \mathbf{u}}& \ x_{N_\text{p}}^\top P_\textup{F} x_{N_\text{p}} + \sum_{k=0}^{N_\text{p}-1} \left( x_k^\top Q x_k + u_k^\top R u_k \right)\\
     \text{s.t.} & \  x_0 = x_{\text{init}}, \ x_{k+1} = f_\text{p}^\text{d}(x_k, u_k), \ k = 0, 1, \dots, N_\text{p}-1,\\
    &\begin{aligned}
    & \sum_{k=0}^{N_\text{p}-1}u_k \leq u_\text{tot}, \\
    &\lvert u_k\rvert \le u_{\max}, \ \lvert u_k - u_{k-1}\rvert \le \Delta u_{\max}, \ k = 1, \dots, N_\text{p}-1,
    \end{aligned}
    \end{aligned}
\end{equation*}
where $P_\textup{F}, Q, R \succ 0$ are the terminal, state, and input weighting matrices, respectively, and $f_\text{p}^\text{d}(x_k, u_k)$ are the equivalent discrete time dynamics obtained via a forward Euler discretization of the dynamics with a time step of $dt = 0.1$. The initial state of the inverted pendulum is $x_\text{init} = [\frac{\pi}{4}, 0]^\top$ and the NMPC is implemented with a prediction horizon of $N_\text{p} = 20$. We also set $Q = \text{diag}(10.0,0.1)$, $R = 0.1$, and $P_\textup{F} =50Q$. The MPC formulation incorporates constraints on the control input, its rate of change, and the total control effort. 

We test a suboptimal MPC approach by reducing $N_\text{iter}$ for real-time applications. The role of $\subscr{k}{d}^{\textup{eq}}$ is highlighted, since a non-zero $\subscr{k}{d}^{\textup{eq}}$ decreases the number of iterations required, and allows the MPC to run at a higher frequency. Having $\subscr{k}{d}^{\textup{eq}} =0$ necessitates setting $\Delta t = 0.0001$, for which convergence is not achieved even with $\subscr{N}{iter} = 500$\begin{arxiv}
    as can be seen in Figure~\ref{fig:NMPC}
\end{arxiv}. On the other hand, we showcase the resulting closed-loop implementation for the suboptimal NMPC as solved using SPPID  with $\Delta t = 0.001$, $k^\text{eq}_\text{p} = k^\text{in}_\text{p} = 1000$, $k^\text{eq}_\text{i} = k^\text{in}_\text{i} = 1500$, and $k^\text{eq}_\text{d} = 10$ in Figure~\ref{fig:NMPC}. When $k^\text{eq}_\text{d} = 0$, the remaining gains are scaled appropriately to compensate for the smaller $\Delta t$, and SPPID is run for more iterations ($\subscr{N}{iter} = 2500$) to drive the solution closer to the optimum. We present the applied input and closed-loop behavior in comparison to the optimal solution obtained via Sequential Least Squares Programming (SLSQP).

Though increasing $\subscr{k}{d}^{\textup{eq}}$ may slow convergence, it reduces the numerical stiffness of the flow, permitting larger integration steps. In applications with limited computational time, such as real-time MPC, this trade-off favors a non-zero $\subscr{k}{d}^{\textup{eq}}$, since it reduces the number of iterations needed to converge. However, given the nonconvexity of the problem, increasing $\subscr{N}{iter}$ does not guarantee convergence of SPPID towards the optimum with a constant $\Delta t$, rather a varying step size is required. 

\begin{arxiv}
\begin{figure}[!h]
\centering
     \centering
    \includegraphics[width =0.6\linewidth, page = 1]{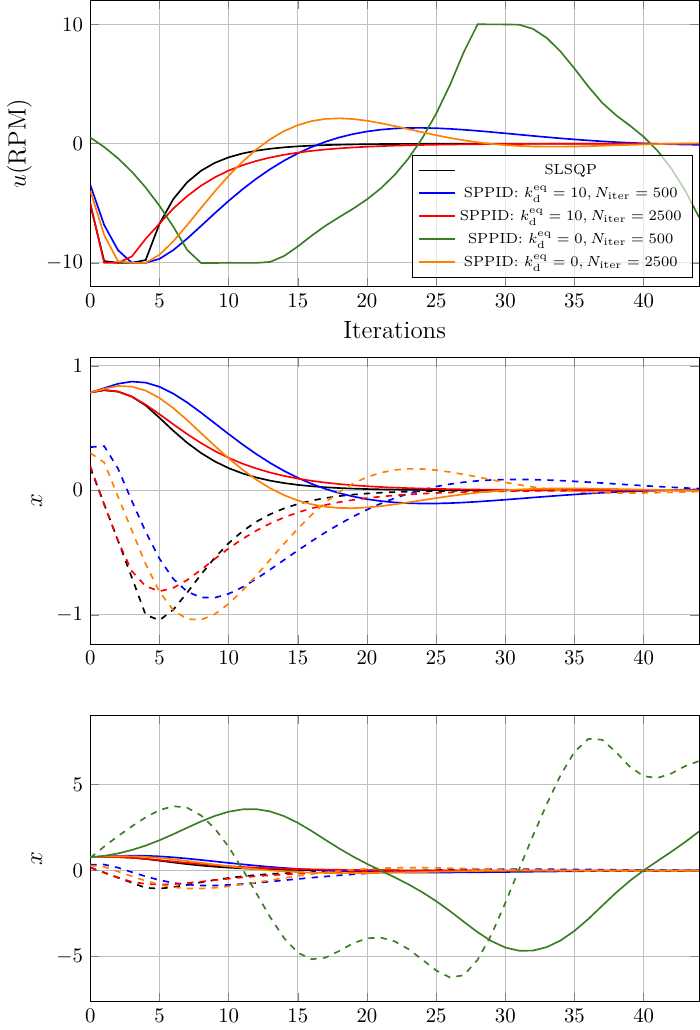}
    \caption{Application of NMPC solved using SPPID in closed-loop on inverted pendulum. Top: Closed-Loop NMPC input applied on inverted pendulum. Middle-Bottom: Inverted pendulum $\dot\theta$  (rad/s) (dashed lines) and $\theta$ (rad) (solid lines) trajectories under NMPC.}
    \label{fig:NMPC}
\end{figure}
\end{arxiv}
\section{Discussion and Conclusion}
\label{sec:discussion_conclusion}
We presented a unified control-theoretic framework for nonlinear constrained optimization based on PID control on the dual variables. The resulting \emph{saddle-point PID dynamics} (SPPID) systematically generate a family of continuous-time saddle-point flows associated with the augmented Lagrangian, whose equilibria coincide exactly with the KKT points of the original optimization problem. Moreover, the proposed framework provides a clear interpretation of the role of each feedback component.
For convex problems, we proved convergence to the KKT set, established conditions for global asymptotic stability, and, under additional regularity, showed linear-exponential convergence to the unique equilibrium. For equality-constrained problems, we recovered the projected gradient flow as the infinite-derivative-gain limit and proved local exponential convergence under non-affine constraints. For affine equality or inequality constraints with strongly convex, $L$-smooth costs, we established strong infinitesimal contractivity of SPPID (resp.\ SPPI), with explicit rate bounds.
Finally, we illustrated the effectiveness of our framework through three applications: non-convex SQP, a bilevel optimization solution for a multi-period Stackelberg game, and a closed-loop NMPC implementation solved via SPPID.

A natural extension is the study of discrete-time counterparts of the proposed dynamics. Other promising directions include extending the convergence analysis to nonconvex optimization problems and resolving the nonsmooth algebraic loop that currently precludes derivative feedback for inequality constraints (Remark~\ref{rem:kd_inequality}).

\bibliographystyle{IEEEtran}
\bibliography{main}

\ifcdc
    \appendices
\else
    \appendix
\fi
\renewcommand{\thelem}{\Alph{section}.\arabic{lem}}
\section{Logarithmic Norm of Saddle Matrices}
\label{apx:log_norm_of_saddle_matrix}
We now derive explicit bounds on the log-norm of the scaled saddle matrices that appear as Jacobians of the SPPID for equality- and inequality- constraints.
The next result extend the log-norm analysis of saddle matrices appearing in primal--dual dynamics in~\cite[Lemma 20]{AD-VC-AG-GR-FB:23f} and recover that result then when $M = I_n$ and $\Theta$ is scalar. The scalar case $\Theta = \tau^{-1}I_p$, $\tau\in\R_{>0}$, recovers the statement given in the conference version~\cite{VC-RH-ECB-JL:26a}.
\begin{arxiv}
\begin{mylemma}[label=lemma:saddle-matrices_general]{Log-norm of scaled saddle matrices}
Given $B=B^\top$, $M = M^\top \in\R^{n \times n}$, $M\succ 0$, $A\in\R^{p \times n}$, with $p\leq n$, {and $\Theta = \Theta^\top \in \R^{p\times p}$, $\Theta \succ 0$}, we consider the \emph{scaled saddle matrix}
\begin{equation*}
S =
\begin{bmatrix}
-M^{-1}B & - M^{-1}A^\top \\
{\Theta A} & 0
\end{bmatrix}
\in \R^{(p+n)\times (p+n)}.
\end{equation*}  
Then, for each matrix {quadruple $(M, B, A, \Theta)$} satisfying $\subscr{m}{min}I_n\preceq M \preceq \subscr{m}{max} I_n$, $\subscr{b}{min} I_n\preceq B \preceq \subscr{b}{max} I_n$, $\amin I_p \preceq A A^\top \preceq \amax I_p$, {and $\subscr{\theta}{min}I_p \preceq \Theta \preceq \subscr{\theta}{max}I_p$}, for $\subscr{m}{min}, \subscr{m}{max}, \subscr{b}{min}, \subscr{b}{max}, \amin,\amax{, \subscr{\theta}{min}, \subscr{\theta}{max}} \in \R_{>0}$, it holds
\begin{equation*}
S^\top P + PS\preceq -2 \subscr{c}{eq} P \quad\iff\quad \lognorm{S}{P} \leq - \subscr{c}{eq},
\end{equation*}
where
\begin{align}
P &=
\begin{bmatrix}
M & \alpha A^\top \\
\alpha A & {\Theta^{-1}}
\end{bmatrix}
\succ 0, \label{lem:def:P}\\
\alpha&=\frac{1}{2}\min\Big\{\frac{\subscr{m}{min}}{\subscr{b}{max}},{\frac{1}{\subscr{\theta}{max}}}\frac{\subscr{b}{min}}{\amax}\Big\}, 
\label{lem:def:alpha} \\
\subscr{c}{eq}&=\frac{1}{2}\alpha{\subscr{\theta}{min}}\frac{\amin}{\subscr{m}{max}}.
\label{lem:def:c}
\end{align}
\end{mylemma}
\begin{proof}
\textbf{Positivity of P.}\\
We start by verifying that $P\succ0$. Using the Schur complement of the $(2,2)$ entry, $P\succ0$ is equivalent to
$
M - \alpha^2 A^\top {\Theta} A\succ 0.
$
{Since $\Theta \preceq \subscr{\theta}{max}I_p$, for every $x\in\R^n$ we have $x^\top A^\top \Theta A x = (Ax)^\top \Theta (Ax) \leq \subscr{\theta}{max}\|Ax\|^2 = \subscr{\theta}{max} x^\top A^\top A x$, i.e., $A^\top \Theta A \preceq \subscr{\theta}{max} A^\top A$. Hence it suffices to show}
$
M - \alpha^2 {\subscr{\theta}{max}} A^\top A\succ 0,
$
which holds if and only if
\begin{align*}
\subscr{\lambda}{min}(M) > \alpha^2 {\subscr{\theta}{max}} \subscr{\lambda}{max}(A^\top A) \iff  \alpha^2 < \frac{\subscr{m}{min}}{{\subscr{\theta}{max}}\amax}.
\end{align*}
The bound $\ds \alpha^2 < \frac{\subscr{m}{min}}{{\subscr{\theta}{max}}\amax}$ follows from the tighter inequality $(2\alpha)^2 \leq \displaystyle \frac{\subscr{m}{min}}{{\subscr{\theta}{max}}\amax}$ which is proved as follows:
\begin{align*}
(2\alpha)^2 = \min\Big\{\frac{\subscr{m}{min}}{\subscr{b}{max}},{\frac{1}{\subscr{\theta}{max}}}\frac{\subscr{b}{min}}{\amax}\Big\}^2&\leq
\min\Big\{\frac{\subscr{m}{min}}{\subscr{b}{max}},{\frac{1}{\subscr{\theta}{max}}}\frac{\subscr{b}{min}}{\amax}\Big\} \cdot
\max\Big\{\frac{\subscr{m}{min}}{\subscr{b}{max}},{\frac{1}{\subscr{\theta}{max}}}\frac{\subscr{b}{min}}{\amax}\Big\}\\
&= \frac{\subscr{m}{min}}{\subscr{b}{max}} \cdot {\frac{1}{\subscr{\theta}{max}}}\frac{\subscr{b}{min}}{\amax} \leq \frac{\subscr{m}{min}}{{\subscr{\theta}{max}}\amax}.
\end{align*}

\textbf{Expansion of the LMI.}\\
Next, we aim to show that $-S^\top P - PS -2 \subscr{c}{eq} P\succeq 0$. {Since $\Theta\Theta^{-1}=I_p$, direct block multiplication gives}
\begin{align}
-  S^\top P&  - PS - 2 \subscr{c}{eq} P
=
\begin{bmatrix}
BM^{-1} & {-A^\top\Theta} \\
A M^{-1} & 0
\end{bmatrix}
\begin{bmatrix}
M & \alpha A^\top \\
\alpha A & {\Theta^{-1}}
\end{bmatrix}
+
\begin{bmatrix}
M & \alpha A^\top \\
\alpha A & {\Theta^{-1}}
\end{bmatrix}
\begin{bmatrix}
M^{-1}B & M^{-1}A^\top \\
{-\Theta A} & 0
\end{bmatrix}
- 2 \subscr{c}{eq}
\begin{bmatrix}
M  & \alpha A^\top \\
\alpha A & {\Theta^{-1}}
\end{bmatrix}
\\
&= 
\begin{bmatrix}
2 B - 2\alpha A^\top {\Theta} A  - 2 \subscr{c}{eq} M 
& \left(\alpha BM^{-1} - 2 \subscr{c}{eq} \alpha  I_n\right)A^\top \\
A \left(\alpha M^{-1}B - 2 \subscr{c}{eq} \alpha I_n\right)
&  2 \alpha A M^{-1}A^\top  - 2 \subscr{c}{eq} {\Theta^{-1}}
\end{bmatrix}.
\end{align}
The (2,2) block satisfies the lower bound 
\begin{align*}
2 \alpha A M^{-1}A^\top  - 2 \subscr{c}{eq} {\Theta^{-1}}
&= 2 \Bigl(\frac{1}{2} \alpha A M^{-1} A^\top - \subscr{c}{eq} {\Theta^{-1}}\Bigr) + \alpha A M^{-1} A^\top \\
&{\succeq 2\Bigl(\frac{1}{2}\alpha \frac{\amin}{\subscr{m}{max}} - \subscr{c}{eq}\subscr{\theta}{max}^{-1}\Bigr)I_p + \alpha AM^{-1}A^\top}\\
& \overset{\eqref{lem:def:c}}{=} \alpha A M^{-1} A^\top\succ 0,
\end{align*}
{where we used $AM^{-1}A^\top \succeq \subscr{m}{max}^{-1}AA^\top \succeq (\amin/\subscr{m}{max})I_p$ and $\Theta^{-1}\preceq \subscr{\theta}{max}^{-1}I_p$.}

Given this lower bound, we can factorize the resulting matrix as follows:
$$
-S^\top P - PS -2 \subscr{c}{eq} P \succeq
\begin{bmatrix} I_n & 0 \\ 0 & A\end{bmatrix}
\underbrace{\begin{bmatrix}
2 B - 2\alpha A^\top {\Theta} A  - 2 \subscr{c}{eq} M 
&\alpha BM^{-1} - 2 \subscr{c}{eq} \alpha  I_n\\
\alpha M^{-1}B - 2 \subscr{c}{eq} \alpha I_n
& \alpha M^{-1}
\end{bmatrix}}_{2n\times 2n}
\begin{bmatrix} I_n & 0 \\ 0 & A^\top \end{bmatrix}.
$$

\textbf{Schur complement verification.}\\
Since $\alpha M^{-1} \succ 0$, it now suffices to show that the Schur complement of the (2,2) block of the $2n\times 2n$ matrix is positive
semidefinite. We compute
\begin{align*}      
& 2 B - 2\alpha A^\top {\Theta} A  - 2 \subscr{c}{eq} M - \alpha \bigl(BM^{-1} - 2 \subscr{c}{eq}  I_n\bigr) M 
\bigl(M^{-1}B - 2 \subscr{c}{eq} I_n\bigr) \succeq 0 \\
& \iff  2 B - 2\alpha A^\top {\Theta} A  - 2 \subscr{c}{eq} M - \alpha \bigl(B - 2 \subscr{c}{eq} M\bigr)\bigl(M^{-1}B - 2 \subscr{c}{eq} I_n\bigr) \succeq 0\\
& \iff  2 B - 2\alpha A^\top {\Theta} A  - 2 \subscr{c}{eq} M - \alpha \bigl(BM^{-1}B - 2 \subscr{c}{eq} B - 2 \subscr{c}{eq} B + 4 \subscr{c}{eq}^2 M\bigr) \succeq 0\\
&\impliedby 2B  - \alpha BM^{-1}B \succeq 2\alpha A^\top {\Theta} A  + 2 \subscr{c}{eq} M
\quad \text{and} \quad 4 \alpha \subscr{c}{eq} B \succeq 4 \alpha \subscr{c}{eq}^2 M.
\end{align*}

\textbf{Prove $\mathbf{2B - \alpha BM^{-1}B \succeq 2\alpha A^\top {\Theta} A  + 2 \subscr{c}{eq} M}$.}\\
To prove
\beq
\label{eq:ineq_1_lem}
2B - \alpha BM^{-1}B \succeq 2\alpha A^\top {\Theta} A  + 2 \subscr{c}{eq} M,
\eeq
we first upper bound its right hand side. {Using $A^\top\Theta A \preceq \subscr{\theta}{max}A^\top A \preceq \subscr{\theta}{max}\amax I_n$ (as in Step 1) and $M\preceq \subscr{m}{max}I_n$,}
\begin{align*}
2\alpha A^\top {\Theta} A  + 2 \subscr{c}{eq} M &\preceq \alpha \left(2{\subscr{\theta}{max}} \amax + 2 \subscr{c}{eq} \subscr{m}{max}\right)I_n \overset{\eqref{lem:def:c}}{\preceq} \alpha\left(2{\subscr{\theta}{max}}\amax + {\subscr{\theta}{min}}\amin\right)I_n \\
& \!\!\!\!\!\!\!\!\! \overset{\alpha\leq\frac{1}{2}\subscr{b}{min}/({\subscr{\theta}{max}}\amax)}{\preceq}
\frac{\subscr{b}{min}}{2{\subscr{\theta}{max}}\amax}\left(2{\subscr{\theta}{max}}\amax + {\subscr{\theta}{min}}\amin\right) I_n = \subscr{b}{min}\Bigl(1+{\frac{\subscr{\theta}{min}\amin}{2\subscr{\theta}{max}\amax}}\Bigr)I_n \preceq \ds \frac{3}{2}\subscr{b}{min} I_n,
\end{align*}
{where the last inequality uses $\subscr{\theta}{min}\leq\subscr{\theta}{max}$ and $\amin\leq\amax$.}
Next, we upper bound the left hand side of~\eqref{eq:ineq_1_lem} exactly as follows
\begin{align*}
2B - \alpha BM^{-1}B &\succeq 2B - \alpha\frac{1}{\subscr{m}{min}} B^2 \succeq 2B - \alpha\frac{\subscr{b}{max}}{\subscr{m}{min}} B \succeq \left(2 - \alpha\frac{\subscr{b}{max}}{\subscr{m}{min}}\right) B \succeq \left(2 - \frac{1}{2}\right) B \succeq \frac{3}{2}\subscr{b}{min} I_n.
\end{align*}

\textbf{Prove $\mathbf{4 \alpha \subscr{c}{eq} B \succeq 4 \alpha \subscr{c}{eq}^2 M}$.}\\
The LMI $4 \alpha \subscr{c}{eq} B \succeq 4 \alpha \subscr{c}{eq}^2 M$ follows from the inequality $\ds \subscr{c}{eq} \leq \frac{\subscr{b}{min}}{\subscr{m}{max}}$. This follows from~\eqref{lem:def:c} and $\alpha \leq \frac{1}{2}\subscr{b}{min}/({\subscr{\theta}{max}}\amax)$:
\[
\subscr{c}{eq} = \frac{1}{2}\alpha{\subscr{\theta}{min}}\frac{\amin}{\subscr{m}{max}} \leq \frac{1}{4}\frac{\subscr{b}{min}}{{\subscr{\theta}{max}}\amax}\cdot{\subscr{\theta}{min}}\frac{\amin}{\subscr{m}{max}} = \frac{1}{4}{\frac{\subscr{\theta}{min}}{\subscr{\theta}{max}}}\frac{\amin\subscr{b}{min}}{\amax\subscr{m}{max}} < \frac{\subscr{b}{min}}{\subscr{m}{max}}.
\]
This concludes the proof.
\end{proof}
\end{arxiv}

\begin{arxiv}%
Next, we derive explicit bounds on the log-norm of the generalized saddle matrices that appear as Jacobians of the SPPI for inequality constraints.
\begin{mylemma}[label=lemma:generalized-saddle_ineq]{Log-norm of generalized saddle matrices}
Given $C\in\R^{m \times n}$, ${B=B^\top\in\R^{n \times n}}$, and
$G=G^\top\in\R^{m\times m}$ with $0\preceq G\preceq I_m$, for
$m\leq n$ and $\alpha,\beta>0$, we consider the \emph{saddle matrix}
$$
\bar{S} =
\begin{bmatrix}
- B - \alpha C^\top G C & - C^\top G \\
\beta\alpha GC & - \beta (I_m - G)
\end{bmatrix} \in \R^{(m+n)\times (m+n)}.
$$
Then, for each $(B,C,G)$ satisfying $\subscr{b}{min} I_n\preceq B \preceq
\subscr{b}{max} I_n$ and $\subscr{\textup{c}}{min} I_m \preceq C C^\top
\preceq \subscr{\textup{c}}{max} I_m$, for $\subscr{b}{min}$,
$\subscr{b}{max}$, $\subscr{\textup{c}}{min}$, $\subscr{\textup{c}}{max}
\in \realpositive$, the following LMI holds:
$$
\bar{S}^\top P + P\bar{S}\preceq -2 \subscr{c}{in} P
\quad\iff\quad
\lognorm{\bar{S}}{P} \leq - \subscr{c}{in},
$$
where
$$
P =
\begin{bmatrix}
\gamma I_n & C^\top \\
C & \dfrac{\gamma}{\alpha\beta} I_m
\end{bmatrix}
\succ 0,
\qquad
\subscr{c}{in} := \frac{3\alpha\beta\subscr{\textup{c}}{min}}{8\gamma},
$$
and $\gamma>0$ is given by
$$
\gamma := \max\left(\beta, \alpha\subscr{\textup{c}}{max},
\frac{2}{3} \subscr{b}{max} + \dfrac{3\alpha\beta\subscr{\textup{c}}{min}}{8 \subscr{b}{min}} +\dfrac{2}{3} \dfrac{x(x+2\subscr{b}{max})}{\subscr{b}{min}}\right)+ \eps, \text{ with } \eps > 0 \text{ and } x := \alpha \subscr{\textup{c}}{max} + \frac{3\beta\subscr{\textup{c}}{min}}{4 \subscr{\textup{c}}{max}} + \beta\sqrt{\frac{\subscr{\textup{c}}{max}}{\subscr{\textup{c}}{min}}}.
$$
\end{mylemma}
\begin{proof}
\textbf{Positivity of P.}\\
We verify $P \succ 0$ via the Schur complement of the $(2,2)$ block. Since $I_m \succ 0$, we need
$
\gamma I_n - \frac{\alpha\beta}{\gamma }C^\top C \succ 0.
$
This LMI holds if and only if $\gamma^2 > \alpha\beta\subscr{\textup{c}}{max}$, which is guaranteed by hypothesis, being $\gamma > \alpha\subscr{\textup{c}}{max}$ and $\gamma > \beta$.

\textbf{Expansion of the LMI.}\\
We aim to show that $-\bar{S}^\top P - P\bar{S} -2 \subscr{c}{in} P\succeq 0$. We compute
\begin{align*}
-\bar{S}^\top P - P\bar{S} -2 \subscr{c}{in} P
&=
\begin{bmatrix}
B + \alpha C^\top G C & -\beta\alpha C^\top G \\
G C & \beta (I_m - G)
\end{bmatrix}
\begin{bmatrix}
\gamma I_n & C^\top \\
C & \dfrac{\gamma}{\alpha\beta} I_m
\end{bmatrix}
+
\begin{bmatrix}
\gamma I_n & C^\top \\
C & \dfrac{\gamma}{\alpha\beta} I_m
\end{bmatrix}
\begin{bmatrix}
B + \alpha C^\top G C & C^\top G \\
-\beta\alpha G C & \beta (I_m - G)
\end{bmatrix}\\
&\quad - 2 \subscr{c}{in}
\begin{bmatrix}
\gamma I_n & C^\top \\
C & \dfrac{\gamma}{\alpha\beta} I_m
\end{bmatrix}
\\&=
\begin{bmatrix}
\gamma(B+\alpha C^\top GC) - \beta\alpha C^\top GC
&
(B+\alpha C^\top GC)C^\top - \gamma C^\top G \\
\gamma GC + \beta C - \beta GC
&
GCC^\top + \dfrac{\gamma}{\alpha}(I_m-G)
\end{bmatrix}
\\&\quad+
\begin{bmatrix}
\gamma(B+\alpha C^\top GC) - \beta\alpha C^\top GC
&
\gamma C^\top G + \beta C^\top - \beta C^\top G \\
CB + \alpha CC^\top GC - \gamma GC
&
CC^\top G + \dfrac{\gamma}{\alpha}(I_m-G)
\end{bmatrix}
-
\begin{bmatrix}
2\subscr{c}{in}\gamma I_n & 2\subscr{c}{in}C^\top \\
2\subscr{c}{in}C & \dfrac{2\subscr{c}{in}\gamma}{\alpha\beta}I_m
\end{bmatrix}
\\&=
\begin{bmatrix}
2\gamma B + 2\alpha(\gamma-\beta)C^\top G C - 2\subscr{c}{in}\gamma I_n
&
\bigl(B+\beta I_n - 2\subscr{c}{in}I_n\bigr)C^\top +
C^\top G\bigl(\alpha CC^\top - \beta I_m\bigr) \\
C\bigl(B+\beta I_n - 2\subscr{c}{in}I_n\bigr) +
\bigl(\alpha CC^\top - \beta I_m\bigr)GC
&
GCC^\top + CC^\top G + \frac{2\gamma}{\alpha}(I_m-G) -
\dfrac{2\subscr{c}{in}\gamma}{\alpha\beta}I_m
\end{bmatrix} \succeq 0.
\end{align*}
First, we find a lower bound for $Q_{22}$. Let $G=U\diag{d}U^\top$ be an eigendecomposition, with $U \in \R^{m \times m}$ orthogonal and $d\in[0,1]^m$. Substituting this decomposition into $Q_{22}$ and multiplying on the left and on the right by $U$ and $U^\top$, respectively, we get
\beq
\label{eq:decomposition}
U^\top Q_{22}U = \diag{d}U^\top CC^\top U + U^\top CC^\top U\diag{d} + \frac{2\gamma}{\alpha}(I_m- \diag{d}) - \dfrac{2\subscr{c}{in}\gamma}{\alpha\beta}I_m
\eeq 
Note that the eigenvalues of $U^\top CC^\top U$ and $CC^\top$ coincide. Then, applying Lemma~\cite[Lemma 6]{GQ-NL:19} to~\eqref{eq:decomposition} with $X:=U^\top CC^\top U$, and multiplying this on the left and on the right by $U$ and $U^\top$, for any $\gamma \geq \alpha \subscr{\textup{c}}{max}$ we get
$$
(C C^\top G + GC C^\top) + 2 \frac{\beta}{\alpha} (I_m - G) - \dfrac{2\subscr{c}{in}\gamma}{\alpha\beta}I_m \succeq \frac{3}{2} C C^\top - \dfrac{2\subscr{c}{in}\gamma}{\alpha\beta}I_m.
$$
Then, the (2,2) block satisfies the lower bound 
\begin{align*}
(C C^\top G + GC C^\top) + 2 \frac{\beta}{\alpha} (I_m - G) - \dfrac{2\subscr{c}{in}\gamma}{\alpha\beta} I_m  &\succeq \frac{3}{2} C C^\top - \dfrac{2\subscr{c}{in}\gamma}{\alpha\beta}I_m\\
&= 2 \Big(\frac{3}{8} C C^\top -\dfrac{\subscr{c}{in}\gamma}{\alpha\beta} I_m\Big) + \frac{3}{4} C C^\top \\
& \succeq 2\Big(\frac{3}{8} \subscr{\textup{c}}{min} - \dfrac{\subscr{c}{in}\gamma}{\alpha\beta} \Big)I_m + \frac{3}{4}C C^\top\\
&\overset{\subscr{c}{in}:=  \frac{3\alpha\beta\subscr{\textup{c}}{min}}{8\gamma}}{=} \frac{3}{4} C C^\top\succ 0.
\end{align*}
Given this lower bound, we compute
$$
-\bar{S}^\top P - P\bar{S} -2 \subscr{c}{in} P \succeq
\begin{bmatrix}
2\gamma B + 2\alpha(\gamma-\beta)C^\top G C - 2\subscr{c}{in}\gamma I_n
&
\bigl(B+\beta I_n - 2\subscr{c}{in}I_n\bigr)C^\top +
C^\top G\bigl(\alpha CC^\top - \beta I_m\bigr) \\
C\bigl(B+\beta I_n - 2\subscr{c}{in}I_n\bigr) +
\bigl(\alpha CC^\top - \beta I_m\bigr)GC
&
\frac{3}{4} C C^\top
\end{bmatrix}.
$$
It now suffices to show that the Schur complement of the (2,2) block of the above matrix is positive semidefinite. To this end, we first bound $Q_{12}Q_{22}^{-1}Q_{12}^\top$.

\textbf{Schur complement verification.}\\
First, we bound the matrix $Q_{12}Q_{22}^{-1}Q_{12}^\top$. Let $Q_{12} = \bigl(B+\beta I_n - 2\subscr{c}{in}I_n\bigr)C^\top + C^\top G\bigl(\alpha CC^\top - \beta I_m\bigr) := \bigl(B + Q_3\bigr) C^\top + \beta C^\top(I_m-G)$, with $Q_3 = \alpha C^\top GC - 2\subscr{c}{in}I_n$. We have
\begin{align*}
\frac{3}{4}Q_{12}&Q_{22}^{-1}Q_{12}^\top \preceq Q_{12}(CC^\top)^{-1}Q_{12}^\top \\
& = \left(\bigl(B + Q_3\bigr) C^\top + \beta C^\top(I_m-G)\right)(CC^\top)^{-1}\left(\bigl(B + Q_3\bigr) C^\top + \beta C^\top(I_m-G)\right)^\top\\
& = \left(\bigl(B + Q_3\bigr) C^\top + \beta C^\top(I_m-G)\right)(CC^\top)^{-1}\left(C \bigl(B + Q_3^\top\bigr) + \beta (I_m-G)C\right)\\
& = \underbrace{\bigl(B + Q_3\bigr) C^\top (CC^\top)^{-1}C \bigl(B + Q_3^\top\bigr)}_{\text{(i)}} + \underbrace{\beta \bigl(B + Q_3\bigr) C^\top (CC^\top)^{-1} (I_m-G)C + 
\beta C^\top(I_m-G)(CC^\top)^{-1}C \bigl(B + Q_3^\top\bigr)}_{\text{(ii)}} \\
& \quad + \underbrace{\beta^2 C^\top(I_m-G)(CC^\top)^{-1}(I_m-G)C}_{\text{(iii)}}.
\end{align*}
We bound each term separately. To do so, we use the bounds $\norm{I_m-G} \leq 1$ and $\norm{C^\top} = \norm{C} = \sqrt{\subscr{\textup{c}}{max}}$. Moreover, we have 
$$
\norm{Q_3} = \norm{\alpha C^\top GC - 2\subscr{c}{in}I_n} \leq \alpha \norm{C^\top GC} + 2\subscr{c}{in} =  \alpha \norm{C^\top GC} + 2 \frac{3\alpha\beta\subscr{\textup{c}}{min}}{8\gamma} \leq 
\alpha \subscr{\textup{c}}{max} + \frac{3\beta\subscr{\textup{c}}{min}}{4 \subscr{\textup{c}}{max}},
$$
where we used the inequality $\gamma \geq \alpha \subscr{\textup{c}}{max}$. We compute

\emph{Term (i):} Since $\subscr{\textup{c}}{min}>0$, $C$ has full row rank, so $C^\top(CC^\top)^{-1}C$ is the orthogonal projector onto $\mathrm{row}(C)$, giving $C^\top(CC^\top)^{-1}C \preceq I_n$. Then
\beq
\label{eq:bound_i}
\begin{aligned}
\bigl(B + Q_3\bigr) C^\top (CC^\top)^{-1}C \bigl(B + Q_3\bigr)^\top &\preceq \bigl(B + Q_3\bigr)^2 = B^2+BQ_3+Q_3B+Q_3^2 \preceq \subscr{b}{max}B +\bigl(2\subscr{b}{max}\norm{Q_3}+\norm{Q_3}^2\bigr)I_n\\
&\preceq \subscr{b}{max}B + \left(2\subscr{b}{max}(\alpha \subscr{\textup{c}}{max} + \frac{3\beta\subscr{\textup{c}}{min}}{4 \subscr{\textup{c}}{max}}) + \Bigl(\alpha \subscr{\textup{c}}{max} + \frac{3\beta\subscr{\textup{c}}{min}}{4 \subscr{\textup{c}}{max}}\Bigr)^2\right)I_n,
\end{aligned}
\eeq
where we used the notation $X^2 =: X^\top X$ as shorthand for the symmetric product.

\emph{Term (ii):} Noting that $C^\top(CC^\top)^{-1}$ is the Moore--Penrose pseudoinverse of $C$, with $\norm{C^\top(CC^\top)^{-1}} = 1/\sqrt{\subscr{\textup{c}}{min}}$, we compute
\beq
\label{eq:bound_iii}
\begin{aligned}
\norm{\beta \bigl(B + Q_3\bigr) C^\top (CC^\top)^{-1} (I_m-G)C + \beta C^\top(I_m-G)(CC^\top)^{-1}C \bigl(B + Q_3^\top\bigr)} &\leq 2 \beta \norm{B + Q_3} \norm{C^\top (CC^\top)^{-1}} \norm{(I_m-G)} \norm{C^\top} \\
&\leq 2 \beta \norm{B + Q_3} \dfrac{\sqrt{\subscr{\textup{c}}{max}}}{\sqrt{\subscr{\textup{c}}{min}}}\\
&\leq 2 \beta \Bigl(\subscr{b}{max} + \alpha \subscr{\textup{c}}{max} + \frac{3\beta\subscr{\textup{c}}{min}}{4 \subscr{\textup{c}}{max}}\Bigr) \dfrac{\sqrt{\subscr{\textup{c}}{max}}}{\sqrt{\subscr{\textup{c}}{min}}}.
\end{aligned}
\eeq

\emph{Term (iii):} Since $\|(CC^\top)^{-1}\| = 1/\subscr{\textup{c}}{min}$, we have bound
\beq
\label{eq:bound_iii}
\norm{\beta^2 C^\top(I_m-G)(CC^\top)^{-1}(I_m-G)C} \leq \beta^2 \norm{C^\top}^2 \norm{I_m-G}^2 \norm{(CC^\top)^{-1}}\leq \dfrac{\beta^2 \subscr{\textup{c}}{max}}{\subscr{\textup{c}}{min}}.
\eeq
Combining bounds~\eqref{eq:bound_i}--\eqref{eq:bound_iii} we have
\begin{align*}
\frac{3}{4}Q_{12}&Q_{22}^{-1}Q_{12}^\top \preceq
\subscr{b}{max}B + \left(\Bigl(2\subscr{b}{max} + \alpha \subscr{\textup{c}}{max} + \frac{3\beta\subscr{\textup{c}}{min}}{4 \subscr{\textup{c}}{max}}\Bigr)\Bigl(\alpha \subscr{\textup{c}}{max} + \frac{3\beta\subscr{\textup{c}}{min}}{4 \subscr{\textup{c}}{max}}\Bigr)
+ 
2 \beta \Bigl(\subscr{b}{max} + \alpha \subscr{\textup{c}}{max} + \frac{3\beta\subscr{\textup{c}}{min}}{4 \subscr{\textup{c}}{max}}\Bigr) \dfrac{\sqrt{\subscr{\textup{c}}{max}}}{\sqrt{\subscr{\textup{c}}{min}}}
+
\dfrac{\beta^2 \subscr{\textup{c}}{max}}{\subscr{\textup{c}}{min}}\right) I_n.
\end{align*}

Next, we note that $\gamma \geq \beta$ and $C^\top GC \succeq 0$, then $\alpha(\gamma-\beta)C^\top GC \succeq 0$ which in turn implies
$$
Q_{11} = 2\gamma B + 2\alpha(\gamma-\beta)C^\top G C - 2\subscr{c}{in}\gamma I_n \succeq 2\gamma B - 2 \gamma \subscr{c}{in}I_n.
$$
Then
\begin{align*}
Q_{11} - Q_{12}&Q_{22}^{-1}Q_{12}^\top \succeq (2\gamma - \frac{4}{3}\subscr{b}{max})B - 2 \gamma \subscr{c}{in}I_n \\
&\quad - \frac{4}{3} \left(\Bigl(2\subscr{b}{max} + \alpha \subscr{\textup{c}}{max} + \frac{3\beta\subscr{\textup{c}}{min}}{4 \subscr{\textup{c}}{max}}\Bigr)\Bigl(\alpha \subscr{\textup{c}}{max} + \frac{3\beta\subscr{\textup{c}}{min}}{4 \subscr{\textup{c}}{max}}\Bigr)
+ 
2 \beta \Bigl(\subscr{b}{max} + \alpha \subscr{\textup{c}}{max} + \frac{3\beta\subscr{\textup{c}}{min}}{4 \subscr{\textup{c}}{max}}\Bigr) \dfrac{\sqrt{\subscr{\textup{c}}{max}}}{\sqrt{\subscr{\textup{c}}{min}}}
+
\dfrac{\beta^2 \subscr{\textup{c}}{max}}{\subscr{\textup{c}}{min}}\right) I_n \\
&\overset{2\gamma\subscr{c}{in}:=  \frac{3\alpha\beta\subscr{\textup{c}}{min}}{4}}{=}
\bigl(2\gamma - \frac{4}{3}\subscr{b}{max}\bigr)B - \frac{3\alpha\beta\subscr{\textup{c}}{min}}{4}I_n \\
&\quad - \frac{4}{3} \left(\Bigl(2\subscr{b}{max} + \alpha \subscr{\textup{c}}{max} + \frac{3\beta\subscr{\textup{c}}{min}}{4 \subscr{\textup{c}}{max}}\Bigr)\Bigl(\alpha \subscr{\textup{c}}{max} + \frac{3\beta\subscr{\textup{c}}{min}}{4 \subscr{\textup{c}}{max}}\Bigr)
+ 
2 \beta \Bigl(\subscr{b}{max} + \alpha \subscr{\textup{c}}{max} + \frac{3\beta\subscr{\textup{c}}{min}}{4 \subscr{\textup{c}}{max}}\Bigr) \dfrac{\sqrt{\subscr{\textup{c}}{max}}}{\sqrt{\subscr{\textup{c}}{min}}}
+
\dfrac{\beta^2 \subscr{\textup{c}}{max}}{\subscr{\textup{c}}{min}}\right) I_n \\ 
& \succeq \Bigl(\bigl(2\gamma - \frac{4}{3}\subscr{b}{max}\bigr)\subscr{b}{min} - \frac{3\alpha\beta\subscr{\textup{c}}{min}}{4}\Bigr)I_n \\
&\quad - \frac{4}{3} \left(\Bigl(2\subscr{b}{max} + \alpha \subscr{\textup{c}}{max} + \frac{3\beta\subscr{\textup{c}}{min}}{4 \subscr{\textup{c}}{max}}\Bigr)\Bigl(\alpha \subscr{\textup{c}}{max} + \frac{3\beta\subscr{\textup{c}}{min}}{4 \subscr{\textup{c}}{max}}\Bigr)
+ 
2 \beta \Bigl(\subscr{b}{max} + \alpha \subscr{\textup{c}}{max} + \frac{3\beta\subscr{\textup{c}}{min}}{4 \subscr{\textup{c}}{max}}\Bigr) \dfrac{\sqrt{\subscr{\textup{c}}{max}}}{\sqrt{\subscr{\textup{c}}{min}}}
+
\dfrac{\beta^2 \subscr{\textup{c}}{max}}{\subscr{\textup{c}}{min}}\right) I_n \\
& \succeq \Bigl(\bigl(2\gamma - \frac{4}{3}\subscr{b}{max}\bigr)\subscr{b}{min}
- \frac{3\alpha\beta\subscr{\textup{c}}{min}}{4}\Bigr)I_n \\
&\quad - \frac{4}{3} \Bigl(\alpha \subscr{\textup{c}}{max} + \frac{3\beta\subscr{\textup{c}}{min}}{4 \subscr{\textup{c}}{max}} + \beta\sqrt{\frac{\subscr{\textup{c}}{max}}{\subscr{\textup{c}}{min}}}\Bigr)\Bigl(\alpha \subscr{\textup{c}}{max} + \frac{3\beta\subscr{\textup{c}}{min}}{4 \subscr{\textup{c}}{max}} + \beta\sqrt{\frac{\subscr{\textup{c}}{max}}{\subscr{\textup{c}}{min}}} + 2\subscr{b}{max}\Bigr) I_n
,
\end{align*}
where in the last inequality we used the fact that $\gamma \geq \dfrac{2\subscr{b}{max}}{3}$.
Let $x := \alpha \subscr{\textup{c}}{max} + \frac{3\beta\subscr{\textup{c}}{min}}{4 \subscr{\textup{c}}{max}} + \beta\sqrt{\frac{\subscr{\textup{c}}{max}}{\subscr{\textup{c}}{min}}}$. 
Then, the condition $\gamma \geq \dfrac{2}{3} \subscr{b}{max} + \dfrac{3\alpha\beta\subscr{\textup{c}}{min}}{8 \subscr{b}{min}} +\dfrac{2}{3} \dfrac{x(x+2\subscr{b}{max})}{\subscr{b}{min}}$ implies $Q_{11} - Q_{12}Q_{22}^{-1}Q_{12}^\top \succeq 0$ and thus the thesis. This concludes the proof.
\end{proof}

\section{Review of Dynamics Solving the Affine Inequality-Constrained Problem}
\label{apx:review_lin}
For completeness we recall the two primal--dual dynamical systems proposed in the literature to solve the affine inequality-constrained problem $\min\setdef{f(x)}{Cx \leq d}$ that are identified with the PI-PSPF in Section~\ref{sec:convergence_ineq}. We also state the convergence results of~\cite{GQ-NL:19} and~\cite[Theorem~11]{AD-VC-AG-GR-FB:23f} that are invoked therein.

\subsection*{Augmented Primal--Dual Gradient Dynamics~\cite{GQ-NL:19}}
The Augmented Primal--Dual Gradient Dynamics (Aug-PDGD) of~\cite{GQ-NL:19} is
\beq
\label{eq:aug_pdgd_qu_li}
\begin{cases}
\dot{x} = -\nabla f(x) - C^\top \relu\!\big(\nu(Cx-d) + \mu\big), \\
\dot{\mu} = \dfrac{\eta}{\nu}\big(- \mu + \relu\big(\nu(Cx-d) + \mu\big)\big).
\end{cases}
\eeq
This coincide with SPPI~\eqref{eq:pi_spf_ineq_affine} when $\nu =  \subscr{k}{p}^{\textup{in}}$ and $\eta =  \subscr{k}{i}^{\textup{in}} \subscr{k}{p}^{\textup{in}}$. Theorem 2 in~\cite{GQ-NL:19} establishes global exponential stability of~\eqref{eq:aug_pdgd_qu_li} under Assumption~\ref{ass:1}, for any $\nu, \eta > 0$.
The proof proceeds via a quadratic Lyapunov function with metric $P\succ 0$ defined by, in terms of the PI gains,
\beq
\label{eq:lyap_metric_QL}
\subscr{P}{QL} = 
\begin{bmatrix}
\subscr{k}{p}^{\textup{in}}\subscr{k}{i}^{\textup{in}}\kappa I_n & \subscr{k}{p}^{\textup{in}}\subscr{k}{i}^{\textup{in}} C^\top \\
\subscr{k}{p}^{\textup{in}}\subscr{k}{i}^{\textup{in}} C              & \kappa I_m
\end{bmatrix} \in \R^{(n+m)\times(n+m)},
\eeq
where $\displaystyle \kappa = 20 L \frac{\amax}{\amin} \left[\max \left(\frac{\subscr{k}{p}^{\textup{in}} \amax}{\rho},\frac{L}{\rho}\right)\right]^{2} \left[\max \left(\frac{\subscr{k}{i}^{\textup{in}}}{L},\frac{L}{\rho}\right)\right]^{2}$.
Expressing the exponential rate in terms of the PI gains yields
\beq
\label{eq:exp_rate_ql}
\subscr{c}{QL} = \frac{1}{2}\frac{\subscr{k}{i}^{\textup{in}}\subscr{k}{p}^{\textup{in}}\amin^2}
{40L \amax
\left[\max \left(\dfrac{\subscr{k}{p}^{\textup{in}}\amax}{\rho},\dfrac{L}{\rho}\right)\right]^2
\left[\max \left(\dfrac{\subscr{k}{i}^{\textup{in}}}{L}, \dfrac{L}{\rho}\right)\right]^2}.
\eeq

\subsection*{Proximal Augmented Lagrangian Primal--Dual Dynamics~\cite{AD-VC-AG-GR-FB:23f}}
The proximal augmented Lagrangian primal--dual dynamics of~\cite{AD-VC-AG-GR-FB:23f} is
\beq
\label{eq:prox-aug-inequality}
\begin{cases}
\dot{x} = -\nabla f(x) - \dfrac{1}{\gamma} C^\top \relu(Cx + \gamma\mu - d), \\
\dot{\mu} = -\gamma\mu + \relu(Cx + \gamma\mu - d).
\end{cases}
\eeq
This coincide with SPPI~\eqref{eq:pi_spf_ineq_affine} when $\subscr{k}{p}^{\textup{in}} = \dfrac{1}{\gamma}$, $\subscr{k}{i}^{\textup{in}} = \gamma$. For $\varepsilon \in {]0,1/\sqrt{\amax}[}$, consider the nonlinear program
\begin{subequations}\label{eq:nonlinear-prog}
\begin{alignat}{2}
&\!\max_{c\geq 0,\alpha \geq 0,\varkappa \geq 0}        &\qquad& c\label{eq:optProb}\\
&\text{s.t.} &      & \alpha \leq \min\Big\{\frac{1}{\sqrt{\amax}} - \varepsilon, \frac{\gamma}{\amax}\Big\},\label{eq:alpha-constraint}\\
&                  &      & \varkappa \geq \frac{2}{3},\label{eq:constraint2} \\
&    &    & c \leq \Big(\frac{3}{4}-\frac{1}{2\varkappa}\Big)\alpha \amin, \label{eq:constraint3}\\
& & & h(c,\alpha,\varkappa) \geq 0, \label{eq:nonlinear-constraint}
\end{alignat}
\end{subequations}
with $\map{h}{\realnonnegative \times \realnonnegative \times \realnonnegative}{\real}$ given by
\begin{align*}
h(c,\alpha,\varkappa) =& \ 2\rho - \relu\Big(2\alpha-\frac{2}{\gamma}\Big)\amax - 2c - \alpha\varkappa\frac{\amax}{\amin} \Big(\gamma^2 \frac{\amax}{\amin} + (L + \frac{\amax}{\gamma} + 2c)^2 + 2\gamma\frac{\amax}{\amin} (L + \frac{\amax}{\gamma} + 2c)\Big). 
\end{align*}
The following result holds~\cite[Theorem~11]{AD-VC-AG-GR-FB:23f}.
\bt[Contractivity of~\eqref{eq:prox-aug-inequality}]
\label{thm:contractivity-pal}
Suppose Assumption~\ref{ass:1} holds and let $\gamma > 0$ be arbitrary. Then the primal--dual dynamics~\eqref{eq:prox-aug-inequality} are strongly infinitesimally contracting with respect to $\|\cdot\|_{P}$ with rate $c^\star > 0$ where
\beq
\label{eq:P-def-pal}
P = \begin{bmatrix}
I_n & \alpha^\star C^\top \\
\alpha^\star C & I_m
\end{bmatrix},
\eeq
and $\alpha^\star > 0, c^\star > 0$ are the optimal arguments of problem~\eqref{eq:nonlinear-prog}.
\et

\subsection{Numerical Comparison of Contraction Rates Estimates}
\label{apx:rate_comparison_numerical}
We empirically compare our estimate of the contraction rate of SPPI~\eqref{eq:pi_spf_ineq_affine} against the rates of~\cite[Theorem~2]{GQ-NL:19} (Aug-PDGD) and~\cite[Theorem~11]{AD-VC-AG-GR-FB:23f} (Prox-PDGD), across a range of gains and problem conditioning.

Aug-PDGD and SPPI are each parametrized by two independently tunable gains, whereas Prox-PDGD is parametrized by a single scalar. Comparing the three rates on a common axis therefore requires an explicit identification of parameters. We use two such identifications, one for each figure below. In both figures, $\subscr{c}{in}$ is computed from~\eqref{eq:P_generalized_saddle}, the Aug-PDGD rate $\subscr{c}{QL}$ is computed from~\eqref{eq:exp_rate_ql} and the Prox-PDGD rate is the numerically-optimized $c^\star$ solving the nonlinear program~\eqref{eq:nonlinear-prog}, approximated via bisection over $c$ and a grid search over the auxiliary variable $\varkappa$.

Figure~\ref{fig:rate_comparison_1} plots all three rates as a function of a single shared parameter $\gamma > 0$, identified with all three dynamics as described above, with $\gamma$ playing the role of Prox-PDGD's own free parameter. We consider four conditioning regimes: the example of~\cite[Fig.~7]{AD-VC-AG-GR-FB:23f} ($\rho=1.03$, $L=27.81$, $\amin=0.1$, $\amax=1$, top left panel), a well-conditioned regime ($\rho=1$, $L=5$, $\amin=0.5$, $\amax=1$, top right panel), an ill-conditioned regime ($\rho=0.5$, $L=50$, $\amin=0.05$, $\amax=1$, bottom left panel), and a regime with a near-isotropic $C$ ($\rho=1$, $L=20$, $\amin=0.8$, $\amax=1$, bottom right panel). Across all four regimes, $\subscr{c}{in}$ exceeds the Aug-PDGD rate by four to seven orders of magnitude. In comparison against Prox-PDGD, $\subscr{c}{in}$ exceeds $c^\star$ throughout the original, milder-$L/\rho$, and harsher-$L/\rho$ regimes (top left, top right, and bottom left panels), but the gap narrows substantially, with the two curves becoming comparable for larger $\gamma$, in near-isotropic matrix $C$ (bottom right panel, $\amax/\amin=100$).

\begin{figure}[!h]
\centering
\includegraphics[width=.75\linewidth]{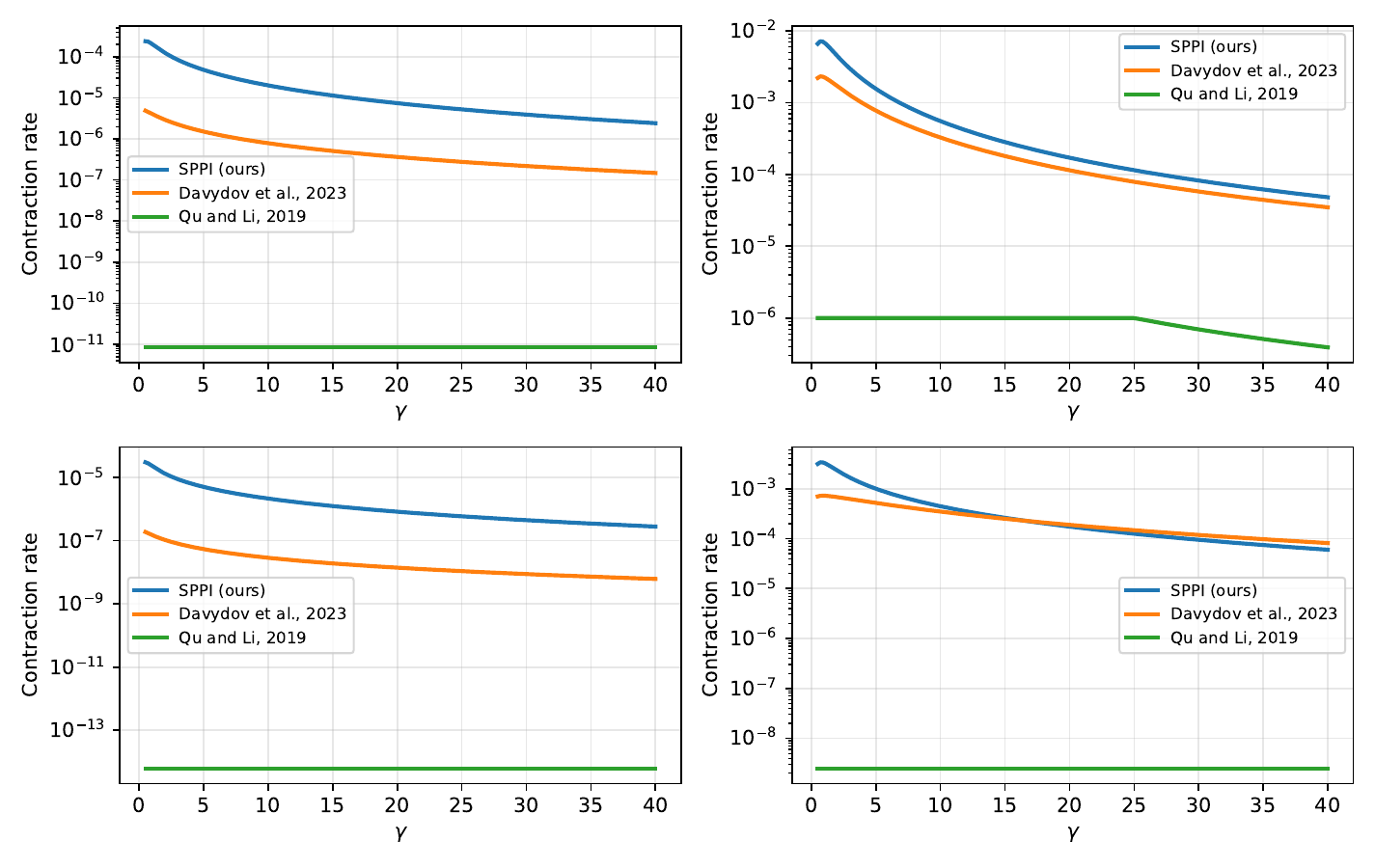}
\caption{Certified rates $\subscr{c}{in}$~\eqref{eq:P_generalized_saddle}, $c^\star$~\cite{AD-VC-AG-GR-FB:23f}, and $\subscr{c}{QL}$~\cite{GQ-NL:19} for the three dynamics under the shared identification over $\gamma$, plotted against $\gamma$, for four choices of $(\rho, L, \amin, \amax)$: the parameters of~\cite[Figure~7]{AD-VC-AG-GR-FB:23f} ($\rho=1.03$, $L=27.81$, $\amin=0.1$, $\amax=1$), a well-conditioned regime ($\rho=1$, $L=5$, $\amin=0.5$, $\amax=1$), an ill-conditioned regime ($\rho=0.5$, $L=50$, $\amin=0.05$, $\amax=1$), and a regime with a near-isotropic $C$, i.e., small $\amax/\amin$ ($\rho=1$, $L=20$, $\amin=0.8$, $\amax=1$).}
\label{fig:rate_comparison_1}
\end{figure}

Figure~\ref{fig:rate_comparison_2} compares $\subscr{c}{in}$, as a function of $\subscr{k}{p}^{\textup{in}}$, against $\subscr{c}{QL}$ using the direct identification stated in the definition of Aug-PDGD, $\nu = \subscr{k}{p}^{\textup{in}}$ and $\eta = \subscr{k}{i}^{\textup{in}} \subscr{k}{p}^{\textup{in}}$. Across the entire tested range of $\subscr{k}{p}^{\textup{in}} \in [10^{-2}, 10^2]$ and for every value of $\subscr{k}{i}^{\textup{in}} \in \{0.5, 2, 10, 50\}$ shown, $\subscr{c}{in}$ exceeds the Aug-PDGD rate by three or more orders of magnitude, including at each curve's respective peak.

\begin{figure}[!h]
\centering
\includegraphics[width=.75\linewidth]{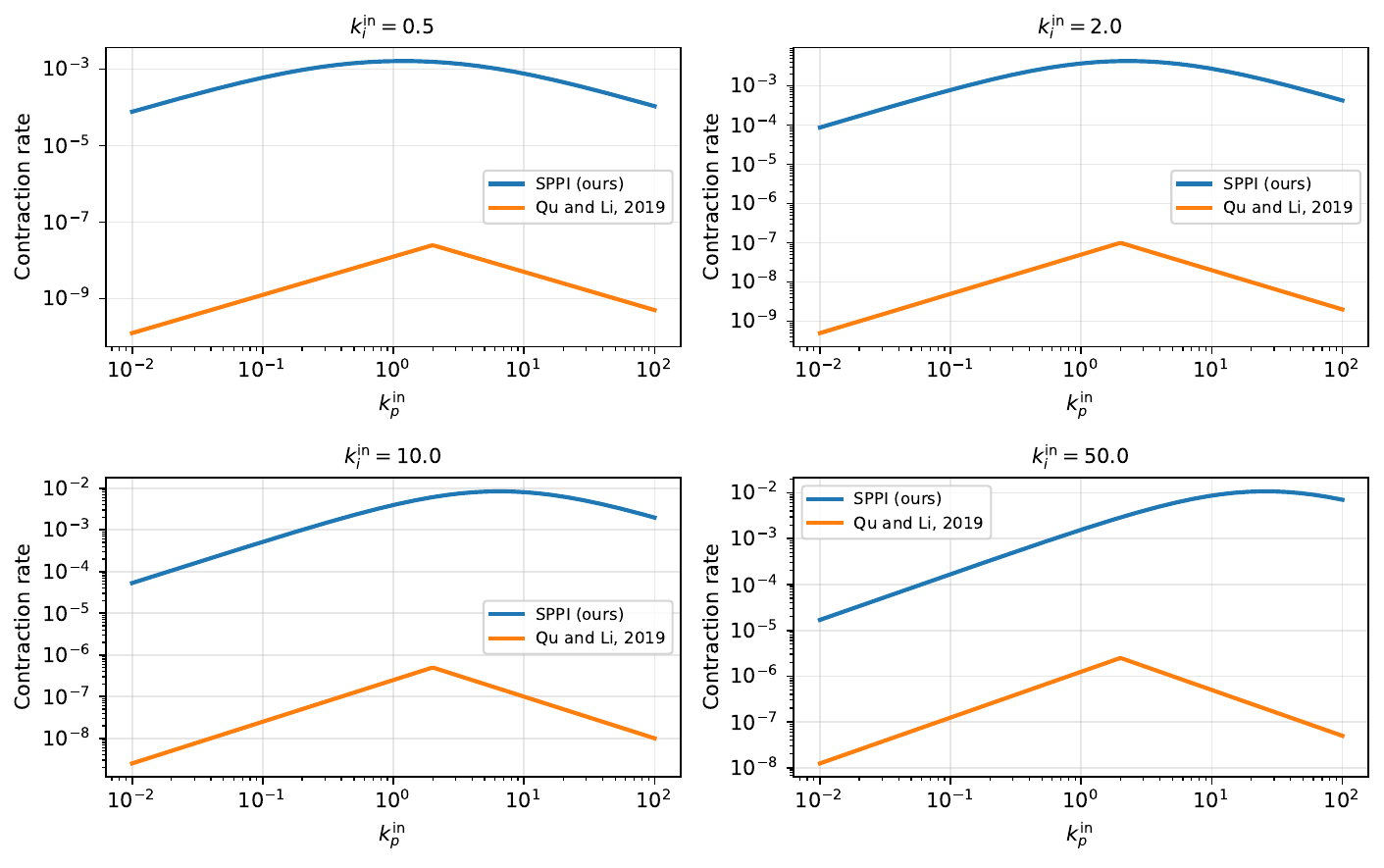}
\caption{Certified rates $\subscr{c}{in}$~\eqref{eq:P_generalized_saddle} and $\subscr{c}{QL}$~\cite{GQ-NL:19}, as a function of $\subscr{k}{p}^{\textup{in}} \in [10^{-2}, 10^2]$, for $\rho=1$, $L=10$, $\amin=1$, $\amax=5$ and four fixed values of $\subscr{k}{i}^{\textup{in}} \in \{0.5, 2, 10, 50\}$, using the direct identification $\nu=\subscr{k}{p}^{\textup{in}}$, $\eta=\subscr{k}{i}^{\textup{in}}\subscr{k}{p}^{\textup{in}}$.}
\label{fig:rate_comparison_2}
\end{figure}
\end{arxiv}
\end{document}